\documentclass[11pt,reqno]{amsart}
\usepackage[T1]{fontenc}
\usepackage{lmodern}
\usepackage{amsmath,amssymb,mathtools,mathrsfs}
\usepackage[a4paper,left=30mm,right=30mm,top=28mm,bottom=29mm]{geometry}
\usepackage{microtype}
\usepackage{enumitem}
\usepackage{needspace}
\usepackage[hidelinks]{hyperref}
\hypersetup{pdftitle={Majorizing measures for canonical processes with log-concave tails},pdfauthor={}}
\allowdisplaybreaks[1]
\numberwithin{equation}{section}
\newtheorem{theorem}{Theorem}[section]
\newtheorem{lemma}[theorem]{Lemma}
\newtheorem{proposition}[theorem]{Proposition}
\newtheorem{corollary}[theorem]{Corollary}
\theoremstyle{definition}

\theoremstyle{remark}

\newcommand{\E}{\textsf E}
\newcommand{\Prob}{\textsf P}
\newcommand{\R}{\mathbb R}

\newcommand{\Z}{\mathbb Z}
\newcommand{\cP}{\mathcal P}
\newcommand{\cR}{\mathcal R}
\newcommand{\cM}{\mathcal M}
\newcommand{\RD}{\operatorname{RD}}
\newcommand{\supp}{\operatorname{supp}}
\newcommand{\Law}{\operatorname{Law}}
\newcommand{\clip}{\operatorname{clip}}

\newcommand{\diam}{\operatorname{diam}}
\newcommand{\cx}{\preceq_{\rm cx}}
\newcommand{\st}{\preceq_{\rm st}}
\newcommand{\ind}{\mathbf 1}
\newcommand{\ip}[2]{\langle #1,#2\rangle}
\newcommand{\norm}[1]{\lVert #1\rVert}
\newcommand{\abs}[1]{\lvert #1\rvert}
\newcommand{\eps}{\varepsilon}
\newcommand{\push}{_{\#}}
\newcommand{\KL}{\operatorname{KL}}
\newcommand{\Ent}{\operatorname{H}}
\newcommand{\asympc}{\asymp}

\title[Majorizing measures for canonical processes]{Majorizing measures for canonical processes\protect\\ with log-concave tails}
\author{Xuanang Hu}
\address{Shandong University,  Jinan, China.}
\email{xuananghu7@gmail.com}

\author{Hanchao Wang}
\address{ Shandong University, Jinan,  China.}
\email{wanghanchao@sdu.edu.cn}

\author{Xinglong Wu}
\address{ Shandong University, Jinan,  China.}
\email{wuxinglong3@gmail.com}
\date{}

\begin{document}
\begin{abstract} Let \(Y_1,\ldots,Y_n\) be independent symmetric random variables with log-concave tails. We give a dimension-free characterization of the expected supremum of the canonical process \[ X_x=\sum_{i=1}^n x_iY_i \] without any \(\Delta_2\) or regular-growth assumption on the coordinate tails. The characterization is governed by scale-dependent intrinsic costs determined directly by the convex tail potentials. It has an information-theoretic rate--distortion formulation as well as two equivalent majorizing-measure formulations on the original index set, one multiscale and one based on a single probability measure. More strongly, the comparison holds for every prescribed law of the index. In particular, the Gaussian majorizing-measure theorem of Talagrand and the Bernoulli theorem of Bednorz--Lata{\l}a are recovered from the same formula, through quadratic and truncated quadratic costs, respectively. 

The new geometry is necessary: once the \(\Delta_2\) growth condition
is removed, the usual chaining functional based on increment moments
can exceed the expected supremum by an unbounded factor. The proof is based on a uniform convex-order truncation principle for symmetric exponentials and a dyadic decomposition of the tail slopes that transfers simultaneously the random process and its intrinsic cost, with universal constants. \end{abstract}

\maketitle

\section{Introduction}

Let \(Y_1,\ldots,Y_n\) be independent symmetric random variables and, for \(T\subset\mathbb R^n\), consider the canonical process
\[
X_x=\ip{x}{Y}=\sum_{i=1}^n x_iY_i,\qquad x\in T.
\]
As usual, to avoid measurability issues we interpret the expected supremum over an arbitrary index set by
\[
S_Y(T):=\sup\left\{
\E\max_{x\in F}X_x:\ F\subset T,\ F\ \text{finite}
\right\}.
\]
Thus the basic problem reduces to finite index sets, and throughout we formulate our results for finite nonempty \(T\). Our aim is to give a deterministic characterization of \(S_Y(T)\) in terms of the geometry of \(T\) and the laws of the coordinates, up to universal multiplicative constants.

For Gaussian coordinates, Fernique's majorizing measure bound and
Talagrand's converse answer this question in terms of Euclidean
distances \cite{Fernique75,Talagrand}. For independent signs, the
Bernoulli theorem of Bednorz and Lata\l a gives a fundamentally
different characterization \cite[Theorem~1.1]{BL}. It decomposes
the index set as \(T\subset T_1+T_2\), with one part controlled by
its \(\ell_1\) geometry and the other by Gaussian chaining, or
equivalently by the \(\gamma_2\) functional. The best decomposition
gives the right balance between the two.

This paper gives one formula for the full class of independent
symmetric coordinates with log-concave tails. It contains both the
Gaussian and Bernoulli theorems with numerical constants. The
coordinates may have different laws. They may be bounded, and their
tails may become arbitrarily steep. Writing
\[
 N_i(u)=-\log\Prob(|Y_i|\ge u),
\]
we do not assume a \(\Delta_2\) condition
\[
 N_i(2u)\le C_0N_i(u).
\]
The formula uses probability measures on \(T\) itself and a cost read
directly from each coordinate tail. Quadratic and truncated quadratic
costs give the Gaussian and Bernoulli cases.

This change of assumptions also changes the geometry of the problem.
The growth condition in the earlier log-concave theory cannot simply
be dropped while keeping the same moment chain. For example, consider
a log-concave tail with
\[
 N(u)=2^u-1.
\]
Then
\[
 \frac{N(2u)}{N(u)}=2^u+1\longrightarrow\infty,
\]
so the \(\Delta_2\) condition fails. For
\(T_m=\{0,e_1,\ldots,e_m\}\), we show that
\[
 S_Y(T_m)\asympc \log\log(m+e),
 \qquad
 \gamma_Y(T_m)\ge
 c\bigl(\log\log(m+e)\bigr)^2.    
\]
Here \(\gamma_Y\) is the usual chain built from
\(L^{2^k}\)-increment distances. Thus the old chain can be too large
by an unbounded factor. Removing the growth condition therefore
requires a new quantity, not just a different proof of the old one.

The same type of supremum appears naturally in several other settings.
If \(u_i\) belong to a finite-dimensional normed space \(F\), then
\begin{equation}\label{eq:intro-banach}
 \left\|\sum_iY_i u_i\right\|_F
 =\sup_{\|f\|_{F^*}\le1}\sum_iY_i f(u_i).
\end{equation}
Thus a bound for a canonical process also gives a bound for a random
sum in a normed space; this connection is classical and plays an
important role in \cite{LedouxTalagrand}. Likewise, for deterministic
real \(a\times b\) matrices \(A_i\),
\begin{equation}\label{eq:intro-matrix}
 \left\|\sum_iY_iA_i\right\|_{2\to2}
 =\sup_{u\in S^{a-1},\,v\in S^{b-1}}
       \sum_iY_i\,u^{\mathsf T}A_i v.
\end{equation}
In finite dimensions, the suprema in these displays can be approximated
by suprema over finite sets, so our results apply to both. Similar
suprema also arise in empirical processes, statistical estimation,
and random matrix problems; see \cite{Dirksen,Liu}.

We prove more than a comparison of suprema: we keep the distribution
of the selected index fixed. For an integrable centered vector
\(A\in\R^d\) and a finitely supported probability measure \(\mu\) on
\(\R^d\), write
\begin{equation}\label{eq:pairing}
 B_A(\mu)
 =
 \sup_{X\sim\mu,\ A\text{ with its prescribed law}}
 \E\ip{X}{A},
 \qquad
 S_A(T)=\E\max_{x\in T}\ip{x}{A}.
\end{equation}
The supremum in \(B_A(\mu)\) is over all couplings of \(X\) and \(A\),
while the law of \(A\) is kept fixed. Choosing a maximizing index gives
\begin{equation}\label{eq:set-pairing}
 S_A(T)=\sup_{\mu\in\cP(T)}B_A(\mu).
\end{equation}
This is the canonical-process form of a functional introduced by
Fernique \cite{Fernique76,Fernique81}; see also
\cite[Definition~1.4 and Section~2]{vH}.
For example, if \(\mu\) gives mass \(1/2\) to each of \(0\) and \(h\)
and \(A\) is symmetric, then
\[
 B_A(\mu)=\tfrac12\E|\ip{h}{A}|.
\]
Choose \(h\) when \(\ip{h}{A}>0\), choose \(0\) when it is negative,
and split the event of equality to obtain the required masses. Thus,
for a two-point index law, \(B_A(\mu)\) is half the expected size of
a single increment. For a general \(\mu\), it also keeps track of how
often each index is selected.

We use \(X\) for a random index and \(Y\) for the source in the main
theorem. We write \(\cP(T)\) for the probability measures on a finite
set \(T\), and \(v\push\mu\) for the law of \(v(X)\) when \(X\sim\mu\).
Throughout the paper, \(c,C\) denote positive universal constants,
whose values may change from line to line, and we write
\(a\asympc b\) when \(cb\le a\le Cb\).

\subsection{The main theorem}

Recall that
\[
 N_i(u)=-\log\Prob(|Y_i|\ge u),\qquad u\ge0.
\]
Assume first that \(Y_1,\ldots,Y_n\) are independent, symmetric,
and nonzero, and that each \(N_i\) is a closed, convex,
nondecreasing extended-valued function with \(N_i(0)=0\).
If \(Y_i\) is bounded, then \(N_i(u)=+\infty\) beyond the upper
endpoint of \(|Y_i|\), while at the endpoint
\[
 N_i(u)=\lim_{v\uparrow u}N_i(v)\in[0,+\infty].
\]
A finite endpoint value means that the law has positive mass there.
Thus bounded laws, including symmetric two-point laws, are allowed.

Put
\[
 q_i(p)=\sup\{u\ge0:N_i(u)\le p\},\qquad p\ge0.
\]
The function \(q_i\) is finite, nonnegative, nondecreasing, and concave.
It is continuous on \((0,\infty)\) and may become constant. Its value
at zero need not vanish. If $E$ is a
mean-one exponential variable, then
\begin{equation}\label{eq:quantile-representation}
 q_i(E)\stackrel{\rm d}=\abs{Y_i}.
\end{equation}
We rescale each coordinate, and the corresponding index coordinate,
so that $q_i(1)=1$. Then $q_i(p)\le\max(1,p)$ and
$b_i:=\lim_{p\to\infty}q_i(p)\in[1,\infty]$.
If $b_i>1$, continuity and inversion give $N_i(1)=1$.
Define the even potential
\begin{equation}\label{eq:potential}
 V_i(u)=
 \begin{cases}
 u^2,&\abs u\le1,\\
 2N_i(\abs u)-1,&1<\abs u\le b_i,\\
 +\infty,&\abs u>b_i.
 \end{cases}
\end{equation}
If $b_i=1$, instead set $V_i(u)=u^2$ on $[-1,1]$ and $+\infty$
outside. In this case $\abs{Y_i}\le1$ and
$\E\abs{Y_i}\ge e^{-1}$, so the convention changes comparisons by
only a numerical constant. Convexity at the join in
\eqref{eq:potential} follows from $N'_{i,+}(1)\ge1$.

The quadratic part of \(V_i\) fixes the behavior near the origin.
Beyond one, \(V_i\) records the original tail, up to a fixed numerical
change. This convention also makes the constants independent of the
tail at the joining point.

We use $+$ and $-$ on derivatives for the right and left
derivatives, respectively. For $0<u<b_i$ write $\sigma_i(u)=V'_{i,+}(u)$.
The scalar cost is
\begin{equation}\label{eq:intrinsic-cost}
 \Psi_i(r)=\int_0^{b_i}\sigma_i(u)\ind_{\{\sigma_i(u)<r\}}\,du,
 \qquad r\ge0.
\end{equation}
In this formula, an interval where $V_i$ has slope $s$ contributes
its slope times its length once $r>s$. Thus $\Psi_i$ records where
the tail becomes steeper. It need not be a distance or a convex
function. Its role is to assign a nonnegative cost to a coordinate
error. The integral is the definition even when it is infinite.
If $f_i(r)=V_i^*(r):=\sup_u\{ru-V_i(u)\}$ is the convex
conjugate of $V_i$, it agrees with $rf'_{i,-}(r)-f_i(r)$ whenever that
expression is defined and finite. The strict inequality in
\eqref{eq:intrinsic-cost} fixes the value at a slope threshold. We have
\begin{equation}\label{eq:core}
 \Psi_i(r)=r^2/4\quad(0\le r\le2),\qquad
 \Psi_i(r)\ge\chi(r/2),\qquad \chi(r)=\min(1,r^2).
\end{equation}

For $t>0$ define
\begin{align}
 D_t(x,y)&=\sum_{i=1}^n\Psi_i(\abs{x_i-y_i}/t),\label{eq:D}\\
 F_t(x,\nu)&=-\log\sum_{y\in T}\nu(y)e^{-D_t(x,y)},
                      \qquad \nu\in\cP(T).\label{eq:F}
\end{align}
The two deterministic functionals in our theorem are
\begin{align}
 \Gamma^\Psi(T)&=\inf_{(\nu_k)\in\cP(T)^\Z}\max_{x\in T}
              \sum_{k\in\Z}2^{-k}F_{2^{-k}}(x,\nu_k),
                         \label{eq:Gamma-definition}\\
 \cM^\Psi(T)&=\inf_{\nu\in\cP(T)}\max_{x\in T}
                  \int_0^\infty F_t(x,\nu)\,dt.
                         \label{eq:single-measure}
\end{align}
In $\Gamma^\Psi$ one may choose a measure at each scale. In
$\cM^\Psi$ the same measure is used at every scale. All these
measures are supported on the original set $T$. The definition
of $\Gamma^\Psi$ uses measures, not partitions; we will compare
it with the usual $\gamma_2$ in the Gaussian case below.

To state the assertion for a fixed index law, set
\begin{align}
 \RD^\Psi_\mu(t)&=\inf_{X,X'\sim\mu}
                \{I(X;X')+\E D_t(X,X')\},\label{eq:RD}\\
 \cR^\Psi(\mu)&=\int_0^\infty\RD^\Psi_\mu(t)\,dt,
 \qquad \cR^\Psi(T)=\sup_{\mu\in\cP(T)}\cR^\Psi(\mu).
                         \label{eq:area}
\end{align}
The variable $X'$ is a randomized approximation to $X$. The term
$I(X;X')$ measures how much the approximation reveals about $X$;
$\E D_t(X,X')$ measures its error at scale $t$. We minimize their
sum. Both marginals in \eqref{eq:RD} equal $\mu$. This last
requirement will be important in the proof. The functional
$\cR^\Psi(\mu)$ is obtained by integrating this minimum over all scales $t>0$.
All logarithms are natural.
More explicitly, if $\mu=\sum_{j=1}^m\mu_j\delta_{x_j}$ has positive
masses at distinct points, the admissible matrices are
\[
 \Pi(\mu,\mu)=\left\{P\in[0,1]^{m\times m}:
           \sum_kP_{jk}=\mu_j,\quad\sum_jP_{jk}=\mu_k\right\},
\]
and
\begin{equation}\label{eq:RD-matrix}
 \RD^\Psi_\mu(t)=\min_{P\in\Pi(\mu,\mu)}
 \left\{\sum_{j,k}P_{jk}\log\frac{P_{jk}}{\mu_j\mu_k}
                    +\sum_{j,k}P_{jk}D_t(x_j,x_k)\right\}.
\end{equation}
The first sum is $I(X;X')$ for the joint law
$\Prob(X=x_j,X'=x_k)=P_{jk}$.
We use $0\log0=0$, $0\cdot\infty=0$, and $e^{-\infty}=0$.
Put $D_\infty=\diam_\infty(\supp\mu)$ and
$D_2=\diam_2(\supp\mu)$. The diagonal and independent couplings give
\begin{equation}\label{eq:elementary-RD-bounds}
 0\le\RD^\Psi_\mu(t)\le\Ent(\mu):=-\sum_x\mu(x)\log\mu(x),
 \qquad
 \RD^\Psi_\mu(t)\le\frac{D_2^2}{4t^2}\quad(t\ge D_\infty/2).
\end{equation}
These bounds show directly that $\cR^\Psi(\mu)<\infty$. Section~\ref{sec:information} gives the information
identities needed later, including their short proofs.

\begin{theorem}\label{thm:main}
Under the assumptions above, for every finitely supported probability
measure \(\mu\) on \(\R^n\),
\begin{equation}\label{eq:main-fixed}
 B_Y(\mu)\asympc \cR^\Psi(\mu).
\end{equation}
For every nonempty finite set \(T\subset\R^n\),
\begin{equation}\label{eq:main-set}
 S_Y(T)
 \asympc
 \Gamma^\Psi(T)
 \asympc
 \cM^\Psi(T)
 \asympc
 \cR^\Psi(T).
\end{equation}
The constants are universal.
\end{theorem}

Deterministic zero coordinates can be deleted. Unnormalized coordinates
are treated by replacing $Y_i$ by $Y_i/q_i(1)$ and $x_i$ by
$q_i(1)x_i$. A singleton index set has value zero throughout.
The single measure in \eqref{eq:single-measure} is supported on the
original finite set.

\subsection{Gaussian and Bernoulli processes}

The two classical cases are recovered from Theorem~\ref{thm:main}.
Let
\[
 G=(g_i)_{i=1}^n
\]
have independent standard Gaussian coordinates, and let
\[
 \eps=(\eps_i)_{i=1}^n
\]
have independent uniform sign coordinates. Put
\[
 \psi_{\rm G}(r)=r^2,
 \qquad
 \psi_{\rm B}(r)=\min(1,r^2).
\]
For either scalar cost, we use the definitions above with the same
cost in every coordinate.

\begin{corollary}\label{cor:gaussian-bernoulli}
For every finitely supported probability measure \(\mu\) on \(\R^n\),
\[
 B_G(\mu)\asympc \cR^{\psi_{\rm G}}(\mu),
 \qquad
 B_{\eps}(\mu)\asympc \cR^{\psi_{\rm B}}(\mu).
\]
For every nonempty finite set \(T\subset\R^n\),
\begin{align}
 S_G(T)
 &\asympc
 \cR^{\psi_{\rm G}}(T)
 \asympc
 \Gamma^{\psi_{\rm G}}(T)
 \asympc
 \cM^{\psi_{\rm G}}(T)
 \asympc
 \gamma_2(T),                                      \label{eq:gaussian-case}\\
 S_{\eps}(T)
 &\asympc
 \cR^{\psi_{\rm B}}(T)
 \asympc
 \Gamma^{\psi_{\rm B}}(T)
 \asympc
 \cM^{\psi_{\rm B}}(T)
 \asympc
 b(T),                                             \label{eq:bernoulli-case}
\end{align}
where
\[
 \gamma_2(T)
 =
 \inf_{(\mathcal A_k)}
 \max_{x\in T}
 \sum_{k\ge0}
 2^{k/2}\diam_2\mathcal A_k(x),
\]
the infimum being over nested partitions with
\[
 \mathcal A_0=\{T\},
 \qquad
 |\mathcal A_k|\le 2^{2^k},
\]
and
\begin{equation}\label{eq:bernoulli-decomposition}
 b(T)
 =
 \inf_{T\subset T_1+T_2}
 \left\{
   \sup_{x\in T_1}\|x\|_1+\gamma_2(T_2)
 \right\}.
\end{equation}
\end{corollary}

The last equivalences in
\eqref{eq:gaussian-case} and \eqref{eq:bernoulli-case}
are the classical Gaussian and Bernoulli theorems
\cite{Talagrand,BL}. For signs, the potential is \(u^2\) on
\([-1,1]\) and \(+\infty\) outside, so the intrinsic cost is
\[
 \Psi(r)=\min\left(1,\frac{r^2}{4}\right)
        =\psi_{\rm B}(r/2).
\]
For Gaussian coordinates, the intrinsic cost is comparable to
\(r^2\) with numerical constants. Thus the same construction gives
the two classical geometries. Neither classical identification is
used in the proof of Theorem~\ref{thm:main}.

\subsection{Earlier results and proof strategy} Talagrand obtained two-sided bounds for canonical processes with tails \(\exp(-u^r)\), \(1\le r<\infty\) \cite{Talagrand94}. Lata\l a extended this to independent log-concave tails satisfying an upper growth condition \[ N_i(2u)\le C_0N_i(u) \] \cite[Theorem~3]{Latala97}. Lata\l a and Tkocz later described canonical processes in terms of chains of increment moments under regular moment growth and a positive lower speed of growth \cite{LT}. The example given above shows that the moment chain from the regular growth theory cannot simply be kept after the growth assumption is removed. The difficulty is therefore not only to extend an earlier proof. One first has to identify the right geometry for a tail whose slope may change arbitrarily. The functions \(\Psi_i\) provide this geometry. They are determined directly by the tail potentials \(V_i\). If \[ C_p=\left\{z\in\R^n:\sum_iV_i(z_i)\le p\right\}, \qquad \rho_p(h)=\sup_{z\in C_p}\ip{h}{z}, \] then we prove \begin{equation}\label{eq:intro-support} \rho_p(h) = \int_0^\infty \min\{p,D_t(0,h)\}\,dt. \end{equation} Thus the scale costs \(D_t\) recover the support geometry of the tail body. For Gaussian coordinates this gives the quadratic geometry, and for signs it gives the truncated quadratic geometry. For a general tail, changes in its slope enter at the scales where they occur. This changes the way the general case is reduced to a simpler one. Instead of imposing regular growth on the whole tail, we split each tail into intervals on which its slope changes by at most a fixed factor. These pieces are represented by truncated exponential coordinates. At the smallest truncation level the corresponding cost is of
Bernoulli type, while without truncation one obtains the symmetric
exponential process. The first step is therefore a uniform comparison for truncated exponentials. Let \[ Z=(\eps_iE_i)_{i=1}^d,\qquad \xi^{[\beta]} = (\eps_i\min(E_i,\beta_i))_{i=1}^d, \] where the signs and the mean-one exponentials are independent. We prove \begin{equation}\label{eq:intro-capping} W\cx Z,\qquad |W_i|\le\beta_i \quad\Longrightarrow\quad W\cx C\xi^{[\beta]}. \end{equation} The coordinates of \(W\) need not be independent and the truncation levels \(\beta_i\) may all be different. The main point of the proof is that the error produced by truncation, after a suitable rescaling, satisfies the same convex-order and coordinatewise bounds as the original vector. It can therefore be returned to the same comparison and absorbed. The second step transfers a general tail to this basic model. We construct a deterministic linear map \(L\), depending only on the coordinate tails, and independent truncated exponential coordinates \(\Xi\) such that \begin{equation}\label{eq:intro-transfer} B_Y(\mu)\asympc B_\Xi(L\push\mu), \qquad \cR^\Psi(\mu)\asympc \cR^{\rm cap}(L\push\mu). \end{equation}
Here \(\cR^{\rm cap}\) denotes the rate--distortion functional
associated with the coordinatewise capped-exponential costs; its
precise definition is given in Section~4. The same map therefore compares both the random process and the quantity used to describe its geometry. This is what allows the reduction to truncated exponentials to preserve the majorizing-measure problem itself, rather than only the random sums. For each original coordinate, the coefficients created by the decomposition have total sum at most three. Hence the constants do not depend on the number or size of the changes of tail slope. For comparison, the Bernoulli theorem of Bednorz and Lata\l a is proved through a decomposition of the index set by admissible partitions \cite{BL}. Liu developed a lifting method that turns bounds for stochastic-process suprema into fixed-law and information-theoretic statements \cite{Liu}. Liu and Zadik recently used an information-theoretic argument based on Bayesian estimation in a Cauchy channel to prove the Bernoulli theorem and obtained a prescribed-law characterization \cite{LZ}. We also use the fixed-law point of view. In the present setting, however, the coordinate geometry is not fixed in advance: the functions \(\Psi_i\) must first be obtained from the individual tails. Once this geometry has been identified, the fixed-law formulation allows the different scales to be combined through the rate--distortion functional. The only external process comparison used in the proof of Theorem~\ref{thm:main} is the exponential case of Lata\l a and Tkocz \cite{LT}. Section~\ref{sec:convex-order} proves the truncation comparison. Section~\ref{sec:geometry} identifies the scale costs and their information form. Section~\ref{sec:general-tails} transfers the result to arbitrary tails, and Section~\ref{sec:one-measure} proves the \(\Gamma\) and single-measure formulas and the examples.

\section*{Tool and computational resource disclosure}

The conceptual framework and the mathematical ideas of this work were
developed by the authors. Generative AI tools were used for language
editing, exposition, and consistency checks during manuscript
preparation. All mathematical arguments and references were
independently verified by the authors, who take full responsibility
for the paper.

\section{Convex order and truncated exponentials}\label{sec:convex-order}

The purpose of this section is to prove the truncation comparison \eqref{eq:intro-capping}. We use \(d\) for the dimension in this section, since the dimension will later change when the original coordinates are expanded. Throughout, \[ Z=(\eps_iE_i)_{i=1}^d \] denotes a vector of independent standard symmetric exponential variables. Thus \[ \Prob(\eps_i=1)=\Prob(\eps_i=-1)=\frac12, \qquad \Prob(E_i>u)=e^{-u},\quad u\ge0, \] and all the signs and exponential variables are independent.

For integrable random vectors \(A,B\in\R^d\), we write \(A\cx B\) if \[ \E f(A)\le \E f(B) \] for every real-valued convex function \(f:\R^d\to\R\). Every such function has an affine minorant. Since \(A\) and \(B\) are integrable, the negative parts of \(f(A)\) and \(f(B)\) are therefore integrable, so these expectations are well defined in \((-\infty,\infty]\). Since affine functions and their negatives are convex, \(A\cx B\) implies \[ \E A=\E B. \] For \(\beta=(\beta_i)_{i=1}^d\in(0,\infty]^d\), define coordinatewise truncation by \[ \bigl(\clip_\beta x\bigr)_i = \max\{-\beta_i,\min\{x_i,\beta_i\}\}. \] When \(\beta\) is a scalar, we use the same notation with \(\beta_i=\beta\) for every \(i\).

The proof has three parts. We first relate convex order to the prescribed-law functional \(B_A(\mu)\) and recall the replication argument that converts comparisons for finite index sets into comparisons for prescribed index laws. We then use the symmetric exponential process as a reference process and prove a contraction estimate for a small residual. Finally, we construct a martingale coupling for which the residual created by truncation satisfies the hypotheses of that estimate. These ingredients are combined at the end of the section to prove the uniform truncation theorem.

We state the contraction estimate first and prove it below after the
required comparison results have been prepared.

\begin{theorem}\label{thm:amplitude} 
Let \((U,R)\) be jointly defined integrable vectors in \(\R^d\). Suppose \begin{equation}\label{eq:amplitude-assumptions} U\cx Z,\qquad R\cx\eta,\qquad \abs{R_i}\le\epsilon(1+\abs{U_i}) \quad\text{almost surely}, \end{equation} where \(0<\epsilon\le1\). Here \[ \eta=(\eta_i)_{i=1}^d, \qquad \eta_i=\delta_iZ'_i, \] \(Z'=(Z'_i)_{i=1}^d\) is an independent copy of \(Z\), and \(\delta_1,\ldots,\delta_d\) are independent \(\{0,1\}\)-valued random variables, independent of \(Z'\), such that \[ \Prob(\delta_i=1)\le q\le1. \] Then \begin{equation}\label{eq:amplitude-result} R\cx C(\sqrt\epsilon+\sqrt q)\,Z. \end{equation}
\end{theorem}

The two assumptions on \(R\) in \eqref{eq:amplitude-assumptions} give different kinds of control. The relation \(R\cx\eta\) compares \(R\) with a sparse exponential vector, while the pointwise bound controls the size of each coordinate relative to the corresponding coordinate of \(U\). Theorem~\ref{thm:amplitude} combines these two bounds into a single convex-order comparison with the reference vector \(Z\).

\subsection{Convex order and prescribed index laws}

We begin with the martingale form of Strassen's theorem \cite{Strassen}. It gives the coupling interpretation of convex order that will be used throughout this section. \begin{theorem}\label{thm:strassen} Let \(A,B\) be integrable random vectors in \(\R^d\). Then \(A\cx B\) if and only if there exist random vectors \((\widetilde A,\widetilde B)\) on a common probability space such that \[ \widetilde A\stackrel{\rm d}=A, \qquad \widetilde B\stackrel{\rm d}=B, \qquad \E[\widetilde B\mid\widetilde A]=\widetilde A \quad\text{almost surely}. \] \end{theorem} 
A coupling with the last property is called a martingale coupling. The implication from a martingale coupling to convex order follows immediately from conditional Jensen: \[ f(\widetilde A) = f\bigl(\E[\widetilde B\mid\widetilde A]\bigr) \le \E[f(\widetilde B)\mid\widetilde A]. \] The converse is the content of Strassen's theorem.

We will also need to add other random variables to such a coupling. Fix a martingale coupling of \(A\) and \(B\), and let \(K(a,db)\) be the conditional law of \(B\) given \(A=a\). Suppose that \((X,A)\) is any other coupling. Conditionally on \((X,A)\), draw a new variable \(\widetilde B\) from \(K(A,db)\), independently of \(X\) once \(A\) is given. Then 
\[ \widetilde B\stackrel{\rm d}=B, \qquad \E[\widetilde B\mid X,A] = \int b\,K(A,db) = A. \]
Thus the original joint law of \((X,A)\) is preserved while the martingale coupling is added. We will use this extension repeatedly.

\begin{lemma}\label{lem:convex-pairing}
Let \(A,B\) be integrable centered vectors.

\begin{enumerate}[label=\textup{(\roman*)}]
\item
\(A\cx B\) if and only if
\[
 B_A(\mu)\le B_B(\mu)
\]
for every finite index law \(\mu\).

\item
Suppose that
\[
 A_j\cx B_j,\qquad 1\le j\le m,
\]
that \(A_1,\ldots,A_m\) are independent, and that
\(B_1,\ldots,B_m\) are independent. Then
\[
 (A_1,\ldots,A_m)\cx(B_1,\ldots,B_m).
\]
Moreover, if \(A\cx B\) and \(L\) is a deterministic linear map, then
\[
 LA\cx LB.
\]

\item
If \(A\) is a symmetric real random variable and \(\eps\) is a
symmetric sign, then
\[
 (\E\abs A)\eps\cx A.
\]

\item
If \(A,B\) are symmetric real random variables and
\[
 \abs A\st\abs B,
\]
then
\[
 A\cx B.
\]
Here \(\abs A\st\abs B\) means
\[
 \Prob(\abs A>t)\le\Prob(\abs B>t),
 \qquad t\ge0.
\]
\end{enumerate}
\end{lemma}

\begin{proof}
\smallskip
\noindent\emph{(i)}
Suppose first that \(A\cx B\). By Theorem~\ref{thm:strassen}, we may
couple \(A\) and \(B\) so that
\[
 \E[B\mid A]=A.
\]
Now start with an arbitrary coupling of \(X\sim\mu\) and \(A\). Using
the martingale kernel described above, construct \(B\) conditionally
on \(A\) so that its prescribed law is preserved and
\[
 \E[B\mid X,A]=A.
\]
Therefore
\[
 \E\ip{X}{A}
 =
 \E\ip{X}{\E[B\mid X,A]}
 =
 \E\ip{X}{B}
 \le B_B(\mu).
\]
Taking the supremum over all couplings of \(X\) and \(A\) gives
\[
 B_A(\mu)\le B_B(\mu).
\]

Conversely, assume that
\[
 B_A(\mu)\le B_B(\mu)
\]
for every finite index law \(\mu\). We first compare finite maxima of
affine functions. Let
\[
 g(v)=\max_{1\le j\le m}
 \bigl(\ip{x_j}{v}+b_j\bigr).
\]
If two vectors \(x_j\) are equal, we may keep only the larger of the
corresponding constants \(b_j\). Thus we may assume that the
\(x_j\)'s are distinct.

We claim that
\begin{equation}\label{eq:affine-max}
 \E g(A)
 =
 \sup_{\lambda\in\cP(\{1,\ldots,m\})}
 \left[
 B_A\left(\sum_{j=1}^m\lambda_j\delta_{x_j}\right)
 +\sum_{j=1}^m\lambda_jb_j
 \right].
\end{equation}
To prove the upper bound, fix
\(\lambda=(\lambda_1,\ldots,\lambda_m)\), and let \(J\) be a random
index with
\[
 \Prob(J=j)=\lambda_j.
\]
For every coupling of \(J\) and \(A\),
\[
 \ip{x_J}{A}+b_J
 \le
 \max_{1\le j\le m}
 \bigl(\ip{x_j}{A}+b_j\bigr)
 =
 g(A).
\]
Hence
\[
 \E\ip{x_J}{A}+\sum_{j=1}^m\lambda_jb_j
 \le
 \E g(A).
\]
Optimizing over all couplings of \(J\) and \(A\) gives
\[
 B_A\left(\sum_{j=1}^m\lambda_j\delta_{x_j}\right)
 +\sum_{j=1}^m\lambda_jb_j
 \le
 \E g(A).
\]

For the reverse inequality, define \(J\) to be the first index at
which the maximum in the definition of \(g(A)\) is attained, and let
\[
 \lambda_j=\Prob(J=j).
\]
Then
\[
 g(A)=\ip{x_J}{A}+b_J
\]
almost surely, and therefore
\[
 \E g(A)
 =
 \E\ip{x_J}{A}
 +\sum_{j=1}^m\lambda_jb_j
 \le
 B_A\left(\sum_{j=1}^m\lambda_j\delta_{x_j}\right)
 +\sum_{j=1}^m\lambda_jb_j.
\]
This proves \eqref{eq:affine-max}.

Applying \eqref{eq:affine-max} to both \(A\) and \(B\), and using
\[
 B_A(\mu)\le B_B(\mu)
\]
for every finite index law \(\mu\), we obtain
\[
 \E g(A)\le \E g(B)
\]
for every finite maximum \(g\) of affine functions.

It remains to pass from such functions to an arbitrary real-valued
convex function \(f:\R^d\to\R\). Choose a countable dense set
\(\{v_k:k\ge1\}\subset\R^d\), and at each \(v_k\) choose a supporting
affine function
\[
 \ell_k(v)\le f(v),
 \qquad
 \ell_k(v_k)=f(v_k).
\]
Since a finite convex function on \(\R^d\) is continuous,
\[
 f(v)=\sup_{k\ge1}\ell_k(v),\qquad v\in\R^d.
\]
Set
\[
 g_m(v)=\max_{1\le k\le m}\ell_k(v).
\]
Then
\[
 g_m(v)\uparrow f(v).
\]
Because \(g_m-\ell_1\ge0\), monotone convergence gives
\[
 \E[f(A)-\ell_1(A)]
 =
 \lim_{m\to\infty}\E[g_m(A)-\ell_1(A)],
\]
and similarly for \(B\). Since \(A\) and \(B\) are centered,
\[
 \E\ell_1(A)=\E\ell_1(B).
\]
Using
\[
 \E g_m(A)\le\E g_m(B)
\]
and letting \(m\to\infty\), we obtain
\[
 \E f(A)\le\E f(B).
\]
Hence \(A\cx B\).

\smallskip
\noindent\emph{(ii)}
For each \(1\le j\le m\), choose by
Theorem~\ref{thm:strassen} a martingale coupling
\((\widetilde A_j,\widetilde B_j)\) such that
\[
 \widetilde A_j\stackrel{\rm d}=A_j,
 \qquad
 \widetilde B_j\stackrel{\rm d}=B_j,
 \qquad
 \E[\widetilde B_j\mid\widetilde A_j]
 =
 \widetilde A_j.
\]
Choose these \(m\) couplings independently. Then
\[
 \E\bigl[
   (\widetilde B_1,\ldots,\widetilde B_m)
   \,\big|\,
   \widetilde A_1,\ldots,\widetilde A_m
 \bigr]
 =
 (\widetilde A_1,\ldots,\widetilde A_m).
\]
Hence, by Theorem~\ref{thm:strassen},
\[
 (A_1,\ldots,A_m)\cx(B_1,\ldots,B_m).
\]

Now suppose that \(A\cx B\), and let \(L\) be a deterministic linear
map. Choose a martingale coupling such that
\[
 \E[B\mid A]=A.
\]
Since \(LA\) is a function of \(A\), conditioning first on \(A\) gives
\[
 \E[LB\mid LA]
 =
 \E\bigl[\E[LB\mid A]\mid LA\bigr]
 =
 \E[LA\mid LA]
 =
 LA.
\]
Therefore, by Theorem~\ref{thm:strassen},
\[
 LA\cx LB.
\]

\smallskip
\noindent\emph{(iii)}
By symmetry, \(A\) has the same distribution as
\[
 \eps\abs A,
\]
where \(\eps\) is a symmetric sign independent of \(\abs A\).
Set
\[
 U=(\E\abs A)\eps,
 \qquad
 V=\eps\abs A.
\]
If \(\E\abs A=0\), then \(A=0\) almost surely and the conclusion is
immediate. Otherwise \(U\) determines the sign \(\eps\), and hence
\[
 \E[V\mid U]
 =
 \eps\,\E\abs A
 =
 U.
\]
Since \(V\stackrel{\rm d}=A\), Theorem~\ref{thm:strassen} gives
\[
 (\E\abs A)\eps\cx A.
\]

\smallskip
\noindent\emph{(iv)}
Let \(f:\R\to\R\) be convex, and define
\[
 \phi(t)=\frac{f(t)+f(-t)}2,
 \qquad t\ge0.
\]
The function \(\phi\) is convex and nondecreasing on
\([0,\infty)\). Since \(A\) and \(B\) are symmetric,
\[
 \E f(A)=\E\phi(\abs A),
 \qquad
 \E f(B)=\E\phi(\abs B).
\]
The assumption
\[
 \abs A\st\abs B
\]
therefore gives
\[
 \E\phi(\abs A)
 \le
 \E\phi(\abs B).
\]
Thus
\[
 \E f(A)\le\E f(B)
\]
for every convex \(f\), and consequently
\[
 A\cx B.
\]
\end{proof}

We will also use a simple compatibility property of \(B_A\) with maps
of the index set. Let \(v:\supp\mu\to\R^d\). Then
\begin{equation}\label{eq:fiber}
 B_A(v\push\mu)
 =
 \sup_{\substack{X\sim\mu\\ A\text{ with its prescribed law}}}
 \E\ip{v(X)}A.
\end{equation}
Indeed, every coupling of \(X\sim\mu\) and \(A\) induces a coupling of
\(v(X)\sim v\push\mu\) and \(A\), so the right-hand side is at most
the left-hand side.

For the converse, start with any coupling \((Y,A)\) such that
\(Y\sim v\push\mu\). For \(y\in\supp(v\push\mu)\), conditionally on
\((Y,A)\) draw \(X\) from the fiber \(v^{-1}(y)\) according to
\[
 \Prob(X=x\mid Y=y,A)
 =
 \frac{\mu(x)\ind_{\{v(x)=y\}}}{(v\push\mu)(y)}.
\]
Then \(X\sim\mu\), \(v(X)=Y\) almost surely, and
\[
 \E\ip{v(X)}A=\E\ip YA.
\]
Thus every coupling on the left can be lifted to one on the right,
which proves \eqref{eq:fiber}. In particular, no injectivity of \(v\)
is required.

The next lemma is the replication principle used to pass from ordinary
suprema on finite sets to the fixed-law functional \(B_A(\mu)\).
We follow the assignment argument in
\cite[Proposition~3.1 and Section~3.3]{vH}, where the principle is
attributed to Liu \cite{Liu}.

\begin{lemma}\label{lem:type}
Let \(A\) be an integrable centered vector and let
\[
 \mu=\sum_{j=1}^m\mu_j\delta_{x_j}
\]
have positive rational weights. For every integer \(N\) such that
\(N\mu_j\in\mathbb N\) for all \(j\), let \(T_N(\mu)\) be the set of
ordered \(N\)-tuples
\[
 (t_1,\ldots,t_N)\in(\R^d)^N
\]
containing exactly \(N\mu_j\) copies of \(x_j\) for each \(j\).
Let
\[
 A^{\oplus N}=(A^{(1)},\ldots,A^{(N)}),
\]
where \(A^{(1)},\ldots,A^{(N)}\) are independent copies of \(A\).
Then
\begin{equation}\label{eq:type-limit}
 \lim_{N\to\infty}\frac1N
 S_{A^{\oplus N}}(T_N(\mu))
 =
 B_A(\mu),
\end{equation}
where the limit is taken along integers for which all \(N\mu_j\) are
integral.

Moreover, suppose that for some \(a\ge0\),
\[
 S_{A^{\oplus N}}(T)
 \le
 a\,S_{Z^{\oplus N}}(T)
\]
for every \(N\) and every finite \(T\subset(\R^d)^N\), with the same
constant \(a\). Then
\[
 B_A(\mu)\le aB_Z(\mu)
\]
for every finite index law \(\mu\), and consequently
\[
 A\cx aZ.
\]
\end{lemma}

\begin{proof}
We first rewrite
\[
 \frac1N S_{A^{\oplus N}}(T_N(\mu))
\]
as an assignment problem. For probability measures \(\alpha,\nu\) on \(\R^d\) with
finite first moments, define
\[
 P_\alpha(\nu)
 =
 \sup_{\substack{X\sim\alpha\\ V\sim\nu}}
 \E\ip XV.
\]
In particular,
\[
 P_\mu(\Law(A))=B_A(\mu).
\]

Fix deterministic vectors \(a_1,\ldots,a_N\), and write
\[
 \nu_N^a=\frac1N\sum_{k=1}^N\delta_{a_k}.
\]
List the \(N\) points prescribed by \(\mu\) as
\[
 \widetilde x_1,\ldots,\widetilde x_N,
\]
with exactly \(N\mu_j\) copies of \(x_j\). After treating repeated
points as labeled atoms of mass \(1/N\), every coupling of
\(\mu\) and \(\nu_N^a\) is represented by a matrix
\(M=(M_{\ell k})\) with
\[
 M_{\ell k}\ge0,\qquad
 \sum_kM_{\ell k}=1,\qquad
 \sum_\ell M_{\ell k}=1.
\]
Hence
\[
 P_\mu(\nu_N^a)
 =
 \frac1N
 \max_{M\in\mathcal B_N}
 \sum_{\ell,k=1}^N
 M_{\ell k}\ip{\widetilde x_\ell}{a_k},
\]
where
\[
 \mathcal B_N
 =
 \{M\ge0:
 M\mathbf1=M^{\mathsf T}\mathbf1=\mathbf1\}.
\]

By the Birkhoff--von Neumann theorem \cite{Birkhoff},
\[
 \mathcal B_N
 =
 \operatorname{conv}
 \{P_\sigma:\sigma\in\mathfrak S_N\},
\]
where
\[
 (P_\sigma)_{\ell k}
 =
 \ind_{\{k=\sigma(\ell)\}}.
\]
Since the objective is linear in \(M\), its maximum is attained at a
permutation matrix. Therefore
\[
 P_\mu(\nu_N^a)
 =
 \frac1N
 \max_{(t_1,\ldots,t_N)\in T_N(\mu)}
 \sum_{k=1}^N\ip{t_k}{a_k}.
\]
Thus \(P_\mu(\nu_N^a)\) is exactly the best average pairing between the
sample points \(a_k\) and an \(N\)-tuple whose empirical distribution
is \(\mu\).

Now let \(a_k=A^{(k)}\), where the \(A^{(k)}\)'s are independent
copies of \(A\), and set
\[
 \nu_N=\frac1N\sum_{k=1}^N\delta_{A^{(k)}}.
\]
Taking expectations in the preceding identity gives
\[
 \frac1N
 S_{A^{\oplus N}}(T_N(\mu))
 =
 \E P_\mu(\nu_N).
\]
It remains to show that the random empirical measure \(\nu_N\) may be
replaced by \(\Law(A)\).

For probability measures \(\nu,\nu'\) with finite first moments, let
\[
 W_1(\nu,\nu')
 =
 \inf_{\substack{V\sim\nu\\V'\sim\nu'}}
 \E\norm{V-V'}_2
\]
be the \(1\)-Wasserstein distance. Put
\[
 R_\mu=\max_{1\le j\le m}\norm{x_j}_2.
\]
If \((X,V)\) is any coupling with \(X\sim\mu\), and
\((V,V')\) is a coupling of \(\nu,\nu'\), construct \(V'\)
conditionally on \(V\) according to this coupling. Then
\[
 \left|
 \E\ip XV-\E\ip X{V'}
 \right|
 \le
 R_\mu\,\E\norm{V-V'}_2.
\]
It follows that, for every coupling of \(\nu\) and \(\nu'\),
\[
 P_\mu(\nu)
 \le
 P_\mu(\nu')+R_\mu\,\E\norm{V-V'}_2.
\]
Taking the infimum over such couplings gives
\[
 P_\mu(\nu)
 \le
 P_\mu(\nu')+R_\mu W_1(\nu,\nu').
\]
Interchanging \(\nu\) and \(\nu'\) yields
\[
 \bigl|P_\mu(\nu)-P_\mu(\nu')\bigr|
 \le
 R_\mu W_1(\nu,\nu').
\]

We next show
\[
 \E W_1(\nu_N,\Law(A))\longrightarrow0.
\]
Fix \(\varepsilon>0\). Since \(A\) is integrable, we may choose a finite-range measurable map \(Q:\R^d\to\R^d\) such that
\[
 \E\norm{A-Q(A)}_2<\varepsilon,
\]
and put \(\widehat A=Q(A)\).
Apply the same map \(Q\) to each \(A^{(k)}\), and let
\[
 \widehat\nu_N
 =
 \frac1N\sum_{k=1}^N\delta_{\widehat A^{(k)}}.
\]
Then
\[
 \E W_1(\nu_N,\widehat\nu_N)<\varepsilon,
 \qquad
 W_1(\Law(\widehat A),\Law(A))<\varepsilon.
\]

Write
\[
 \Law(\widehat A)
 =
 \sum_{j=1}^r p_j\delta_{a_j},
\]
and let \(\widehat p_{N,j}\) be the empirical mass of \(a_j\) in
\(\widehat\nu_N\). If
\[
 D_A=\max_{j,k}\norm{a_j-a_k}_2,
\]
then matching the common mass at each atom gives
\[
 W_1(\widehat\nu_N,\Law(\widehat A))
 \le
 \frac{D_A}{2}
 \sum_{j=1}^r
 |\widehat p_{N,j}-p_j|.
\]
Since
\[
 \E(\widehat p_{N,j}-p_j)^2
 =
 \frac{p_j(1-p_j)}N,
\]
Cauchy--Schwarz gives
\[
 \E|\widehat p_{N,j}-p_j|
 \le
 \sqrt{\frac{p_j(1-p_j)}N}.
\]
Therefore
\[
 \E W_1(\nu_N,\Law(A))
 \le
 2\varepsilon
 +
 \frac{D_A}{2\sqrt N}
 \sum_{j=1}^r\sqrt{p_j(1-p_j)}.
\]
For fixed \(\widehat A\), the last term tends to zero as
\(N\to\infty\). Hence
\[
 \limsup_{N\to\infty}
 \E W_1(\nu_N,\Law(A))
 \le
 2\varepsilon.
\]
Letting \(\varepsilon\downarrow0\) proves
\[
 \E W_1(\nu_N,\Law(A))\longrightarrow0.
\]

Combining this with the continuity estimate for \(P_\mu\), we obtain
\[
 \frac1N
 S_{A^{\oplus N}}(T_N(\mu))
 =
 \E P_\mu(\nu_N)
 \longrightarrow
 P_\mu(\Law(A))
 =
 B_A(\mu).
\]
This proves \eqref{eq:type-limit}.

We now prove the last assertion. Assume first that \(\mu\) has rational weights. For every admissible
\(N\), the hypothesis applied to \(T_N(\mu)\) gives
\[
 S_{A^{\oplus N}}(T_N(\mu))
 \le
 a\,S_{Z^{\oplus N}}(T_N(\mu)).
\]
Dividing by \(N\) and using \eqref{eq:type-limit} for both \(A\) and
\(Z\), we obtain
\[
 B_A(\mu)\le aB_Z(\mu).
\]

It remains only to remove the rationality assumption. Let \(\mu\) and \(\mu'\) be probability measures on the same finite set \(T\), and put \[ \norm{\mu-\mu'}_{\rm TV} = \frac12\sum_{x\in T}|\mu(x)-\mu'(x)|. \] 
Let
\[
 c(x)=\min\{\mu(x),\mu'(x)\},\qquad x\in T.
\]
Couple the common part by setting
\[
 \Prob(X=X'=x)=c(x).
\]
The remaining measures
\[
 \mu-c,\qquad \mu'-c
\]
have the same total mass and disjoint supports, so their remaining
mass can be coupled with \(X\ne X'\). Hence
\[
 \Prob(X\ne X')
 =
 1-\sum_{x\in T}c(x)
 =
 \norm{\mu-\mu'}_{\rm TV}.
\]

Let
\[
 D_T=\max_{x,x'\in T}\norm{x-x'}_2.
\]
Starting from any coupling of \(X\sim\mu\) and \(A\), extend the
probability space by generating \(X'\) conditionally on \(X\)
according to the coupling above. For every \(M>0\),
\begin{align*}
 \left|
 \E\ip XA-\E\ip{X'}A
 \right|
 &\le
 D_T\,\E\left[
 \norm A_2\ind_{\{X\ne X'\}}
 \right]\\
 &\le
 D_T\left(
 M\norm{\mu-\mu'}_{\rm TV}
 +
 \E\left[
 \norm A_2\ind_{\{\norm A_2>M\}}
 \right]
 \right).
\end{align*}
Taking suprema and then reversing \(\mu,\mu'\) gives the same bound
for
\[
 |B_A(\mu)-B_A(\mu')|.
\]
For fixed \(M\), the first term tends to zero when
\(\mu'\to\mu\), while the second tends to zero as \(M\to\infty\)
because \(A\) is integrable. Hence \(B_A(\mu)\) is continuous in the
masses of \(\mu\). The same argument applies to \(B_Z(\mu)\).

Approximating arbitrary masses by rational ones therefore gives
\[
 B_A(\mu)\le aB_Z(\mu)
 =
 B_{aZ}(\mu)
\]
for every finite index law \(\mu\). Lemma~\ref{lem:convex-pairing}
then yields
\[
 A\cx aZ.
\]
\end{proof}

\subsection{The exponential reference process}

This subsection provides the process comparison needed in the proof of
Theorem~\ref{thm:amplitude}. We first record moment estimates for the
reference vector \(Z\) and for the sparse exponential vector \(\eta\).
We then use a comparison theorem of Lata\l a and Tkocz to pass from
bounds on increment moments to bounds on expected suprema.

\begin{lemma}\label{lem:laplace-moments}
For \(p\ge2\) and \(h\in\R^d\),
\begin{align}
 \norm{\ip hZ}_p
 &\asympc
 \sqrt p\norm h_2+p\norm h_\infty,
                  \label{eq:laplace-moments}\\
 \norm{\sum_i\abs{h_i}E_i}_p
 &\le
 C\bigl(\norm h_1+\sqrt p\norm h_2+p\norm h_\infty\bigr),
                  \label{eq:positive-moments}\\
 \norm{\ip h\eta}_p
 &\le
 C\bigl(\sqrt{qp}\norm h_2+p\norm h_\infty\bigr).
                  \label{eq:sparse-moments}
\end{align}
For \(1\le p\le2\),
\[
 \norm{\ip hZ}_p\asympc\norm h_2.
\]
\end{lemma}

\begin{proof}
We begin with the sparse vector \(\eta\). Put
\[
 p_i=\Prob(\delta_i=1).
\]
Since
\[
 \eta_i=\delta_iZ_i'
\]
and \(Z_i'\) is a standard symmetric exponential variable, for
\(\abs t\norm h_\infty\le1/2\),
\[
 \E e^{t h_i\eta_i}
 =
 1-p_i+\frac{p_i}{1-t^2h_i^2}.
\]
Independence of the coordinates therefore gives
\[
 \log\E e^{t\ip h\eta}
 =
 \sum_i
 \log\left(
 1-p_i+\frac{p_i}{1-t^2h_i^2}
 \right)
 \le
 2t^2\sum_i p_i h_i^2.
\]
Since \(p_i\le q\),
\[
 \sum_i p_i h_i^2\le q\norm h_2^2.
\]
Set
\[
 v^2=q\norm h_2^2,
 \qquad
 b=\norm h_\infty.
\]
For \(0<t\le(2b)^{-1}\), Markov's inequality gives
\[
 \Prob(\ip h\eta>x)
 \le
 \exp(-tx+2t^2v^2).
\]
Choosing
\[
 t=\min\left\{\frac{x}{8v^2},\frac1{4b}\right\},
\]
with the evident interpretation when \(v=0\) or \(b=0\), yields
\[
 \Prob(\ip h\eta>x)
 \le
 \exp\left[
 -c\min\left\{\frac{x^2}{v^2},\frac{x}{b}\right\}
 \right].
\]
Applying the same argument to \(-\ip h\eta\), we obtain
\[
 \Prob\left\{
 \abs{\ip h\eta}>
 C\bigl(\sqrt{qu}\norm h_2+u\norm h_\infty\bigr)
 \right\}
 \le2e^{-u},
 \qquad u\ge0.
\]

For a nonnegative random variable \(V\),
\[
 \E V^p
 =
 p\int_0^\infty
 x^{p-1}\Prob(V>x)\,dx.
\]
Applying this identity to \(V=\abs{\ip h\eta}\) gives
\begin{align*}
 \E V^p
 &\le
 2p\int_0^\infty
 x^{p-1}
 \left(
 e^{-cx^2/v^2}+e^{-cx/b}
 \right)\,dx\\
 &=
 2(v/\sqrt c)^p\Gamma(1+p/2)
 +
 2(b/c)^p\Gamma(1+p),
\end{align*}
where \(\Gamma\) is Euler's gamma function. The elementary estimate
\[
 \Gamma(1+u)\le2(2u/e)^u,
 \qquad u>0,
\]
then gives
\[
 \norm V_p
 \le
 C\bigl(\sqrt p\,v+pb\bigr).
\]
Substituting the definitions of \(v\) and \(b\) proves
\eqref{eq:sparse-moments}.

If all selectors are identically equal to one, then
\(\eta\) has the same law as \(Z\). The same estimate with \(q=1\)
therefore gives the upper bound in
\eqref{eq:laplace-moments}.

For \eqref{eq:positive-moments}, note that for
\(\abs t\le1/2\),
\[
 \log\E e^{t(E_i-1)}
 =
 -t-\log(1-t)
 \le2t^2.
\]
Applying the preceding moment argument to
\[
 \sum_i\abs{h_i}(E_i-1)
\]
and then adding its mean
\[
 \sum_i\abs{h_i}=\norm h_1
\]
gives \eqref{eq:positive-moments}.

It remains to prove the lower bound in
\eqref{eq:laplace-moments}. For every \(i\),
\[
 \E[\ip hZ\mid Z_i]=h_iZ_i,
\]
because the other coordinates are centered. Conditional Jensen
therefore gives
\[
 \norm{\ip hZ}_p
 \ge
 |h_i|\norm{Z_i}_p.
\]
Taking the largest coordinate and using
\[
 \norm{Z_i}_p
 =
 \Gamma(p+1)^{1/p}
 \ge cp,
 \qquad p\ge2,
\]
we obtain
\[
 \norm{\ip hZ}_p
 \ge
 cp\norm h_\infty.
\]

To obtain the Euclidean term, let \(G_i\) be independent standard
Gaussian variables, independent of the exponentials \(E_i\). Since
\[
 \E e^{it\sqrt{2E_i}G_i}
 =
 \E e^{-t^2E_i}
 =
 \frac1{1+t^2}
 =
 \E e^{itZ_i},
\]
we have
\[
 Z_i\stackrel{\rm d}=\sqrt{2E_i}G_i.
\]
Conditioning on the exponentials gives
\[
 \norm{\ip hZ}_p
 =
 \norm{G_1}_p
 \norm{\sqrt{2\sum_i h_i^2E_i}}_p.
\]
For \(p\ge2\),
\[
 \norm{\sqrt{\sum_i h_i^2E_i}}_p
 \ge
 \left(
 \E\sum_i h_i^2E_i
 \right)^{1/2}
 =
 \norm h_2,
\]
and
\[
 \norm{G_1}_p\ge c\sqrt p.
\]
Hence
\[
 \norm{\ip hZ}_p
 \ge
 c\sqrt p\norm h_2.
\]
Together with the preceding \(\ell_\infty\) bound, this proves the
lower estimate in \eqref{eq:laplace-moments} for \(p\ge2\).

Finally, let \(1\le p\le2\) and put
\[
 V=\ip hZ.
\]
Independence and symmetry give
\[
 \E V^2=2\norm h_2^2,
\]
and
\[
 \E V^4
 =
 12\norm h_2^4+12\sum_i h_i^4
 \le
 24\norm h_2^4.
\]
Applying the Paley--Zygmund inequality to \(V^2\) yields
\[
 \Prob\{\abs V\ge\norm h_2\}\ge\frac1{24}.
\]
Consequently,
\[
 \norm V_1\ge\frac1{24}\norm h_2,
 \qquad
 \norm V_2=\sqrt2\norm h_2.
\]
Monotonicity of \(L^p\)-norms completes the proof for
\(1\le p\le2\).
\end{proof}

We next recall the functional used to pass from moment
bounds to expected suprema. For a process \((X_x)_{x\in T}\) on a
finite set \(T\), define
\[
 d_p(x,y)=\norm{X_x-X_y}_p.
\]
Set
\begin{equation}\label{eq:moment-gamma}
 \gamma_X(T)
 =
 \inf_{(\mathcal A_k)}
 \max_{x\in T}
 \sum_{k\ge0}
 \diam_{d_{2^k}}\mathcal A_k(x).
\end{equation}
Here the infimum is over nested partitions satisfying
\[
 \mathcal A_0=\{T\},
 \qquad
 |\mathcal A_k|\le2^{2^k},
 \quad k\ge1,
\]
and \(\mathcal A_k(x)\) denotes the unique cell of
\(\mathcal A_k\) containing \(x\). Since \(T\) is finite, the
partitions may be taken to be singletons from some level onward.

The functional \(\gamma_X\) is used here only as an auxiliary chaining functional built from the increment distances \(d_{2^k}\); it should not be confused with \(\Gamma^\Psi\) in Theorem~\ref{thm:main}. For a random vector \(A\), we write \(\gamma_A(T)\) for the functional associated with the canonical process \[ X_x=\ip xA,\qquad x\in T. \]

We use the following finite-set form of a result of Lata\l a and
Tkocz \cite[Theorem~6]{LT}. To distinguish their growth parameter
from the truncation levels \(\beta_i\), we denote it by
\(\beta_{\rm LT}\). In their terminology, the moments of independent
standardized random variables \(V_i\) grow \(\alpha\)-regularly with
speed \(\beta_{\rm LT}\) if, for some \(\alpha\ge1\) and
\(\beta_{\rm LT}>1\),
\begin{equation}\label{eq:LT-hypotheses}
 \norm{V_i}_p
 \le
 \alpha\frac pq\norm{V_i}_q
 \quad(p\ge q\ge2),
 \qquad
 \norm{V_i}_{\beta_{\rm LT} p}
 \ge
 2\norm{V_i}_p
 \quad(p\ge2).
\end{equation}

\begin{theorem}\label{thm:LT-input}
Let
\[
 X_t=\sum_{i=1}^{\infty}t_iV_i,\qquad t\in\ell^2,
\]
be the canonical process based on independent standardized random
variables \(V_i\) whose moments grow \(\alpha\)-regularly with speed
\(\beta_{\rm LT}\). Then, for every finite nonempty
\(T\subset\ell^2\),
\[
 \frac1{C(\alpha,\beta_{\rm LT})}\gamma_X(T)
 \le
 \E\max_{s,t\in T}(X_s-X_t)
 \le
 C\gamma_X(T).
\]
\end{theorem}

We will also use the following consequence
\cite[Corollary~7]{LT}. Suppose that \(X\) is as in
Theorem~\ref{thm:LT-input}, and that another process
\((Y_t)_{t\in T}\) satisfies
\[
 \norm{Y_s-Y_t}_p
 \le
 \norm{X_s-X_t}_p,
 \qquad
 s,t\in T,\quad p\ge1.
\]
Then
\[
 \E\max_{s,t\in T}(Y_s-Y_t)
 \le
 C(\alpha,\beta_{\rm LT})
 \E\max_{s,t\in T}(X_s-X_t).
\]

We now verify that the symmetric exponential process satisfies the
hypotheses of Theorem~\ref{thm:LT-input}. Since
\[
 \E Z_i=0,
 \qquad
 \E Z_i^2=2,
\]
the variables
\[
 V_i=\frac{Z_i}{\sqrt2}
\]
are standardized. Moreover,
\[
 \E\abs{Z_i}^p=\Gamma(p+1).
\]
For \(p\ge2\),
\[
 \Gamma(p+1)
 =
 \E E_i^p
 \ge
 p^p\Prob(E_i\ge p)
 =
 (p/e)^p,
\]
while the elementary upper bound used above gives
\[
 \Gamma(p+1)\le2(2p/e)^p.
\]
Hence, for \(p\ge q\ge2\),
\[
 \frac{\norm{Z_i}_p}{\norm{Z_i}_q}
 \le
 2e\,\frac pq.
\]
Also,
\[
 \norm{Z_i}_{4ep}
 \ge
 \frac{4ep}{e}
 =
 4p,
\]
whereas
\[
 \norm{Z_i}_p
 \le
 \frac{2p}{e}\,2^{1/p}
 \le
 2p
 \qquad (p\ge2).
\]
Therefore
\[
 \norm{Z_i}_{4ep}
 \ge
 2\norm{Z_i}_p.
\]
Thus \eqref{eq:LT-hypotheses} holds with numerical values of
\(\alpha\) and \(\beta_{\rm LT}\), and the same is true for
\(V_i=Z_i/\sqrt2\).

Both \(\gamma_X(T)\) and the expected range scale linearly when the
process is multiplied by a constant. Since \(Z\) is symmetric,
\[
 \E\max_{x,y\in T}\ip{x-y}Z
 =
 2S_Z(T).
\]
Applying Theorem~\ref{thm:LT-input} to \(Z/\sqrt2\) and rescaling
therefore gives
\begin{equation}\label{eq:reference-comparison}
 \gamma_Z(T)\le CS_Z(T).
\end{equation}

We will use Corollary~7 of \cite{LT} in the following form. Let \(A\)
be a centered random vector; its coordinates need not be independent.
If
\begin{equation}\label{eq:moment-to-set}
 \norm{\ip hA}_p
 \le
 a\norm{\ip hZ}_p,
 \qquad
 h\in\R^d,\quad p\ge1,
\end{equation}
then, for every finite \(T\subset\R^d\),
\[
 S_A(T)\le CaS_Z(T).
\]
Indeed, apply the cited corollary to
\[
 Y_x=\ip xA,
 \qquad
 X_x=a\ip xZ.
\]
The assumption \eqref{eq:moment-to-set} gives
\[
 \norm{Y_x-Y_y}_p
 \le
 \norm{X_x-X_y}_p,
 \qquad
 x,y\in T,\quad p\ge1.
\]
Hence
\[
 \E\max_{x,y\in T}\ip{x-y}A
 \le
 Ca\,\E\max_{x,y\in T}\ip{x-y}Z.
\]
Since \(A\) is centered, for any fixed \(x_0\in T\),
\[
 S_A(T)
 =
 \E\max_{x\in T}\ip{x-x_0}A
 \le
 \E\max_{x,y\in T}\ip{x-y}A.
\]
Using the symmetry of \(Z\),
\[
 \E\max_{x,y\in T}\ip{x-y}Z
 =
 2S_Z(T),
\]
which proves
\[
 S_A(T)\le CaS_Z(T).
\]

\subsection{A contraction estimate}

We now prove Theorem~\ref{thm:amplitude}. The convex-order bound
\(R\cx\eta\) alone is not enough: the estimate
\eqref{eq:sparse-moments} contains the term
\(p\norm h_\infty\), which need not be small when \(q\) is small.
We therefore separate the largest coefficients of \(h\). On these
coordinates we use the pointwise bound on \(R_i\); on the remaining
coordinates, \(\norm h_\infty\) is small enough for the sparse
exponential estimate to be effective.

\begin{proof}[Proof of Theorem~\ref{thm:amplitude}]
Fix \(p\ge2\) and \(h\in\R^d\). Set
\[
 k=\min\left\{d,\left\lceil\frac p\epsilon\right\rceil\right\},
\]
and let \(I\subset\{1,\ldots,d\}\) index the \(k\) largest values of
\(\abs{h_i}\). For \(J\subset\{1,\ldots,d\}\), write
\[
 h_J=(h_i\ind_{\{i\in J\}})_{i=1}^d.
\]

If \(I^c\ne\varnothing\), then
\begin{equation}
 \norm{h_{I^c}}_\infty
 \le
 \sqrt{\frac{\epsilon}{p}}\norm h_2.
\end{equation}
Indeed, every coefficient selected in \(I\) is at least as large as
every coefficient outside \(I\). Hence
\[
 k\norm{h_{I^c}}_\infty^2
 \le
 \sum_{i\in I}|h_i|^2
 \le
 \norm h_2^2.
\]
Since \(I^c\ne\varnothing\), we have
\(k=\lceil p/\epsilon\rceil\ge p/\epsilon\), which proves the claim.

We first estimate the contribution from \(I\). By
\eqref{eq:amplitude-assumptions},
\[
 \abs{\ip{h_I}R}
 \le
 \epsilon\norm{h_I}_1
 +
 \epsilon\sum_{i\in I}|h_i||U_i|.
\]
Therefore
\[
 \norm{\ip{h_I}R}_p
 \le
 \epsilon\norm{h_I}_1
 +
 \epsilon
 \norm{\sum_{i\in I}|h_i||U_i|}_p.
\]
The function
\[
 u\longmapsto
 \left(\sum_{i\in I}|h_i||u_i|\right)^p
\]
is convex. Since \(U\cx Z\),
\[
 \norm{\sum_{i\in I}|h_i||U_i|}_p
 \le
 \norm{\sum_{i\in I}|h_i||Z_i|}_p
 =
 \norm{\sum_{i\in I}|h_i|E_i}_p.
\]
Using \eqref{eq:positive-moments}, we obtain
\[
 \norm{\ip{h_I}R}_p
 \le
 C\epsilon
 \left(
 \norm{h_I}_1
 +
 \sqrt p\,\norm h_2
 +
 p\norm h_\infty
 \right).
\]
Moreover,
\[
 \norm{h_I}_1
 \le
 \sqrt{k}\norm h_2,
\]
and, since \(p/\epsilon\ge2\),
\[
 k\le\left\lceil\frac p\epsilon\right\rceil
 \le
 \frac{2p}{\epsilon}.
\]
Thus
\[
 \norm{\ip{h_I}R}_p
 \le
 C\left(
 \sqrt{\epsilon p}\,\norm h_2
 +
 \epsilon p\norm h_\infty
 \right).
\]

We now estimate the remaining coordinates. Since \(R\cx\eta\), the
linear image property in Lemma~\ref{lem:convex-pairing} gives
\[
 R_{I^c}\cx\eta_{I^c}.
\]
Hence
\[
 \norm{\ip{h_{I^c}}R}_p
 \le
 \norm{\ip{h_{I^c}}\eta}_p.
\]
By \eqref{eq:sparse-moments},
\[
 \norm{\ip{h_{I^c}}R}_p
 \le
 C\left(
 \sqrt{qp}\,\norm{h_{I^c}}_2
 +
 p\norm{h_{I^c}}_\infty
 \right).
\]
Using
\[
 \norm{h_{I^c}}_2\le\norm h_2
\]
and the preceding bound on \(\norm{h_{I^c}}_\infty\), we get
\[
 \norm{\ip{h_{I^c}}R}_p
 \le
 C(\sqrt q+\sqrt\epsilon)\sqrt p\,\norm h_2.
\]

Combining the two parts gives
\begin{equation}
 \norm{\ip hR}_p
 \le
 C\left[
 (\sqrt\epsilon+\sqrt q)\sqrt p\,\norm h_2
 +
 \epsilon p\norm h_\infty
 \right].
\end{equation}
Since \(0<\epsilon\le1\),
\[
 \epsilon\le\sqrt\epsilon
 \le
 \sqrt\epsilon+\sqrt q.
\]
Together with the lower estimate in
\eqref{eq:laplace-moments}, this yields
\[
 \norm{\ip hR}_p
 \le
 C(\sqrt\epsilon+\sqrt q)
 \norm{\ip hZ}_p,
 \qquad p\ge2.
\]

For \(1\le p\le2\), monotonicity of \(L^p\)-norms and the case
\(p=2\) give
\[
 \norm{\ip hR}_p
 \le
 \norm{\ip hR}_2
 \le
 C(\sqrt\epsilon+\sqrt q)\norm h_2.
\]
By Lemma~\ref{lem:laplace-moments},
\[
 \norm{\ip hZ}_p\ge c\norm h_2,
 \qquad 1\le p\le2.
\]
Consequently,
\[
 \norm{\ip hR}_p
 \le
 C(\sqrt\epsilon+\sqrt q)
 \norm{\ip hZ}_p,
 \qquad h\in\R^d,\quad p\ge1.
\]

Since \(R\cx\eta\) and \(\eta\) is centered, \(R\) is centered.
The process comparison proved above therefore gives, for every finite
\(T\subset\R^d\),
\[
 S_R(T)
 \le
 C(\sqrt\epsilon+\sqrt q)\,S_Z(T).
\]

It remains to pass from this set comparison to convex order.
For \(N\ge1\), take independent copies
\[
 (U^{(1)},R^{(1)}),\ldots,(U^{(N)},R^{(N)})
\]
of the joint pair \((U,R)\), and write
\[
 U^{\oplus N}
 =
 (U^{(1)},\ldots,U^{(N)}),
 \qquad
 R^{\oplus N}
 =
 (R^{(1)},\ldots,R^{(N)}).
\]
By Lemma~\ref{lem:convex-pairing},
\[
 U^{\oplus N}\cx Z^{\oplus N},
 \qquad
 R^{\oplus N}\cx\eta^{\oplus N}.
\]
The coordinatewise bound
\[
 |R_i^{(j)}|
 \le
 \epsilon(1+|U_i^{(j)}|)
\]
also holds in every block. Thus the argument above applies in
dimension \(Nd\), with the same values of \(\epsilon\), \(q\), and
the same numerical constants. Hence
\[
 S_{R^{\oplus N}}(T)
 \le
 C(\sqrt\epsilon+\sqrt q)\,
 S_{Z^{\oplus N}}(T)
\]
for every \(N\) and every finite
\(T\subset(\R^d)^N\).

Applying Lemma~\ref{lem:type} now gives
\[
 R
 \cx
 C(\sqrt\epsilon+\sqrt q)\,Z,
\]
which is \eqref{eq:amplitude-result}.
\end{proof}

\subsection{Conditional tails and a martingale coupling}

We now construct the martingale coupling that will be used in the
truncation argument. The main point is to choose the coupling so that
the conditional distribution of each coordinate \(Z_i\), given
\(W=w\), still has a uniformly controlled tail. We first prove the
one-dimensional estimate needed for this purpose.

Recall that a nonnegative function \(h\) is log-concave if
\[
 h((1-\lambda)x+\lambda y)
 \ge
 h(x)^{1-\lambda}h(y)^\lambda,
 \qquad 0<\lambda<1.
\]
Zeros are allowed. In particular, the indicator of a convex set is
log-concave.

We use the following form of Pr\'ekopa's theorem
\cite[Theorem~6]{Prekopa}.

\begin{theorem}\label{thm:prekopa-input}
Suppose \(f:\R^{a+b}\to[0,\infty)\) is log-concave and
\(A\subset\R^b\) is convex. Then
\[
 g(x)=\int_A f(x,y)\,dy
\]
is log-concave.
\end{theorem}

In particular, marginals of log-concave densities are log-concave.
We will also use the theorem when \(A\) is a half-line or an
unbounded convex set.

\begin{lemma}\label{lem:conditional-tail}
Suppose \(S\) has density proportional to
\[
 e^{-\abs s}h(s),
\]
where \(h\) is nonnegative and log-concave, and let
\[
 m=\E\abs S<\infty.
\]
Then
\begin{equation}\label{eq:posterior-moments}
 m\le\abs{\E S}+2,\qquad
 \E(\abs S-t)_+\le Cm e^{-ct/m},
 \qquad t\ge0.
\end{equation}
If
\[
 \abs{\E S}\le\beta,\qquad \beta\ge1,
\]
then
\begin{equation}\label{eq:posterior-clip}
 \abs{\E[S-\clip_{K\beta}(S)]}
 \le
 \epsilon_K(1+\abs{\E S}),
 \qquad
 \epsilon_K=C_0e^{-c_0K},
 \quad K\ge1.
\end{equation}
\end{lemma}

\begin{proof}
We first prove
\[
 \E|S|\le|\E S|+2.
\]
If the support of \(S\) is contained in one half-line, then
\(\E|S|=|\E S|\), so there is nothing to prove.

Suppose therefore that the support meets both sides of zero. On the
interior of the support write
\[
 v=-\log h.
\]
The function \(v\) is convex. If
\[
 v'_+(0)\ge0,
\]
then \(v\) is nondecreasing on \([0,\infty)\). If \(f\) denotes the
density of \(S\), then for \(s,t\ge0\),
\[
 f(s+t)\le e^{-t}f(s).
\]
Indeed, on the interior of the support,
\[
 f(s+t)
 =
 f(s)e^{-t}e^{-[v(s+t)-v(s)]}
 \le
 e^{-t}f(s),
\]
and outside the support the same inequality is automatic.
Integrating over \(s\ge0\) gives
\[
 \Prob(S>t)
 \le
 e^{-t}\Prob(S>0),
\]
and hence
\[
 \E S_+
 =
 \int_0^\infty\Prob(S>t)\,dt
 \le1.
\]

If \(v'_+(0)<0\), convexity implies that \(v\) is nonincreasing on
the negative half-line. Applying the same argument to \(-S\) gives
\[
 \E S_-\le1.
\]
Thus in either case
\[
 \min\{\E S_+,\E S_-\}\le1.
\]
Since
\[
 \E|S|
 =
 |\E S|+2\min\{\E S_+,\E S_-\},
\]
the first estimate in \eqref{eq:posterior-moments} follows.

We next prove the tail estimate. Let \(f\) again denote the density
of \(S\), and define
\[
 F_+(u)=\Prob(S\ge u)
       =\int_0^\infty f(u+v)\,dv.
\]
The function
\[
 (u,v)\longmapsto f(u+v)\ind_{\{v\ge0\}}
\]
is log-concave. Therefore
Theorem~\ref{thm:prekopa-input} shows that \(F_+\) is log-concave.

Since \(m=\E|S|\), Markov's inequality gives
\[
 F_+(4m)\le\frac14,
 \qquad
 F_+(-4m)\ge\frac34.
\]
For \(u>4m\) such that \(F_+(u)>0\), concavity of
\(\log F_+\) yields
\begin{align*}
 \frac{\log F_+(u)-\log F_+(4m)}{u-4m}
 &\le
 \frac{\log F_+(4m)-\log F_+(-4m)}{8m}\\
 &\le
 -\frac{\log3}{8m}.
\end{align*}
Consequently,
\[
 \Prob(S\ge u)
 \le
 \frac14
 \exp\left[
 -\frac{\log3}{8m}(u-4m)
 \right],
 \qquad u\ge4m.
\]
If \(F_+(u)=0\), the same inequality is automatic.

Applying the same argument to \(-S\), and enlarging the numerical
constant to cover \(0\le u\le4m\), gives
\[
 \Prob(|S|>u)
 \le
 \sqrt3
 \exp\left[-\frac{\log3}{8m}u\right],
 \qquad u\ge0.
\]
Therefore
\begin{align*}
 \E(|S|-t)_+
 &=
 \int_t^\infty\Prob(|S|>u)\,du\\
 &\le
 \frac{8\sqrt3}{\log3}\,
 m\exp\left[-\frac{\log3}{8m}t\right].
\end{align*}
This proves the second estimate in
\eqref{eq:posterior-moments}.

Finally suppose that
\[
 |\E S|\le\beta,\qquad \beta\ge1.
\]
The first part gives
\[
 m\le|\E S|+2
 \le2(1+|\E S|)
\]
and also
\[
 m\le\beta+2\le3\beta.
\]
Moreover,
\[
 |S-\clip_{K\beta}(S)|
 =
 (|S|-K\beta)_+.
\]
Hence
\begin{align*}
 \abs{\E[S-\clip_{K\beta}(S)]}
 &\le
 \E(|S|-K\beta)_+\\
 &\le
 \frac{8\sqrt3}{\log3}\,
 m\exp\left[
 -\frac{\log3}{8m}K\beta
 \right]\\
 &\le
 \frac{16\sqrt3}{\log3}\,
 e^{-(\log3)K/24}
 (1+|\E S|).
\end{align*}
Thus \eqref{eq:posterior-clip} holds with, for example,
\[
 C_0=\frac{16\sqrt3}{\log3},
 \qquad
 c_0=\frac{\log3}{24}.
\]
\end{proof}

We next choose a martingale coupling for which
Lemma~\ref{lem:conditional-tail} applies to every conditional
coordinate distribution. The construction uses a finite-dimensional
convex optimization. We first replace each target mean \(w_j\) by
\(\alpha w_j\), \(0<\alpha<1\). This places the corresponding
mass-and-moment constraints in the interior of the feasible set.
After constructing the coupling for these smaller means, we let
\(\alpha\uparrow1\).

The following compactness fact will be used in this last step.

\begin{lemma}\label{lem:kernel-compact}
Let \(P\) be a Borel probability measure on \(\R^d\). Suppose
\[
 \kappa^{(N)}
 =
 (\kappa_1^{(N)},\ldots,\kappa_m^{(N)})
\]
are measurable functions satisfying
\[
 \kappa_j^{(N)}\ge0,
 \qquad
 \sum_{j=1}^m\kappa_j^{(N)}=1.
\]
Then some subsequence converges to a kernel
\(\kappa=(\kappa_1,\ldots,\kappa_m)\) in the sense that
\begin{equation}\label{eq:kernel-limit}
 \int g\kappa_j^{(N)}\,dP
 \longrightarrow
 \int g\kappa_j\,dP
 \qquad
 (g\in L^1(P),\ 1\le j\le m).
\end{equation}
The same subsequence works for every \(g\in L^1(P)\).
\end{lemma}

\begin{proof}
Choose a countable dense set
\[
 \{g_r:r\ge1\}\subset L^1(P).
\]
Since \(0\le\kappa_j^{(N)}\le1\), the sequence
\[
 \int g_r\kappa_j^{(N)}\,dP
\]
is bounded for every \(r\) and \(j\). A diagonal argument therefore
gives a subsequence along which all these integrals converge.

For arbitrary \(g,h\in L^1(P)\),
\[
 \left|
 \int(g-h)\kappa_j^{(N)}\,dP
 \right|
 \le
 \int|g-h|\,dP.
\]
The limits on the dense set consequently extend uniquely to a
bounded linear functional \(\Lambda_j\) on \(L^1(P)\). If \(g\ge0\),
then
\[
 0\le\Lambda_j(g)\le\int g\,dP,
 \qquad
 \sum_{j=1}^m\Lambda_j(g)=\int g\,dP.
\]

Define
\[
 Q_j(A)=\Lambda_j(\ind_A).
\]
The set function \(Q_j\) is finitely additive and satisfies
\[
 0\le Q_j(A)\le P(A).
\]
If \(A_n\downarrow\varnothing\), then
\[
 0\le Q_j(A_n)\le P(A_n)\longrightarrow0.
\]
Hence \(Q_j\) is countably additive and therefore is a finite measure,
with
\[
 0\le Q_j\le P.
\]
By the Radon--Nikodym theorem there is a measurable
\(\kappa_j\), \(0\le\kappa_j\le1\), such that
\[
 Q_j=\kappa_jP.
\]
The identity
\[
 \sum_jQ_j=P
\]
gives
\[
 \sum_j\kappa_j=1
 \qquad P\text{-almost surely}.
\]
Approximating an arbitrary \(g\in L^1(P)\) by simple functions now
gives \eqref{eq:kernel-limit}.
\end{proof}

\begin{proposition}\label{prop:selected}
Let \(W\) be finitely supported, \(W\cx Z\), and suppose
\[
 |W_i|\le\beta_i,
 \qquad
 1\le\beta_i<\infty.
\]
There is a coupling of \(W\) and \(Z\) satisfying
\[
 \E[Z\mid W]=W
\]
such that, simultaneously for every \(K\ge1\) and every support point
\(w\) of \(W\),
\begin{align}
 \E[\abs{Z_i}\mid W=w]
 &\le
 \abs{w_i}+2,
 \label{eq:selected-moment}\\
 r_{K,i}(w)
 :=
 \E[Z_i-\clip_{K\beta_i}(Z_i)\mid W=w],
 \qquad
 \abs{r_{K,i}(w)}
 &\le
 \epsilon_K(1+\abs{w_i}).
 \label{eq:selected-residual}
\end{align}
\end{proposition}

\begin{proof}
Write
\[
 \Prob(W=w_j)=\mu_j>0,
 \qquad 1\le j\le m.
\]
If \(m=1\), then \(W=0\) because \(W\cx Z\) and both vectors are
centered. In this case the independent coupling satisfies the
conclusion.

Assume \(m\ge2\). By Theorem~\ref{thm:strassen}, there is a
martingale coupling of \(W\) and \(Z\). Let
\[
 \kappa_j^*(z)
 =
 \Prob(W=w_j\mid Z=z)
\]
be its conditional label probabilities. They satisfy
\[
 \kappa_j^*\ge0,
 \qquad
 \sum_j\kappa_j^*=1,
\]
and
\[
 \E\kappa_j^*(Z)=\mu_j,
 \qquad
 \E[Z\kappa_j^*(Z)]
 =
 \mu_jw_j.
\]

Fix \(0<\alpha<1\), and define
\[
 \kappa_j^0(z)
 =
 \alpha\kappa_j^*(z)+(1-\alpha)\mu_j.
\]
Since \(\E Z=0\),
\[
 \E\kappa_j^0(Z)=\mu_j,
 \qquad
 \E[Z\kappa_j^0(Z)]
 =
 \alpha\mu_jw_j.
\]
Moreover,
\[
 \kappa_j^0(z)
 \ge
 (1-\alpha)\mu_j>0.
\]
Thus the mass constraints remain unchanged, while the required
conditional means are reduced from \(w_j\) to \(\alpha w_j\).

We regard
\[
 \kappa^0(z)
 =
 \bigl(\kappa_1^0(z),\ldots,\kappa_m^0(z)\bigr)
\]
as a probability kernel from \(z\in\R^d\) to the labels
\(1,\ldots,m\): \(\kappa_j^0(z)\) is the conditional probability
assigned to the point \(w_j\) when \(Z=z\). Since
\[
 \kappa_j^0(z)\ge(1-\alpha)\mu_j>0,
\]
every component stays uniformly away from zero. We may therefore make
small perturbations of the functions \(\kappa_j^0\) while preserving
their nonnegativity. These perturbations will be used below to vary
the mass-and-moment constraints.

Put
\[
 \zeta(z)=(1,z_1,\ldots,z_d).
\]

The constraints for the last label follow from those for the first
\(m-1\), so we work in
\[
 \mathsf E=(\R^{d+1})^{m-1}.
\]
Define
\[
 \mathcal C
 =
 \left\{
   \bigl(\E[\kappa_j(Z)\zeta(Z)]\bigr)_{j<m}:
   \kappa_j\ge0,\ \sum_j\kappa_j=1
 \right\},
\]
and set
\[
 t_j=(\mu_j,\alpha\mu_jw_j),
 \qquad
 t=(t_j)_{j<m}.
\]
The set \(\mathcal C\) is convex, because kernels may be mixed, and
it is compact by Lemma~\ref{lem:kernel-compact}, applied to the
integrable functions
\[
 1,Z_1,\ldots,Z_d.
\]

We next show that \(t\) is an interior point of \(\mathcal C\).
Set
\[
 g_0(z)=1,
 \qquad
 g_i(z)=\clip_1(z_i),
 \quad 1\le i\le d,
\]
and define
\[
 M_{ab}
 =
 \E[\zeta_a(Z)g_b(Z)],
 \qquad 0\le a,b\le d.
\]
Independence and symmetry of the coordinates of \(Z\) give
\[
 M
 =
 \operatorname{diag}(1,m_1,\ldots,m_d),
 \qquad
 m_i=\E[Z_i\clip_1(Z_i)]>0.
\]
In particular, \(M\) is invertible.

For \(j<m\), let
\[
 \theta_j=(\theta_{j,0},\ldots,\theta_{j,d}),
\]
and put
\[
 \kappa_j^\theta
 =
 \kappa_j^0+\sum_{b=0}^d\theta_{j,b}g_b,
 \qquad j<m,
\]
while
\[
 \kappa_m^\theta
 =
 \kappa_m^0
 -
 \sum_{j<m}\sum_{b=0}^d\theta_{j,b}g_b.
\]
For sufficiently small \(\theta\), these functions remain
nonnegative and still sum to one. Their constraint vectors are
\[
 \bigl(
   \E[\kappa_j^\theta(Z)\zeta(Z)]
 \bigr)_{j<m}
 =
 (t_j+M\theta_j)_{j<m}.
\]
Since \(M\) is invertible, these vectors contain a Euclidean
neighborhood of \(t\). Thus
\[
 t+aB_2\subset\mathcal C
\]
for some \(a>0\), where \(B_2\) is the Euclidean unit ball in
\(\mathsf E\).

We now solve the mass-and-moment constraints by a finite-dimensional
convex optimization. For
\[
 \lambda=(\lambda_1,\ldots,\lambda_{m-1})
 \in\mathsf E,
\]
define
\[
 F(\lambda)
 =
 \E\log\left(
 1+\sum_{j<m}
 e^{\ip{\lambda_j}{\zeta(Z)}}
 \right)
 -
 \sum_{j<m}\ip{\lambda_j}{t_j}.
\]
The support function of \(\mathcal C\) is
\[
 h_{\mathcal C}(\lambda)
 =
 \E\max\left\{
 0,
 \ip{\lambda_1}{\zeta(Z)},
 \ldots,
 \ip{\lambda_{m-1}}{\zeta(Z)}
 \right\}.
\]
Indeed, for each fixed \(z\), maximizing over
\[
 \kappa_j(z)\ge0,
 \qquad
 \sum_j\kappa_j(z)=1
\]
amounts to assigning all the probability to an index corresponding to
the largest of
\[
 0,
 \ip{\lambda_1}{\zeta(z)},
 \ldots,
 \ip{\lambda_{m-1}}{\zeta(z)}.
\]
Since
\[
 \log\left(1+\sum_{j<m}e^{u_j}\right)
 \ge
 \max\{0,u_1,\ldots,u_{m-1}\},
\]
the inclusion \(t+aB_2\subset\mathcal C\) gives
\[
 F(\lambda)
 \ge
 h_{\mathcal C}(\lambda)-\ip\lambda t
 \ge
 a\norm\lambda_2.
\]
Hence
\[
 F(\lambda)\longrightarrow\infty
 \qquad\text{as }\norm\lambda_2\to\infty,
\]
so \(F\) attains its minimum.

Write a minimizing multiplier as
\[
 \lambda_j=(b_j,a_j),
 \qquad
 b_j\in\R,\quad a_j\in\R^d,
\]
and set
\[
 (b_m,a_m)=(0,0).
\]
Define
\begin{equation}\label{eq:softmax}
 \kappa_j(z)
 =
 \frac{
 e^{\ip{a_j}z+b_j}
 }{
 \sum_{\ell=1}^m
 e^{\ip{a_\ell}z+b_\ell}
 }.
\end{equation}
The derivative of the integrand in \(F\) is bounded by a constant
multiple of
\[
 \norm{\zeta(Z)}_2\in L^1.
\]
We may therefore differentiate under the expectation. The first-order
conditions at the minimizer give
\[
 \E[\kappa_j(Z)\zeta(Z)]
 =
 t_j,
 \qquad j<m.
\]
Using
\[
 \sum_j\kappa_j=1
 \qquad\text{and}\qquad
 \E Z=0
\]
gives the same identities for \(j=m\). Hence
\[
 \E\kappa_j(Z)=\mu_j,
 \qquad
 \E[Z\kappa_j(Z)]
 =
 \alpha\mu_jw_j,
 \qquad 1\le j\le m.
\]

The important additional property is that each \(\kappa_j\) is
log-concave. Indeed,
\[
 \log\kappa_j(z)
 =
 \ip{a_j}z+b_j
 -
 \log\sum_\ell e^{\ip{a_\ell}z+b_\ell}.
\]
The last term is convex because
\[
 \nabla^2
 \log\sum_\ell e^{\ip{a_\ell}z+b_\ell}
 =
 \sum_\ell
 \kappa_\ell(z)
 (a_\ell-\bar a(z))
 (a_\ell-\bar a(z))^{\mathsf T}
 \succeq0,
\]
where
\[
 \bar a(z)=\sum_\ell\kappa_\ell(z)a_\ell.
\]
Thus \(\log\kappa_j\) is concave.

Under this kernel, the conditional density of \(Z_i\) given
\(W=w_j\) is
\[
 \frac{e^{-\abs s}}{2\mu_j}
 \int_{\R^{d-1}}
 \kappa_j(s,z_{-i})
 \prod_{\ell\ne i}
 \frac{e^{-\abs{z_\ell}}}{2}\,
 dz_{-i}.
\]
The integrand is log-concave in
\((s,z_{-i})\). Therefore
Theorem~\ref{thm:prekopa-input} shows that the integral is a
log-concave function of \(s\). The conditional density is consequently
of the form required in Lemma~\ref{lem:conditional-tail}.

Its conditional mean is
\[
 \E[Z_i\mid W=w_j]
 =
 \alpha w_{j,i}.
\]
Since
\[
 |\alpha w_{j,i}|
 \le
 |w_{j,i}|
 \le
 \beta_i,
\]
Lemma~\ref{lem:conditional-tail} gives, for every \(K\ge1\),
\begin{align*}
 \frac{\E[|Z_i|\kappa_j(Z)]}{\mu_j}
 &\le
 \alpha|w_{j,i}|+2,\\
 \left|
 \frac{
 \E[
 (Z_i-\clip_{K\beta_i}(Z_i))
 \kappa_j(Z)
 ]
 }{\mu_j}
 \right|
 &\le
 \epsilon_K(1+\alpha|w_{j,i}|).
\end{align*}

It remains to return from the means \(\alpha w_j\) to the desired
means \(w_j\). Choose
\[
 \alpha_N\uparrow1
\]
and let \(\kappa^{(N)}\) be the corresponding kernels.
By Lemma~\ref{lem:kernel-compact}, after passing to a subsequence
there is a kernel \(\kappa\) such that
\[
 \int g\kappa_j^{(N)}\,d\Law(Z)
 \longrightarrow
 \int g\kappa_j\,d\Law(Z)
\]
for every \(g\in L^1(\Law(Z))\) and every \(j\).

Taking successively
\[
 g=1,\quad
 g=Z_i,\quad
 g=|Z_i|,\quad
 g=Z_i-\clip_{K\beta_i}(Z_i)
\]
gives
\begin{align*}
 \E\kappa_j(Z)
 &=\mu_j,\\
 \E[Z_i\kappa_j(Z)]
 &=\mu_jw_{j,i},\\
 \E[|Z_i|\kappa_j(Z)]
 &\le
 \mu_j(|w_{j,i}|+2),\\
 \left|
 \E[
 (Z_i-\clip_{K\beta_i}(Z_i))
 \kappa_j(Z)
 ]
 \right|
 &\le
 \mu_j\epsilon_K(1+|w_{j,i}|).
\end{align*}
The convergence in
Lemma~\ref{lem:kernel-compact} holds for every integrable test
function, so the last inequality holds simultaneously for every
\(K\ge1\).

Finally define the joint law of \(W\) and \(Z\) by
\[
 \Prob(W=w_j,\,Z\in dz)
 =
 \kappa_j(z)\Law(Z)(dz).
\]
Then
\[
 \Prob(W=w_j)=\mu_j
\]
and
\[
 \E[Z\mid W=w_j]=w_j.
\]
Thus
\[
 \E[Z\mid W]=W.
\]
Dividing the last two estimates above by \(\mu_j\) gives
\eqref{eq:selected-moment} and
\eqref{eq:selected-residual}.
\end{proof}

\subsection{The uniform truncation theorem}

We now combine the preceding results to prove the uniform truncation
comparison. For a vector of caps
\[
 1\le\beta_i<\infty,
\]
let
\[
 \xi^{[\beta]}
 =
 (\eps_i\min(E_i,\beta_i))_{i=1}^d
\]
have independent coordinates.

For a finite index law \(\mu\), set
\begin{align}
 \mathcal K_Z(\mu)
 &=
 \{m:\supp\mu\to\R^d:m\push\mu\cx Z\},
 \label{eq:K}\\
 \delta_\beta(\mu)
 &=
 \sup_{\substack{m\in\mathcal K_Z(\mu)\\
                  |m_i(x)|\le\beta_i}}
 \E_\mu\ip X{m(X)}.
 \label{eq:delta}
\end{align}
The quantity \(\delta_\beta(\mu)\) is the largest pairing obtained
from a conditional mean of \(Z\) that also satisfies the prescribed
coordinate bounds.

More precisely, \(\mathcal K_Z(\mu)\) consists exactly of the maps
\[
 m(x)=\E[Z\mid X=x]
\]
that arise from couplings of \(X\sim\mu\) and \(Z\).
Indeed, any such conditional mean satisfies
\[
 m(X)\cx Z
\]
by Jensen inequality. Conversely, if
\[
 m(X)\cx Z,
\]
Theorem~\ref{thm:strassen} gives a coupling of \(m(X)\) and \(Z\)
such that
\[
 \E[Z\mid m(X)]=m(X).
\]
Conditionally on \(m(X)\), one may then generate \(X\) according to
its original conditional distribution given \(m(X)\). This gives a
coupling of \(X\sim\mu\) and \(Z\) with
\[
 \E[Z\mid X]=m(X).
\]

Since \(\mu\) has finite support, a map
\(m:\supp\mu\to\R^d\) may be identified with a vector in a
finite-dimensional Euclidean space. The set
\(\mathcal K_Z(\mu)\) is bounded because
\[
 \sum_x\mu(x)\norm{m(x)}_2
 =
 \E\norm{m(X)}_2
 \le
 \E\norm Z_2.
\]
It is closed: if \(m^{(n)}(x)\to m(x)\) at every support point and
\(m^{(n)}\push\mu\cx Z\), then for every real-valued convex function
\(f\),
\[
 \sum_x\mu(x)f(m^{(n)}(x))
 \le
 \E f(Z),
\]
and passage to the limit gives the same inequality for \(m\).
It is also convex, since for \(m^0,m^1\in\mathcal K_Z(\mu)\) and
\(0\le\lambda\le1\),
\[
 f\bigl((1-\lambda)m^0(x)+\lambda m^1(x)\bigr)
 \le
 (1-\lambda)f(m^0(x))+\lambda f(m^1(x)).
\]
Thus \(\mathcal K_Z(\mu)\) is compact and convex. The coordinate
constraints in \eqref{eq:delta} are closed, so the supremum there is
finite and attained.

\begin{theorem}\label{thm:capping}
There is a numerical constant \(C\) such that
\begin{equation}\label{eq:uniform-capping}
 W\cx Z,\qquad |W_i|\le\beta_i
 \quad\Longrightarrow\quad
 W\cx C\xi^{[\beta]}.
\end{equation}
No independence among the coordinates of \(W\) is assumed.

Equivalently, for every finite index law \(\mu\),
\begin{equation}\label{eq:capped-pairing}
 B_{\xi^{[\beta]}}(\mu)
 \le
 \delta_\beta(\mu)
 \le
 C B_{\xi^{[\beta]}}(\mu).
\end{equation}
\end{theorem}

\begin{proof}
We first prove the upper bound in
\eqref{eq:capped-pairing}. Fix an admissible map \(m\) in
\eqref{eq:delta}, let
\[
 W=m(X),
 \qquad X\sim\mu,
\]
and choose for \(W\) the coupling with \(Z\) constructed in
Proposition~\ref{prop:selected}. For \(K\ge1\), define
\[
 R_K
 =
 \E[Z-\clip_{K\beta}Z\mid W]
 =
 r_K(W).
\]

We first record the two estimates on \(R_K\) that will allow us to use
Theorem~\ref{thm:amplitude}. Since \(R_K\) is the conditional
expectation of \(Z-\clip_{K\beta}Z\) given \(W\), conditional Jensen
gives
\[
 R_K\cx Z-\clip_{K\beta}Z.
\]
The coordinates of the vector on the right remain independent.
Moreover, for \(s\ge0\),
\[
 \Prob\bigl((|Z_i|-K\beta_i)_+>s\bigr)
 =
 e^{-K\beta_i}e^{-s}.
\]
By the memoryless property of the exponential distribution,
\[
 Z_i-\clip_{K\beta_i}(Z_i)
 \stackrel{\rm d}=
 \delta_i Z_i',
\]
where \(Z_i'\) is a standard symmetric exponential variable and
\(\delta_i\) is a \(\{0,1\}\)-valued variable with
\[
 \Prob(\delta_i=1)=e^{-K\beta_i}\le e^{-K}.
\]
These variables may be chosen independently across coordinates.
Hence
\[
 R_K\cx\eta,
 \qquad
 \eta_i=\delta_iZ_i'.
\]

On the other hand, Proposition~\ref{prop:selected} gives the
pointwise estimate
\[
 |R_{K,i}|
 \le
 \epsilon_K(1+|W_i|).
\]
Thus Theorem~\ref{thm:amplitude} applies, once \(K\) is large enough
that \(\epsilon_K\le1\). Choose a numerical constant so that
\begin{equation}\label{eq:theta}
 \theta_K
 =
 C_1\bigl(\sqrt{\epsilon_K}+e^{-K/2}\bigr)
 \le
 C_2e^{-c_2K},
 \qquad
 \theta_K\ge2\epsilon_K.
\end{equation}
Then
\[
 R_K\cx\theta_K Z.
\]
Since \(|W_i|\le\beta_i\) and \(\beta_i\ge1\),
\[
 |R_{K,i}|
 \le
 \epsilon_K(1+|W_i|)
 \le
 2\epsilon_K\beta_i
 \le
 \theta_K\beta_i.
\]

This is the key point: after division by \(\theta_K\), the residual
satisfies the same two constraints that define
\(\delta_\beta(\mu)\). Indeed, define
\[
 \widetilde m(x)
 =
 \frac{r_K(m(x))}{\theta_K}.
\]
Then
\[
 \widetilde m\push\mu
 =
 \Law(R_K/\theta_K)
 \cx Z,
\]
and
\[
 |\widetilde m_i(x)|\le\beta_i.
\]
Therefore \(\widetilde m\) is admissible in \eqref{eq:delta}, and
\[
 \E\ip X{R_K}
 =
 \theta_K\E\ip X{\widetilde m(X)}
 \le
 \theta_K\delta_\beta(\mu).
\]

It remains to estimate the truncated part. The coupling in
Proposition~\ref{prop:selected} was constructed conditionally on
\(W=m(X)\). We may use this same conditional law of \(Z\) while
keeping the original pair \((X,W)\). Then
\[
 \E[Z\mid X,W]=W=m(X).
\]
Since \(W\) is a function of \(X\),
\[
 \E[Z\mid X]=m(X),
\]
and similarly
\[
 \E[Z-\clip_{K\beta}Z\mid X]=R_K.
\]
Consequently,
\begin{align*}
 \E\ip X{m(X)}
 &=
 \E\ip XZ\\
 &=
 \E\ip X{\clip_{K\beta}Z}
 +
 \E\ip X{R_K}\\
 &\le
 B_{\xi^{[K\beta]}}(\mu)
 +
 \theta_K\delta_\beta(\mu).
\end{align*}

For \(K\ge1\),
\[
 \min(E_i,K\beta_i)
 \le
 K\min(E_i,\beta_i).
\]
For each coordinate, Lemma~\ref{lem:convex-pairing}(iv) gives
\[
 \eps_i\min(E_i,K\beta_i)
 \cx
 K\eps_i\min(E_i,\beta_i).
\]
Applying part~(ii) of the same lemma to the independent coordinates
gives
\[
 \xi^{[K\beta]}
 \cx
 K\xi^{[\beta]}.
\]
Hence
\[
 B_{\xi^{[K\beta]}}(\mu)
 \le
 K B_{\xi^{[\beta]}}(\mu).
\]
Choose one numerical \(K_0\) large enough that
\[
 \theta_{K_0}\le\frac12.
\]
Taking the supremum over all admissible \(m\) gives
\[
 \delta_\beta(\mu)
 \le
 K_0B_{\xi^{[\beta]}}(\mu)
 +
 \theta_{K_0}\delta_\beta(\mu),
\]
and therefore
\[
 \delta_\beta(\mu)
 \le
 2K_0 B_{\xi^{[\beta]}}(\mu).
\]

We now prove the reverse inequality in
\eqref{eq:capped-pairing}. First note that
\[
 \xi^{[\beta]}\cx Z.
\]
Indeed, coordinatewise,
\[
 |\xi_i^{[\beta]}|
 =
 \min(E_i,\beta_i)
 \le
 E_i
 =
 |Z_i|,
\]
and both variables are symmetric. Thus
Lemma~\ref{lem:convex-pairing}(iv) gives the coordinatewise convex-order
comparison, and part~(ii) gives
\[
 \xi^{[\beta]}\cx Z.
\]

Take any coupling of \(X\sim\mu\) and \(\xi^{[\beta]}\), and define
\[
 m(X)
 =
 \E[\xi^{[\beta]}\mid X].
\]
Jensen's inequality gives
\[
 m(X)\cx\xi^{[\beta]}\cx Z,
\]
while
\[
 |m_i(X)|
 \le
 \E[|\xi_i^{[\beta]}|\mid X]
 \le
 \beta_i.
\]
Thus \(m\) is admissible in \eqref{eq:delta}, and
\[
 \E\ip X{\xi^{[\beta]}}
 =
 \E\ip X{m(X)}
 \le
 \delta_\beta(\mu).
\]
Taking the supremum over all couplings proves
\[
 B_{\xi^{[\beta]}}(\mu)
 \le
 \delta_\beta(\mu).
\]
This establishes \eqref{eq:capped-pairing}.

Finally, let \(W\) satisfy the assumptions of
\eqref{eq:uniform-capping}. To prove the convex-order conclusion,
fix any finite index law \(\mu\) and any coupling of \(X\sim\mu\)
with \(W\). Define
\[
 m(X)=\E[W\mid X].
\]
Then
\[
 m(X)\cx W\cx Z
\]
by conditional Jensen and transitivity of convex order, and
\[
 |m_i(X)|
 \le
 \E[|W_i|\mid X]
 \le
 \beta_i.
\]
Hence \(m\) is admissible in \eqref{eq:delta}. Moreover,
\[
 \E\ip XW
 =
 \E\ip X{m(X)}
 \le
 \delta_\beta(\mu)
 \le
 C B_{\xi^{[\beta]}}(\mu).
\]
Taking the supremum over couplings of \(X\) and \(W\) yields
\[
 B_W(\mu)
 \le
 C B_{\xi^{[\beta]}}(\mu)
 =
 B_{C\xi^{[\beta]}}(\mu).
\]
Since this holds for every finite index law \(\mu\),
Lemma~\ref{lem:convex-pairing} gives
\[
 W\cx C\xi^{[\beta]}.
\]
This proves \eqref{eq:uniform-capping}.
\end{proof}

We next express the coordinate constraint in
\(\delta_\beta(\mu)\) by an \(\ell^1\) penalty. The resulting formula
will be used in the comparison with the geometric cost in the next
section.

We use Sion's minimax theorem in the following form
\cite[Corollary~3.5, p.~174]{Sion}. If \(M\) is compact and convex,
\(N\) is convex, and \(f:M\times N\to\R\) is continuous, concave in
its first variable, and convex in its second, then
\begin{equation}\label{eq:Sion-input}
 \sup_{m\in M}\inf_{a\in N}f(m,a)
 =
 \inf_{a\in N}\sup_{m\in M}f(m,a).
\end{equation}
The hypotheses above are stronger than the semicontinuity and
quasi-convexity assumptions in Sion's theorem, but they are exactly
the ones needed here.

\begin{proposition}\label{prop:decomposition-dual}
For every finite index law \(\mu\),
\begin{equation}\label{eq:decomposition-dual}
 \delta_\beta(\mu)
 =
 \inf_{a:\supp\mu\to\R^d}
 \left\{
 \E\sum_i\beta_i|a_i(X)|
 +
 B_Z((\operatorname{id}-a)\push\mu)
 \right\}.
\end{equation}
\end{proposition}

\begin{proof}
Put
\[
 T=\supp\mu.
\]
For
\[
 m\in\mathcal K_Z(\mu),
 \qquad
 a:T\to\R^d,
\]
define
\[
 f(m,a)
 =
 \sum_{x\in T}\mu(x)
 \left[
 \ip x{m(x)}
 +
 \sum_i
 \bigl(
 \beta_i|a_i(x)|
 -
 a_i(x)m_i(x)
 \bigr)
 \right].
\]

We first identify the supremum over \(m\). By
\eqref{eq:fiber} and the characterization of
\(\mathcal K_Z(\mu)\) above,
\[
 B_Z((\operatorname{id}-a)\push\mu)
 =
 \sup_{m\in\mathcal K_Z(\mu)}
 \E\ip{X-a(X)}{m(X)}.
\]
Therefore
\begin{align*}
 &\inf_a
 \left\{
 \E\sum_i\beta_i|a_i(X)|
 +
 B_Z((\operatorname{id}-a)\push\mu)
 \right\}\\
 &\qquad=
 \inf_a\sup_{m\in\mathcal K_Z(\mu)}f(m,a).
\end{align*}

Now take
\[
 M=\mathcal K_Z(\mu),
 \qquad
 N=(\R^d)^T.
\]
The set \(M\) is compact and convex, while \(N\) is convex.
The function \(f\) is continuous and affine in \(m\), and convex in
\(a\). Hence Sion's theorem gives
\begin{align*}
 \inf_a\sup_{m\in\mathcal K_Z(\mu)}f(m,a)
 &=
 \sup_{m\in\mathcal K_Z(\mu)}\inf_a f(m,a)\\
 &=
 \sup_{m\in\mathcal K_Z(\mu)}
 \left[
 \E\ip X{m(X)}
 +
 \sum_{x\in T}\sum_i
 \mu(x)
 \inf_{u\in\R}
 \bigl\{
 \beta_i|u|-u\,m_i(x)
 \bigr\}
 \right].
\end{align*}

For every \(i\) and every \(x\in T\),
\[
 \inf_{u\in\R}
 \bigl\{
 \beta_i|u|-u\,m_i(x)
 \bigr\}
 =
 \begin{cases}
 0,& |m_i(x)|\le\beta_i,\\
 -\infty,& |m_i(x)|>\beta_i.
 \end{cases}
\]
Indeed, if \(|m_i(x)|\le\beta_i\), then
\[
 \beta_i|u|-u\,m_i(x)\ge0
\]
for every \(u\), with equality at \(u=0\). If
\(m_i(x)>\beta_i\), let \(u\to+\infty\); if
\(m_i(x)<-\beta_i\), let \(u\to-\infty\).

Thus the inner infimum is finite exactly for those maps \(m\) that
satisfy
\[
 |m_i(x)|\le\beta_i
 \qquad\text{for every }x,i.
\]
For such maps its value is zero. The remaining supremum is therefore
exactly the quantity in \eqref{eq:delta}, proving
\eqref{eq:decomposition-dual}.
\end{proof}

\section{The cost of truncated exponentials}\label{sec:geometry}

Section~\ref{sec:convex-order} gave the probabilistic comparison for
truncated symmetric exponentials. We now identify the deterministic
cost associated with these variables in the rate--distortion
functional. Recall
\[
 \xi^{[\beta]}
 =
 \bigl(\varepsilon_i\min(E_i,\beta_i)\bigr)_{i=1}^d,
 \qquad
 1\le \beta_i\le\infty.
\]
For the corresponding coordinate costs \(\psi_{\beta_i}\),
Theorem~\ref{thm:capped-geometry} will show, uniformly over the caps,
that
\[
 \mathcal R^{\psi_\beta}(\mu)
 \asymp
 B_{\xi^{[\beta]}}(\mu)
\]
for every finite index law \(\mu\). This is the bridge between the
convex-order comparison of the preceding section and the reduction to
general log-concave tails in Section~\ref{sec:general-tails}.

The proof proceeds in several steps. We first prove the area identity
\eqref{eq:area-identity}, which identifies the geometric meaning of
the intrinsic costs that later enter the rate--distortion functional.
We then collect the information-theoretic and coding facts needed to
work with randomized approximations while preserving the prescribed
marginal law. Next we treat the uncapped exponential process, using
chaining and the replication principle to obtain the required estimate
for every fixed index law. Finally, for finite coordinatewise caps, we
construct a decomposition by clipping the differences between
successive dyadic approximations. Mixed finite and infinite caps then
follow by approximation, and the case \(\beta_i=1\) recovers the
Bernoulli cost.

\subsection{An area identity}

For \(p>0\), define
\[
 C_p=\left\{z\in\R^n:\sum_iV_i(z_i)\le p\right\},
 \qquad
 \rho_p(h)=\sup_{z\in C_p}\ip hz.
\]
Thus \(\rho_p(h)\) is the support function of \(C_p\) in the
direction \(h\). The following identity shows that the same quantity
is recovered by integrating the scale costs \(D_t\).

\begin{lemma}\label{lem:area-identity}
For every \(h\in\R^n\) and \(p>0\),
\begin{equation}\label{eq:area-identity}
 \rho_p(h)=\int_0^\infty\min\{p,D_t(0,h)\}\,dt.
\end{equation}
Moreover,
\[
 C_1=\{z\in\R^n:\norm z_2\le1\},
\]
and hence \(\rho_1(h)=\norm h_2\).
\end{lemma}

\begin{proof}
We first establish a one-dimensional identity. Fix one coordinate and
write
\[
 V=V_i,\qquad b=b_i,\qquad s=V'_+,\qquad f=V^*,\qquad \Psi=\Psi_i.
\]
Since \(V(0)=0\) and \(V\) is convex,
\[
 V(u)=\int_0^u s(v)\,dv,
 \qquad 0\le u<b,
\]
with the endpoint value obtained by monotone convergence when
\(b<\infty\). Therefore
\[
 ru-V(u)=\int_0^u(r-s(v))\,dv.
\]
Because \(s\) is nondecreasing, the integrand can change sign only
from positive to nonpositive. Maximizing over \(0\le u\le b\) gives
\begin{equation}\label{eq:conjugate-integral}
 f(r)
 =
 \sup_{0\le u\le b}\{ru-V(u)\}
 =
 \int_0^b(r-s(u))_+\,du,
 \qquad r\ge0.
\end{equation}
The same formula remains valid when \(s\) has flat intervals or when
the endpoint is not attained.

For \(a\ge0\) and \(t>0\), Tonelli's theorem and the definition of
\(\Psi\) give
\begin{align}
 \int_t^\infty\Psi(a/v)\,dv
 &=
 \int_0^b s(u)
 \int_t^\infty\ind_{\{s(u)<a/v\}}\,dv\,du \notag\\
 &=
 \int_0^b(a-ts(u))_+\,du
 =
 t f(a/t).
 \label{eq:perspective-area}
\end{align}
Applying this identity coordinatewise, set
\[
 H(t)=t\sum_i f_i(\abs{h_i}/t).
\]
Then
\[
 H(t)=\int_t^\infty D_v(0,h)\,dv.
\]

We next relate \(H\) to the support function. Put
\[
 \mathcal V(z)=\sum_iV_i(z_i).
\]
For every \(t>0\), the definition of the convex conjugates gives
\[
 \ip hz
 \le
 t\mathcal V(z)+t\sum_i f_i(\abs{h_i}/t).
\]
Hence, for \(z\in C_p\),
\[
 \ip hz\le pt+H(t),
\]
and therefore
\[
 \rho_p(h)\le pt+H(t),
 \qquad t>0.
\]

For the reverse inequality, the map \(q\mapsto\rho_q(h)\) is finite,
nondecreasing, and concave, with \(\rho_0(h)=0\). Fix \(p>0\) and
choose a supporting line at \(p\): for some \(\lambda\ge0\),
\[
 \rho_q(h)\le\rho_p(h)+\lambda(q-p),
 \qquad q\ge0.
\]
If \(\mathcal V(z)<\infty\), then
\[
 \ip hz\le\rho_{\mathcal V(z)}(h),
\]
so
\[
 \ip hz-\lambda\mathcal V(z)
 \le
 \rho_p(h)-\lambda p.
\]
If \(\lambda>0\), taking the supremum over \(z\) yields
\[
 H(\lambda)
 =
 \sup_z\{\ip hz-\lambda\mathcal V(z)\}
 \le
 \rho_p(h)-\lambda p.
\]
Together with the preceding upper bound, this gives equality at
\(t=\lambda\). If \(\lambda=0\), the supporting-line inequality and
monotonicity imply
\[
 \sup_{\mathcal V(z)<\infty}\ip hz=\rho_p(h).
\]
Therefore
\[
 H(t)=\sup_z\{\ip hz-t\mathcal V(z)\}\le\rho_p(h),
\]
and letting \(t\downarrow0\) in
\(\rho_p(h)\le pt+H(t)\) again gives equality. Thus
\begin{equation}\label{eq:support-perspective}
 \rho_p(h)
 =
 \inf_{t>0}\left\{pt+H(t)\right\}.
\end{equation}

It remains to rewrite this infimum. Put
\[
 d(v)=D_v(0,h),
 \qquad
 t_*:=\sup\{v>0:d(v)>p\},
 \qquad \sup\varnothing=0.
\]
The function \(d\) is nonincreasing. For all sufficiently large
\(v\), \eqref{eq:core} gives
\[
 d(v)=\frac{\norm h_2^2}{4v^2},
\]
so \(d(v)\to0\) and its tail is integrable. For every \(t>0\),
\[
 \int_0^\infty\min\{p,d(v)\}\,dv
 \le
 pt+\int_t^\infty d(v)\,dv.
\]
If \(t_*>0\), monotonicity gives
\[
 d(v)\ge p\quad(0<v<t_*),
 \qquad
 d(v)\le p\quad(v>t_*),
\]
up to the value at the single point \(t_*\). Hence
\[
 \int_0^\infty\min\{p,d(v)\}\,dv
 =
 pt_*+\int_{t_*}^\infty d(v)\,dv.
\]
If \(t_*=0\), the same equality follows by letting \(t\downarrow0\).
No continuity of \(d\) is required. Combining this identity with
\eqref{eq:support-perspective} proves \eqref{eq:area-identity}.

Finally, suppose \(z\in C_1\). Since \(V_i(u)>1\) whenever
\(|u|>1\), we must have \(|z_i|\le1\) for every \(i\). On this
range \(V_i(z_i)=z_i^2\), and hence
\[
 C_1
 =
 \left\{z\in\R^n:\sum_i z_i^2\le1\right\}.
\]
Thus \(C_1\) is the Euclidean unit ball and
\(\rho_1(h)=\norm h_2\).
\end{proof}

\subsection{Information-theoretic preliminaries}\label{sec:information}

We collect the elementary information facts used in the proof.
They are standard; see \cite[Chapters 2 and 10]{CT}. Stating them
here also fixes the meaning of a randomized approximation, which
we sometimes call a code.

For probability measures $P,Q$, their relative entropy is
\[
 \KL(P\Vert Q)=\int\log\frac{dP}{dQ}\,dP
\]
when $P$ is absolutely continuous with respect to $Q$, and is
$+\infty$ otherwise. It is nonnegative, by Jensen's inequality.
For a finite random variable $U$, its entropy is
$\Ent(U)=-\sum_u\Prob(U=u)\log\Prob(U=u)$.
Mutual information and conditional mutual information are
\begin{align*}
 I(A;B)&=\KL(\Law(A,B)\Vert\Law(A)\otimes\Law(B)),\\
 I(A;B\mid C)&=\E_C\KL\bigl(\Law(A,B\mid C)
             \Vert\Law(A\mid C)\otimes\Law(B\mid C)\bigr).
\end{align*}
These definitions also apply when $A$ has a density. In the
applications below, the other variable is finite, so its entropy
bounds the mutual information. In particular,
\[
 I(A;U)=\Ent(U)-\Ent(U\mid A),\qquad 0\le I(A;U)\le\Ent(U),
\]
where $\Ent(U\mid A)$ is the average entropy of the conditional
law of $U$ given $A$.

The chain rule follows by averaging the identity
\[
 \log\frac{p(a,b,c)}{p(a)p(b,c)}
 =\log\frac{p(a,b)}{p(a)p(b)}
                  +\log\frac{p(c\mid a,b)}{p(c\mid b)}:
 \qquad I(A;B,C)=I(A;B)+I(A;C\mid B).
\]
Here $p$ denotes the relevant probabilities or densities, and
only positive-probability terms enter the expectation.
If $A$ and $C$ are conditionally independent given $B$, this gives
\[
 I(A;C)\le I(A;B,C)=I(A;B).
\]
This is the data-processing inequality. If finite $U,V$ are
conditionally independent given $X$, the entropy formula gives
\begin{align*}
 I(X;U,V)
 &=\Ent(U,V)-\Ent(U,V\mid X)\\
 &\le\Ent(U)+\Ent(V)-\Ent(U\mid X)-\Ent(V\mid X)\\
 &=I(X;U)+I(X;V).
\end{align*}

For a probability measure $\mu$, a self-coupling of $\mu$ is a joint
law of $(X,X')$ such that both marginals are equal to $\mu$. We denote
the set of all self-couplings by
\[
\Pi(\mu,\mu).
\]
For a self-coupling $P\in\Pi(\mu,\mu)$, we write
\[
I(P)=I(X;X')
\]
for the mutual information of the corresponding pair.

We will also use the relative-entropy variational inequality
\begin{equation}\label{eq:entropy-variational}
 \E_P f\le\KL(P\Vert Q)+\log\E_Q e^f.
\end{equation}

Finally, let a finite observation $O$ be drawn from $X\sim\mu$
by a kernel $K_x(o)=\Prob(O=o\mid X=x)$. Put
$\nu(o)=\sum_x\mu(x)K_x(o)$. Given $O=o$ with $\nu(o)>0$,
draw $X'$ from
\[
 L_o(x')=\frac{\mu(x')K_{x'}(o)}{\nu(o)}.
\]
The joint law is
\begin{equation}\label{eq:posterior-kernel}
 \Prob(X=x,O=o,X'=x')
       =\frac{\mu(x)K_x(o)\mu(x')K_{x'}(o)}{\nu(o)}.
\end{equation}
Labels with $\nu(o)=0$ can be omitted. Summing over $x$ or
$x'$ gives
\begin{equation}\label{eq:posterior-copy}
 X'\sim\mu,\qquad \Law(X',O)=\Law(X,O),\qquad
 I(X;X')\le I(X;O).
\end{equation}
The last inequality is data processing: \eqref{eq:posterior-kernel}
makes $X,X'$ conditionally independent given $O$.
This construction is also used in \cite[equation~(26)]{Liu}.
We call $X'$ a conditional copy given $O$. It restores the
second marginal in \eqref{eq:RD} without changing the joint law
of the original index and its observation.

\subsection{Codes and limits}
We will occasionally work with a general coordinate cost $d=(d_i)_i$, not necessarily coming from a truncated exponential. Here, $d_i:[0,\infty)\to[0,\infty]$ is a nondecreasing function with
$d_i(0)=0$. It need not be a metric. 
We set $d_t(x,y)=\sum_i d_i(\abs{x_i-y_i}/t).$ These costs satisfy the scale triangle inequality
\begin{equation}\label{eq:scale-triangle}
 d_{2t}(x,y)\le d_t(x,z)+d_t(z,y).
\end{equation}
Indeed,
\[
 \frac{|x_i-y_i|}{2t}\le
 \max\left\{\frac{|x_i-z_i|}{t},\frac{|z_i-y_i|}{t}\right\}.
\]
Apply the nondecreasing function $d_i$ and bound the maximum
by the sum, then sum over $i$. 

Motivated in part by recent fixed-law and information-theoretic formulations for Bernoulli processes \cite{Liu,LZ}, we use the following self-coupling version of the rate--distortion functional. For a finite index law \(\mu\), define 
\[ \RD_\mu^d(t) = \inf_{P\in\Pi(\mu,\mu)} \left\{ I(P)+\E_P d_t(X,X') \right\}, \]
where \(\Pi(\mu,\mu)\) denotes the set of couplings whose two marginals are both \(\mu\). We then set
\[
 \cR^d(\mu)=\int_0^\infty \RD_\mu^d(t)\,dt.
\]
For a finite set $T$, we also write
\[
 \cR^d(T)=\sup_{\mu\in\cP(T)}\cR^d(\mu).
\]

\begin{lemma}\label{lem:coding}
Let \(d=(d_i)\) be a family of nondecreasing costs with \(d_i(0)=0\).
\begin{enumerate}[label=\textup{(\roman*)}]
\item
For \(a>0\), define \(d_i^{(a)}(r)=d_i(ar)\). Then
\[
 \cR^{d^{(a)}}(\mu)=a\,\cR^d(\mu).
\]

\item
For \(c\ge1\),
\[
 \cR^d(\mu)\le\cR^{cd}(\mu)\le c\,\cR^d(\mu),
\]
where \(cd=(cd_i)\).

\item
For two cost families \(d=(d_i)\) and \(e=(e_i)\),
\begin{equation}\label{eq:add-cost}
 \cR^{d+e}(\mu)
 \le
 4\bigl\{\cR^d(\mu)+\cR^e(\mu)\bigr\}.
\end{equation}

\item
Suppose that
\[
 x=a(x)+v(x),
 \qquad x\in\supp\mu,
\]
for maps \(a,v:\supp\mu\to\R^n\). Then
\begin{equation}\label{eq:split-maps}
 \cR^d(\mu)
 \le
 8\bigl\{\cR^d(a\push\mu)+\cR^d(v\push\mu)\bigr\}.
\end{equation}

\item
Let \(\phi:\supp\mu\to S\) be injective, and let
\(\widetilde d_t\) be a distortion on \(S\) such that
\[
 \widetilde d_t(\phi(x),\phi(y))=d_t(x,y)
 \qquad(x,y\in\supp\mu).
\]
Then
\[
 \RD^{\widetilde d}_{\phi\push\mu}(t)=\RD^d_\mu(t),
 \qquad t>0,
\]
and consequently
\[
 \cR^{\widetilde d}(\phi\push\mu)=\cR^d(\mu).
\]
\end{enumerate}
\end{lemma}

\begin{proof}
For part~(i),
\[
 \RD_\mu^{d^{(a)}}(t)=\RD_\mu^d(t/a),
\]
so the change of variables \(t=as\) gives the assertion.
For part~(ii), for every self-coupling \(P\),
\[
 I(P)+\E_P d_t
 \le
 I(P)+c\E_P d_t
 \le
 c\bigl(I(P)+\E_P d_t\bigr).
\]
Taking infima and then integrating in \(t\) proves the claim.

For \eqref{eq:add-cost}, fix \(t>0\) and choose a self-coupling
\((X,U)\) within \(\varepsilon\) of \(\RD^d_\mu(t)\), and a
self-coupling \((X,V)\) within \(\varepsilon\) of
\(\RD^e_\mu(t)\). Generate \(U\) and \(V\) conditionally
independently given \(X\) using the corresponding kernels. Then
\[
 I(X;U,V)\le I(X;U)+I(X;V).
\]
Let \(X'\) be drawn from the conditional law of \(X\) given
\(O=(U,V)\), independently of \(X\) conditional on \(O\). By
\eqref{eq:posterior-copy}, \(X'\sim\mu\) and
\[
 I(X;X')\le I(X;U,V).
\]
Moreover, \((X,U)\) and \((X',U)\) have the same law. Hence
\[
 \E d_{2t}(X,X')
 \le
 \E d_t(X,U)+\E d_t(U,X')
 =
 2\E d_t(X,U).
\]
The same argument with \(V\) gives
\[
 \E e_{2t}(X,X')\le2\E e_t(X,V).
\]
Letting \(\varepsilon\downarrow0\),
\[
 \RD^{d+e}_\mu(2t)
 \le
 2\RD^d_\mu(t)+2\RD^e_\mu(t).
\]
Integration proves \eqref{eq:add-cost}.

For \eqref{eq:split-maps}, let \(X\sim\mu\) and fix \(t>0\).
Choose near-minimizing self-couplings for the laws of \(a(X)\) and
\(v(X)\) at scale \(t\). Using the corresponding conditional
kernels, generate \(U\) from \(a(X)\) and \(V\) from \(v(X)\),
conditionally independently given \(X\). Thus
\[
 U\stackrel{d}=a(X),
 \qquad
 V\stackrel{d}=v(X),
\]
and the conditional law of \(U\) given \(X\) depends only on
\(a(X)\), while that of \(V\) depends only on \(v(X)\). Consequently,
\[
 I(X;U)=I(a(X);U),
 \qquad
 I(X;V)=I(v(X);V),
\]
and
\[
 I(X;U,V)
 \le
 I(a(X);U)+I(v(X);V).
\]
Now let \(X'\) be drawn from the conditional law of \(X\) given
\((U,V)\), independently of \(X\) conditional on \((U,V)\). Then
\(X'\sim\mu\) and
\[
 X\longrightarrow(U,V)\longrightarrow X'
\]
is a Markov chain. Hence
\[
 I(X;X')
 \le
 I(a(X);U)+I(v(X);V).
\]
Since
\[
 X-X'
 =
 \bigl(a(X)-a(X')\bigr)
 +
 \bigl(v(X)-v(X')\bigr),
\]
two applications of \eqref{eq:scale-triangle} give
\[
 d_{4t}(X,X')
 \le
 d_{2t}\bigl(a(X),a(X')\bigr)
 +
 d_{2t}\bigl(v(X),v(X')\bigr).
\]
Applying \eqref{eq:scale-triangle} once more to each term,
\begin{align*}
 d_{2t}\bigl(a(X),a(X')\bigr)
 &\le d_t(a(X),U)+d_t(U,a(X')),\\
 d_{2t}\bigl(v(X),v(X')\bigr)
 &\le d_t(v(X),V)+d_t(V,v(X')).
\end{align*}
Conditionally on \((U,V)\), the variables \(X\) and \(X'\) have
the same law. Therefore
\[
 \E d_t(U,a(X'))=\E d_t(U,a(X)),
 \qquad
 \E d_t(V,v(X'))=\E d_t(V,v(X)).
\]
It follows that
\[
 \E d_{4t}(X,X')
 \le
 2\E d_t(a(X),U)+2\E d_t(v(X),V).
\]
Taking infima yields
\[
 \RD^d_\mu(4t)
 \le
 2\RD^d_{a\push\mu}(t)+2\RD^d_{v\push\mu}(t),
\]
and integration gives \eqref{eq:split-maps}.

Finally, an injective map \(\phi\) gives a bijection between
self-couplings of \(\mu\) and self-couplings of \(\phi\push\mu\).
Under this bijection the joint probability masses, and hence the
mutual information, are unchanged, while the assumed pullback identity
preserves the distortion exactly. This proves part~(v).
\end{proof}

Every rate--distortion value is nonincreasing in $t$. Consequently,
for $r_k=2^{-k}$,
\begin{equation}\label{eq:dyadic-area}
 \frac12\sum_{k\in\Z}r_k\RD^d_\mu(r_k)
 \le\cR^d(\mu)\le\sum_{k\in\Z}r_k\RD^d_\mu(r_k).
\end{equation}

The following lemma allows us to pass from finite caps to mixed finite and
infinite caps in Theorem~\ref{thm:capped-geometry}.
\begin{lemma}\label{lem:limits}
Let \(\mu\) be a finite index law and let
\(d^{(M)}=(d_i^{(M)})_i\) be coordinate costs such that
\[
 d_i^{(M)}(r)\uparrow d_i(r)
\]
for every \(i\) and \(r\ge0\). Define \(d_t^{(M)}\) and \(d_t\)
coordinatewise as above. Then, for every \(t>0\),
\[
 \RD_\mu^{d^{(M)}}(t)\uparrow\RD_\mu^d(t),
\]
and consequently
\[
 \cR^{d^{(M)}}(\mu)\uparrow\cR^d(\mu).
\]

In addition, let \(T\) be a fixed finite set and suppose that
\(\mu_M,\mu\in\cP(T)\) with \(\mu_M\to\mu\). For every fixed cost
family \(d\),
\[
 \RD_\mu^d(t)
 \le
 \liminf_{M\to\infty}\RD_{\mu_M}^d(t),
 \qquad t>0,
\]
and hence
\begin{equation}\label{eq:R-law-lsc}
 \cR^d(\mu)
 \le
 \liminf_{M\to\infty}\cR^d(\mu_M).
\end{equation}

Finally, if \(\xi^{(M)}\to\xi\) almost surely and
\(\{\|\xi^{(M)}\|_1:M\ge1\}\) is uniformly integrable, then
\[
 B_{\xi^{(M)}}(\mu)\longrightarrow B_\xi(\mu).
\]
\end{lemma}

\begin{proof}
Fix \(t>0\). Monotonicity of the costs gives
\[
 \RD_\mu^{d^{(M)}}(t)
 \le
 \RD_\mu^{d^{(M+1)}}(t)
 \le
 \RD_\mu^d(t).
\]
Let \(L\) denote the limit of the left-hand side, and choose a
minimizing self-coupling \(P_M\) for each \(M\). Since
\(\Pi(\mu,\mu)\) is compact, after passing to a subsequence we may
assume \(P_M\to P\in\Pi(\mu,\mu)\). Fix \(M_0\). For \(M\ge M_0\),
\[
 I(P_M)+\E_{P_M}d_t^{(M)}
 \ge
 I(P_M)+\E_{P_M}d_t^{(M_0)}.
\]
On the finite state space, mutual information is continuous in the
coupling, while the expectation of the fixed nonnegative extended
cost \(d_t^{(M_0)}\) is lower semicontinuous. Therefore
\[
 L
 \ge
 I(P)+\E_Pd_t^{(M_0)}.
\]
Letting \(M_0\to\infty\) and now applying monotone convergence under
the fixed coupling \(P\),
\[
 L
 \ge
 I(P)+\E_Pd_t
 \ge
 \RD_\mu^d(t).
\]
The reverse inequality follows from the first display. Thus
\[
 \RD_\mu^{d^{(M)}}(t)\uparrow\RD_\mu^d(t),
\]
and monotone convergence in \(t\) gives the corresponding assertion
for \(\cR\).

We next prove lower semicontinuity in the prescribed law. Choose
nearly minimizing couplings \(P_M\in\Pi(\mu_M,\mu_M)\). Passing to a
subsequence realizing the liminf, compactness gives \(P_M\to P\).
The marginal constraints pass to the limit, so
\(P\in\Pi(\mu,\mu)\). Since \(T\) is finite,
\[
 I(P_M)=2\Ent(\mu_M)-\Ent(P_M)\longrightarrow
 2\Ent(\mu)-\Ent(P)=I(P).
\]
The distortion term is lower semicontinuous, hence
\[
 \RD_\mu^d(t)
 \le
 I(P)+\E_Pd_t
 \le
 \liminf_M\RD_{\mu_M}^d(t).
\]
Fatou's lemma proves \eqref{eq:R-law-lsc}.

For the pairing functional, almost sure convergence together with the
uniform integrability assumption implies
\[
 \E\|\xi^{(M)}-\xi\|_1\longrightarrow0.
\]
Since \(\mu\) has finite support,
\[
 \left|
 \sup_{x\in\supp\mu}\langle x,\xi^{(M)}\rangle
 -
 \sup_{x\in\supp\mu}\langle x,\xi\rangle
 \right|
 \le
 \max_{x\in\supp\mu}\|x\|_\infty
 \|\xi^{(M)}-\xi\|_1.
\]
Taking expectations proves the last assertion.
\end{proof}

\subsection{The exponential case}

We first treat the uncapped symmetric exponential vector
\[
 Z=\xi^{[\infty]}=(\varepsilon_iE_i)_{i=1}^n.
\]
Its intrinsic cost is
\[
 \psi_\infty(r)=
 \begin{cases}
 r^2/4,&0\le r\le2,\\
 +\infty,&r>2,
 \end{cases}
\]
and throughout this subsection we denote the corresponding distortion by
\[
 d_t(x,y)=\sum_{i=1}^n
 \psi_\infty\!\left(\frac{|x_i-y_i|}{t}\right).
\]
\begin{proposition}\label{prop:uncapped-RD}
For every finite index law,
\begin{equation}\label{eq:uncapped-RD}
 \cR^{\psi_\infty}(\mu)\le CB_Z(\mu).
\end{equation}
\end{proposition}
\begin{proof}
We first bound $\cR^{\psi_\infty}(T)$ in terms of a chain. Let $T$ be finite.
Choose maps $\pi_j:T\to T$, where $\pi_0$ is constant,
$|\pi_j(T)|\le e^{2^j}$, and $\pi_j(x)=x$ for all sufficiently
large $j$. For the distortion $d_t$ above, put
\[
 \tau_j(x)=\inf\{t>0:d_t(x,\pi_j(x))\le2^j\}.
\]
If $0<t<\tau_j(x)$, then $d_t(x,\pi_j(x))>2^j$.
Lemma~\ref{lem:area-identity} therefore gives
\[
 2^j\tau_j(x)\le\int_0^\infty\min\{2^j,d_t(x,\pi_j(x))\}\,dt
             =\rho_{2^j}(x-\pi_j(x)).
\]

Fix $X\sim\mu$ on $T$ and a scale $t>0$. Let $J=J_t(X)$ be
the first $j$ with $\tau_j(X)<t$, and observe $O_t=(J,\pi_J(X))$.
Such a $J$ exists because the maps eventually equal the identity.
By monotonicity,
\[
 d_t(X,\pi_J(X))\le2^J.
\]
To bound the entropy of $O_t$, give each possible value $(j,y)$
weight $e^{-c2^j}$, with $c$ a fixed sufficiently large number.
Since there are at most $e^{2^j}$ such values at level $j$,
the sum of these weights is at most one. Normalize the weights
to a probability law $Q$. Their normalizing sum is at most one,
so $-\log Q(j,y)\le c2^j$. Nonnegativity of
$\KL(\Law(O_t)\Vert Q)$ gives
\[
 \Ent(O_t)\le\E[-\log Q(O_t)]\le c\E2^J.
\]
Let $X'$ be drawn from the conditional law of $X$ given $O_t$, independently of $X$ conditional on $O_t$. Both $X$ and $X'$
have law $\mu$, and \eqref{eq:scale-triangle} gives
\begin{align*}
 I(X;X')&\le\Ent(O_t),\\
 \E d_{2t}(X,X')
 &\le\E d_t(X,\pi_J(X))+\E d_t(X',\pi_J(X))
 \le2\E2^J.
\end{align*}
The last equality of the two expected errors follows because
$(X,O_t)$ and $(X',O_t)$ have the same law. Thus
$\RD^{\psi_\infty}_\mu(2t)\le C\E2^{J_t(X)}$.

We integrate without paying a constant at infinitely many large
scales. Set $R_*=\max_x\tau_0(x)$. For $0<t\le R_*$,
\[
 2^{J_t(x)}=1+\sum_{j<J_t(x)}2^j,
 \qquad J_t(x)>j\ \Longrightarrow\ t\le \tau_j(x).
\]
Hence
\[
 \int_0^{R_*}\E2^{J_t(X)}\,dt
 \le R_*+\E\sum_{j\ge0}2^j\tau_j(X).
\]
For $t>R_*$, use the constant observation $\pi_0(X)=x_0$.
Its entropy is zero. Since $d_t(x,x_0)\le1$ in this range,
\[
 \int_{R_*}^\infty\RD^{\psi_\infty}_\mu(2t)\,dt
 \le2\E\int_{R_*}^\infty d_t(X,x_0)\,dt
 \le2\E\rho_1(X-x_0).
\]
Also $R_*\le\max_x\rho_1(x-x_0)$. Combining these estimates
and taking the supremum over $\mu$ gives
\begin{equation}\label{eq:RD-chain}
 \cR^{\psi_\infty}(T)\le C\inf_{(\pi_j)}\max_{x\in T}
                   \sum_{j\ge0}\rho_{2^j}(x-\pi_j(x)).
\end{equation}
The strict choice $\tau_J(x)<t$ ensures finite distortion at the
scale used, including when the cost has jumps or infinite values.

For the uncapped exponential potential,
$V(u)=u^2$ for $|u|\le1$ and $V(u)=2|u|-1$ otherwise.
If $\sum_iV(u_i)\le p$, Cauchy--Schwarz and
$2|u_i|-1\ge|u_i|$ on $|u_i|>1$ give
\[
 \sum_i h_i u_i
 \le\sqrt p\norm h_2+p\norm h_\infty.
\]
Conversely, $V(u)\le u^2$ everywhere, so the Euclidean ball
of radius $\sqrt p$ is feasible. For $p\ge1$, a single
coordinate of size $(p+1)/2$ is feasible as well. Consequently,
\[
 \rho_p(h)\asympc\sqrt p\norm h_2+p\norm h_\infty,
 \qquad p\ge1.
\]
Choose $\pi_j(x)$ in the partition cell containing $x$.
The distance to this representative is bounded by the diameter
of the cell. Lemma~\ref{lem:laplace-moments},
\eqref{eq:RD-chain}, and \eqref{eq:reference-comparison} yield
\[
 \cR^{\psi_\infty}(T)\le C\gamma_Z(T)\le C'S_Z(T).
\]

We then lift this bound to an arbitrary index law using Liu's lifting method \cite{Liu}, which replaces a prescribed index law by long tuples with fixed empirical frequencies. 

Retain a prescribed law $\mu=\sum_{j=1}^m\mu_j\delta_{x_j}$.
First suppose its masses are positive and rational. Let $\nu_N$
be uniform on the set $T_N(\mu)$ from Lemma~\ref{lem:type}.
Each coordinate block of a vector $X^N=(X_1,\ldots,X_N)$
with this law has marginal $\mu$, by permutation symmetry.
Its entropy and entropy deficit are
$$
H_N=\Ent(X^N)=\log|T_N(\mu)|,
\qquad
\Delta_N=N\Ent(\mu)-H_N.
$$
By the standard method-of-types estimate,
$$
0\le \Delta_N\le m\log(N+1);
$$
see \cite[Section~2.1]{CK} and \cite[Chapter~11]{CT}.

Take any self-coupling $(X^N,X'{}^N)$ of $\nu_N$. Conditional
entropy decreases under additional conditioning. It follows that
\begin{align*}
 I(X^N;X'{}^N)
 &=H_N-\Ent(X^N\mid X'{}^N)\\
 &\ge H_N-\sum_{k=1}^N\Ent(X_k\mid X'{}^N)\\
 &\ge H_N-\sum_{k=1}^N\Ent(X_k\mid X'_k)\\
 &=\sum_{k=1}^NI(X_k;X'_k)-\Delta_N.
\end{align*}
The distortions add over the blocks, and each $(X_k,X'_k)$
is a self-coupling of $\mu$. Hence
\[
 \RD^{\psi_\infty}_{\nu_N}(t)
 \ge N\RD^{\psi_\infty}_\mu(t)-\Delta_N.
\]
Integrating the constant error over all positive scales would
not be valid. Instead choose
$M_N=\sqrt N/\log(N+2)$ and integrate over $0<t<M_N$:
\begin{align*}
 \int_0^{M_N}\RD^{\psi_\infty}_\mu(t)\,dt
 &\le\frac1N\cR^{\psi_\infty}(\nu_N)+\frac{M_N\Delta_N}{N}\\
 &\le\frac C N S_{Z^{\oplus N}}(T_N(\mu))
                           +\frac{M_N\Delta_N}{N}.
\end{align*}
Here the second line uses the set bound in dimension $Nd$.
As $N\to\infty$ along the allowed integers, $M_N\to\infty$
and $M_N\Delta_N/N\to0$. Lemma~\ref{lem:type} gives the
limit of the first term. We obtain
$\cR^{\psi_\infty}(\mu)\le CB_Z(\mu)$.
Finally, approximate an arbitrary law on the same finite support by laws with
positive rational masses. The lower semicontinuity
\eqref{eq:R-law-lsc} and the continuity of the pairing functional in
the masses give the same inequality for the limiting law.
\end{proof}

\subsection{Truncated exponentials}\label{sec:capped-geometry}

We now prove two-sided bounds for general caps. The proof constructs
the decomposition from Proposition~\ref{prop:decomposition-dual}
by clipping differences between successive approximations of the
index. 
Recall that $
\clip_\beta(x)
=
\bigl(
\operatorname{sgn}(x_i)\min\{|x_i|,\beta_i\}
\bigr)_{i=1}^n.
$

For $1\le\beta\le\infty$ set
\begin{equation}\label{eq:capped-cost}
 \psi_\beta(r)=\begin{cases}
 r^2/4,&0\le r\le2,\\
 2\beta-1,&r>2.
 \end{cases}
\end{equation}
These are the intrinsic costs of $\eps\min(E,\beta)$,
including the value at the slope threshold. For a cap vector \(\beta=(\beta_i)_{i=1}^n\in[1,\infty]^n\),
write
\[
  \psi_\beta=(\psi_{\beta_i})_{i=1}^n,
\]
and let \(\cR^{\psi_\beta}\) denote the corresponding
rate--distortion functional.

\begin{theorem}\label{thm:capped-geometry}
For all finite index laws and all coordinate caps in $[1,\infty]$,
\begin{equation}\label{eq:capped-geometry}
 \cR^{\psi_\beta}(\mu)\asympc B_{\xi^{[\beta]}}(\mu).
\end{equation}
For finite caps both quantities are also comparable to
$\delta_\beta(\mu)$, with numerical constants.
\end{theorem}
\begin{proof}
First suppose every cap is finite, and write
\[
 d_i(r)=\chi(r)+\beta_i\ind_{\{r>1\}}.
\]
Directly from \eqref{eq:capped-cost},
$\tfrac12d_i(r)\le\psi_{\beta_i}(2r)\le2d_i(r)$.
It is therefore enough to compare $\cR^d$ and $\delta_\beta$.

We first prove
\[
 \delta_\beta(\mu)\lesssim \cR^d(\mu).
\]
Starting from near-optimal rate--distortion codes at a sequence of
dyadic scales, we construct a decomposition
\[
 x=a(x)+v(x)
\]
of the type appearing in Proposition~\ref{prop:decomposition-dual}.
The discarded part controls
$\E\sum_i\beta_i|a_i(X)|$, while the retained part controls
$B_Z(v\push\mu)$.
Take $r_k=2^{-k}$ and integers
$m<N$ with $r_m\ge2\diam_\infty(\supp\mu)$.
For $m<k<N$, choose minimizing self-coupling codes $U_k$
at scale $r_k$, conditionally independently given $X\sim\mu$.
Let $U_m\sim\mu$ be independent of $X$, and $U_N=X$.
Put
\[
 h_k=I(X;U_k),\quad e_k=\E d_{r_k}(X,U_k),\quad
 g_k=\E\sum_i\chi(\abs{X_i-U_{k,i}}/r_k).
\]
For $m<k<N$, the choice of $U_k$ gives
\[
 h_k+e_k=\RD_\mu^d(r_k),\qquad g_k\le e_k.
\]
We then define
\begin{equation}\label{eq:clipped-map}
 \Delta_k=\clip_{3r_k/2}(U_{k+1}-U_k),\quad
 \widetilde v=U_m+\sum_{k=m}^{N-1}\Delta_k,\quad
 v(x)=\E[\widetilde v\mid X=x],\quad a(x)=x-v(x).
\end{equation}
Without clipping, the increments would telescope to
$U_N-U_m=X-U_m$. Thus
\[
 X-\widetilde v
 =\sum_{k=m}^{N-1}\bigl[(U_{k+1}-U_k)-\Delta_k\bigr],
 \qquad a(X)=\E[X-\widetilde v\mid X].
\]
We first bound this discarded sum. Fix a coordinate $i$ and call
an index $k$ bad if $\abs{X_i-U_{k,i}}>r_k$. Both endpoints $m,N$ are good by our choices of $U_m,U_N$.
If $k,k+1$ are good, then
\[
 |U_{k+1,i}-U_{k,i}|
 \le |U_{k+1,i}-X_i|+|X_i-U_{k,i}|
 \le r_{k+1}+r_k=3r_k/2.
\]
Such an edge discards nothing. For a
maximal run $p,\ldots,q$ of bad indices, the adjacent indices
$p-1,q+1$ are good. The raw increments telescope over its edge
block to $U_{q+1,i}-U_{p-1,i}$. Therefore
\[
 \left|\sum_{k=p-1}^q
     \bigl[(U_{k+1,i}-U_{k,i})-\Delta_{k,i}\bigr]\right|
 \le |U_{q+1,i}-X_i|+|X_i-U_{p-1,i}|
                      +\sum_{k=p-1}^q|\Delta_{k,i}|.
\]
The right side is at most
\[
 r_{p-1}+r_{q+1}+\tfrac32\sum_{k=p-1}^q r_k
                                                \le\tfrac{17}{2}r_p.
\]
The edge blocks belonging to different runs are disjoint.
Each run starts at one bad index $p$, so
\[
 |X_i-\widetilde v_i|
 \le\tfrac{17}{2}\sum_{k=m+1}^{N-1}
       r_k\ind_{\{|X_i-U_{k,i}|>r_k\}}.
\]
Multiply by $\beta_i$, sum, and take expectations.
Conditional Jensen bounds $\E\sum_i\beta_i|a_i(X)|$ by
$\E\sum_i\beta_i|X_i-\widetilde v_i|$. The plateau term in
$d_i$ bounds each weighted probability above. We obtain
\begin{equation}\label{eq:bad-runs}
 \E\sum_i\beta_i\abs{a_i(X)}
                     \le\tfrac{17}{2}\sum_{k=m+1}^{N-1}r_ke_k.
\end{equation}

We next bound $B_Z(v\push\mu)$.
For $a,b\ge0$,
\[
 \min\{3/2,a+b/2\}^2\le\frac94\{\chi(a)+\chi(b)\}.
\]
If $\max(a,b)>1$, use the bound $9/4$ on the left. Otherwise
use $(a+b/2)^2\le(5/4)(a^2+b^2)$. Taking
$a=|X_i-U_{k,i}|/r_k$ and $b=|X_i-U_{k+1,i}|/r_{k+1}$
gives
\begin{equation}\label{eq:clipped-l2}
 \norm{\Delta_k}_2^2\le\tfrac94r_k^2\sum_i
 \left[\chi\left(\frac{\abs{X_i-U_{k,i}}}{r_k}\right)
       +\chi\left(\frac{\abs{X_i-U_{k+1,i}}}{r_{k+1}}\right)\right].
\end{equation}
The Laplace moment generating function gives the entropy bound
\begin{equation}\label{eq:laplace-entropy}
 \E\ip{w(O)}Z\le sI(Z;O)+\frac2s\E\norm{w(O)}_2^2
            \quad\text{if }\norm{w(O)}_\infty\le s/2.
\end{equation}
Indeed, $\log\E e^{\ip\lambda Z}\le2\norm\lambda_2^2$
when $\norm\lambda_\infty\le1/2$. Apply the relative-entropy
variational inequality \eqref{eq:entropy-variational} to
$P=\Law(Z\mid O=o)$ and $Q=\Law(Z)$, with test
$z\mapsto\ip{w(o)}z/s$. Multiply by $s$ and average in $o$.
The averaged relative entropy is exactly $I(Z;O)$, which
proves \eqref{eq:laplace-entropy}.

In any coupling of $X$ with $Z$, generate all the codes by
the same fixed kernel conditionally independently of $Z$
given $X$. For $O_k=(U_k,U_{k+1})$, data processing and
conditional independence give
$I(Z;O_k)\le I(X;O_k)\le h_k+h_{k+1}$.
Using \eqref{eq:laplace-entropy} with $s=3r_k$ and
\eqref{eq:clipped-l2}, we get
\begin{align*}
 \sum_{k=m}^{N-1}\E\ip{\Delta_k}Z
 &\le\sum_{k=m}^{N-1}
       [3r_k(h_k+h_{k+1})+\tfrac32r_k(g_k+g_{k+1})]\\
 &=9\sum_{k=m+1}^{N-1}r_kh_k+6r_N\Ent(\mu)
       +\tfrac92\sum_{k=m+1}^{N-1}r_kg_k+\tfrac32r_mg_m.
\end{align*}
Here $h_m=0$, $h_N=\Ent(\mu)$, $g_N=0$, and
$r_mg_m=\E\norm{X-U_m}_2^2/r_m$.
The root $U_m$ is independent of $Z$ in this construction,
so $\E\ip{U_m}Z=0$. The conditional law of the codes given
$X$ is fixed, independently of which coupling of $X,Z$ was
chosen. It follows that
\[
 \E[\widetilde v\mid X,Z]=\E[\widetilde v\mid X]=v(X),
 \qquad \E\ip{\widetilde v}Z=\E\ip{v(X)}Z.
\]
Thus the same map $v$ works for every coupling. Taking the
supremum and using \eqref{eq:fiber} gives the asserted bound
for $B_Z(v\push\mu)$.
Combining this with \eqref{eq:bad-runs}, $g_k\le e_k$,
and Proposition~\ref{prop:decomposition-dual} gives the
explicit finite-scale bound
\begin{equation}\label{eq:finite-clipping}
 \delta_\beta(\mu)\le13\sum_{k=m+1}^{N-1}r_k\RD^d_\mu(r_k)
       +\frac{3\E\norm{X-U_m}_2^2}{2r_m}+6r_N\Ent(\mu).
\end{equation}
Put $D_2=\diam_2(\supp\mu)$. The two endpoint terms obey
\[
 \frac{3\E\norm{X-U_m}_2^2}{2r_m}\le\frac{3D_2^2}{2r_m}
       \longrightarrow0\quad(m\to-\infty),
 \qquad
 6r_N\Ent(\mu)\longrightarrow0\quad(N\to\infty).
\]
All terms in the sum are nonnegative. Passing to the limit in
\eqref{eq:finite-clipping} and using \eqref{eq:dyadic-area} gives
\[
 \delta_\beta(\mu)\le13\sum_{k\in\Z}r_k\RD^d_\mu(r_k)
                         \le26\cR^d(\mu).
\]

We next prove the reverse bound
\[
 \cR^d(\mu)\lesssim\delta_\beta(\mu).
\]
Starting from an arbitrary decomposition
\(x=a(x)+v(x)\) in Proposition~\ref{prop:decomposition-dual},
we estimate the rate--distortion costs of the two components
separately and then combine them using
Lemma~\ref{lem:coding}.

Take any map $a$ in \eqref{eq:decomposition-dual} and write
\[
v(x)=x-a(x).
\]
We first estimate the cost of the $a$-part. Let $X,X'$ be independent
copies with law $\mu$. Using the independent coupling in the definition
of $\cR^d(a\push\mu)$, Tonelli's theorem, and the identity
\[
\int_0^\infty d_i(r)r^{-2}\,dr=\beta_i+2,
\]
we obtain
\[
\cR^d(a\push\mu)
\le
\E\sum_i(\beta_i+2)
|a_i(X)-a_i(X')|.
\]
By the triangle inequality and symmetry of the independent copies,
\[
\E|a_i(X)-a_i(X')|
\le
2\E|a_i(X)|,
\]
and since $\beta_i\ge1$,
\[
\cR^d(a\push\mu)
\le
6\E\sum_i\beta_i|a_i(X)|.
\]

For the remainder part, observe that
\[
d_i(r)\le\psi_\infty(2r),
\]
so the uncapped estimate from Proposition~\ref{prop:uncapped-RD}
implies
\[
\cR^d(v\push\mu)
\le
C B_Z(v\push\mu).
\]
Combining the two estimates with the decomposition
\[
x=a(x)+v(x)
\]
and applying \eqref{eq:split-maps}, we obtain the desired upper bound
after taking the infimum over all admissible maps $a$. Hence
\[
\cR^d(\mu)\le C\delta_\beta(\mu).
\]
Theorem~\ref{thm:capping} together with the scaling properties of the
costs then yields \eqref{eq:capped-geometry} for finite caps.

\emph{Infinite caps.}
For mixed finite and infinite caps, replace each $\beta_i$
by $\min(\beta_i,M)$ and let $M\uparrow\infty$.
The costs increase to the stated costs, including their
infinite values. Under the natural coupling the source
vectors converge in $L^1$. Lemma~\ref{lem:limits} passes
both sides to the limit with the same constants.
\end{proof}
The Bernoulli formula follows immediately from the preceding theorem.
\begin{corollary}\label{cor:bernoulli}
For every finite index law,
\begin{equation}\label{eq:bernoulli-RD}
 \cR^\chi(\mu)\asympc B_\eps(\mu).
\end{equation}
\end{corollary}
\begin{proof}
Take all caps equal to one. Then
\[
 \psi_1(r)=\chi(r/2),
\]
and
\[
 (1-e^{-1})\varepsilon\preceq_{\rm cx}\xi^{[1]}
 \preceq_{\rm cx}\varepsilon
\]
in each coordinate. Applying Lemma~\ref{lem:convex-pairing}(ii) to the independent
coordinates and then part~(i) gives
\[
 B_{\xi^{[1]}}(\mu)\asymp B_\varepsilon(\mu).
\]
Since \(\psi_1(r)=\chi(r/2)\), Theorem~\ref{thm:capped-geometry}
and Lemma~\ref{lem:coding}(i) give
\[
 \cR^\chi(\mu)\asymp B_\varepsilon(\mu).
\]
\end{proof}

\section{General log-concave tails}\label{sec:general-tails}

We now pass from truncated exponentials to general log-concave tails.
For each coordinate we construct an auxiliary sum of independent
truncated exponential variables, with coefficients determined by the
slopes of its tail potential. We then show that this expansion
preserves, up to universal constants, both the pairing functional and
the associated rate--distortion cost. At the end of the section we
prove the fixed-law estimate
\[
 B_Y(\mu)\asymp\cR^\Psi(\mu),
\]
which is the first assertion of Theorem~\ref{thm:main}.

For $j\ge0$, put $s_j=2^{j+1}$ and $a_j=2/s_j=2^{-j}$.
For each coordinate, decompose the tail region according to the size
of the slope by setting
\begin{equation}\label{eq:bins}
 A_{ij}=\{u\in[1,b_i):s_j\le V'_{i,+}(u)<2s_j\},
 \qquad
 w_{ij}=s_j\abs{A_{ij}}.
\end{equation}
Since $V'_{i,+}$ is nondecreasing, every nonempty $A_{ij}$ is an
interval. On this interval the slope varies by at most a factor two,
and hence $w_{ij}$ represents, up to a factor two, the contribution
of $A_{ij}$ to the intrinsic cost $\Psi_i$ once the corresponding
slope level has been reached.

We retain those intervals for which $w_{ij}\ge1$ and set
$\beta_{ij}=(w_{ij}+1)/2$, allowing $\beta_{ij}=\infty$.
If $b_i=1$, there are no tail intervals. The expanded coordinate set is
\[
 \mathcal I=\{(i,*):1\le i\le n\}
 \mathbin{\dot\cup}
 \{(i,j):j\ge0,\ w_{ij}\ge1\}.
\]
The label $*$ corresponds to the quadratic core, while the labels
$(i,j)$ correspond to the retained tail intervals.

For $\alpha=(i,*)$ put $\widehat\beta_\alpha=1$, and for
$\alpha=(i,j)$ put $\widehat\beta_\alpha=\beta_{ij}$. Let
$\Xi=(\Xi_\alpha)_{\alpha\in\mathcal I}$ have independent coordinates
\[
 \Xi_{i,*}=\xi_{i,*}^{[1]},\qquad
 \Xi_{i,j}=\xi_{ij}^{[\beta_{ij}]}.
\]
We write $\cR^{\rm cap}$ for the rate--distortion functional on the
expanded coordinate space whose coordinate costs are
$\psi_{\widehat\beta_\alpha}$. Define
\begin{equation}\label{eq:source-expansion}
 H_i=\Xi_{i,*}
     +\sum_{j:w_{ij}\ge1}a_j\Xi_{i,j},
 \qquad
 (Lx)_{i,*}=x_i,\quad (Lx)_{ij}=a_jx_i.
\end{equation}
Thus the original $i$-th coordinate is replaced by one core coordinate
and the retained tail coordinates.

For each fixed $i$, the coefficients satisfy
\[
 1+\sum_{j:w_{ij}\ge1}a_j
 \le1+\sum_{j\ge0}2^{-j}=3.
\]
This bound is uniform in the number of retained intervals and in the
size of the slope jumps. Since $\E\abs{\xi^{[\beta]}}\le1$ and
$\norm{\xi^{[\beta]}}_p\le\norm{\varepsilon E}_p$, the series defining $H_i$
converges absolutely almost surely, in $L^1$, and in every fixed $L^p$.

The map $L$ takes values in $\ell^1(\mathcal I)$ and is injective
because of the core coordinates. For a finite index law, the pairing
with the expanded source is well defined, since
\[
 \E\sum_{\alpha\in\mathcal I}
       |(LX)_\alpha\Xi_\alpha|
 \le3\sum_i\max_{x\in\supp\mu}|x_i|<\infty.
\]
Consequently,
\[
 \langle LX,\Xi\rangle=\sum_iX_iH_i
\]
almost surely.

The expansion replaces the original coordinate $x_i$ by the family of
expanded coordinates
\[
 \bigl(x_i,(a_jx_i)_{j:w_{ij}\ge1}\bigr).
\]
The corresponding cost on the original coordinate, obtained by pulling
back the coordinate costs of the expanded capped exponential vector
through the embedding $L$, is
\begin{equation}\label{eq:expanded-cost}
\Phi_i(r)
=
\chi(r/2)
+
\sum_{j:w_{ij}\ge1}
\psi_{\beta_{ij}}(a_jr).
\end{equation}
To compare the original cost $\Psi_i$ with $\Phi_i$, we also introduce the auxiliary cost
\[
E_0(r)=\sum_{j\ge0}\chi(r/s_j).
\]
This is the cost associated with the full dyadic expansion of an
uncapped exponential coordinate. It will be used only as an auxiliary Bernoulli-type cost controlling
the discarded bins and the additional quadratic terms introduced by
the expansion.

For $r\le2$, every term in $\Phi_i(r)$ lies in the quadratic regime of
the corresponding coordinate cost. Since
\[
\sum_j a_j^2<\infty,
\]
the quadratic contributions remain summable. The same argument applies
to $E_0$. Hence both the induced costs $\Phi_i$ and the auxiliary cost
$E_0$ satisfy a quadratic bound near the origin.

\subsection{Comparing the costs}

The decomposition of the original coordinate into truncated exponential
components introduces an error in the induced cost. We control this
error by the auxiliary cost $E_0$ and show that it is absorbed by the
quadratic part already present in the original and expanded costs.

\begin{proposition}\label{prop:cost-transfer}
For all $i$ and $r\ge0$,
\begin{equation}\label{eq:cost-transfer}
 \Psi_i(r)\le2\Phi_i(r)+2E_0(r),\qquad
 \Phi_i(r/2)\le\Psi_i(r)+E_0(r/2),\qquad
 \int_0^\infty\frac{E_0(r)}{r^2}\,dr=2.
\end{equation}
Consequently, for every finite index law,
\begin{equation}\label{eq:area-transfer}
 \cR^\Psi(\mu)\asympc\cR^\Phi(\mu)
                      =\cR^{\rm cap}(L\push\mu).
\end{equation}
For a countable expansion the expression on the right is the
increasing limit of finite-coordinate costs on the same finite
label space.
\end{proposition}
\begin{proof}
Put
$\Theta_i(r)=\chi(r/2)+\sum_{j\ge0}w_{ij}\ind_{\{r>s_j\}}$.
For a point $u\in A_{ij}$, we have $s_j\le V'_{i,+}(u)<2s_j$.
If $r>2s_j$, the whole interval contributes to $\Psi_i(r)$,
and its contribution is at least $w_{ij}$. If $r\le s_j$,
it contributes nothing. In all cases its contribution is at
most $2w_{ij}\ind_{\{r>s_j\}}$. Including the quadratic
part gives
\[
 \Theta_i(r/2)\le\Psi_i(r)\le2\Theta_i(r).
\]
For a kept bin,
\[
 \psi_{\beta_{ij}}(a_jr)
 =w_{ij}\ind_{\{r>s_j\}}+(r/s_j)^2\ind_{\{r\le s_j\}}.
\]
Separate the discarded steps and the added quadratic terms:
\[
 A_i(r)=\sum_{j:w_{ij}<1}w_{ij}\ind_{\{r>s_j\}},\qquad
 Q_i(r)=\sum_{j:w_{ij}\ge1}(r/s_j)^2\ind_{\{r\le s_j\}}.
\]
Both are bounded by $E_0(r)$, and
$\Theta_i=\Phi_i-Q_i+A_i$. Therefore
\begin{align*}
 \Psi_i(r)&\le2\Theta_i(r)\le2\Phi_i(r)+2E_0(r),\\
 \Phi_i(r/2)&\le\Theta_i(r/2)+Q_i(r/2)
                         \le\Psi_i(r)+E_0(r/2).
\end{align*}
All thresholds are strict where the step term is used. If a
bin has infinite length, the same inequalities hold as inequalities
of nonnegative extended values; alternatively, increase finite
partial lengths to the full length. Tonelli and
$\int_0^\infty\chi(r/s)r^{-2}dr=2/s$ give the last identity.

We next bound the error for each prescribed law $\mu$. For a finite sum in $E_0$,
pull back the Bernoulli cost by $x\mapsto(x_i/s_j)_{i,j}$.
The corresponding source in the original coordinates is
$W_i^{(J)}=\sum_{j=0}^J\eps_{ij}/s_j$. Its coordinates are
independent and centered and $\abs{W_i^{(J)}}\le1$.
For a centered scalar $W\in[-1,1]$ and a convex test $f$,
\[
 \E f(W)\le\E\left[\frac{1+W}{2}f(1)+\frac{1-W}{2}f(-1)\right]
                  =\frac{f(1)+f(-1)}2.
\]
Thus $W_i^{(J)}\cx\eps_i$. The pullback identity,
Corollary~\ref{cor:bernoulli}, and independent products give
\[
 \cR^{E_0^{(J)}}(\mu)\asympc B_{W^{(J)}}(\mu)
                                      \le B_\eps(\mu).
\]
Increase $J$ and use Lemma~\ref{lem:limits}.
Both \(\Psi_i\) and \(\Phi_i\) dominate \(\chi(r/2)\).
Hence Lemma~\ref{lem:coding}(i), monotonicity of the
rate--distortion value in the cost, and Corollary~\ref{cor:bernoulli}
give
\[
  \cR^\Psi(\mu),\ \cR^\Phi(\mu)
  \ge c\,B_\varepsilon(\mu).
\]
Lemma~\ref{lem:coding} gives
\begin{align*}
 \cR^\Psi(\mu)&\le C[\cR^\Phi(\mu)+B_\eps(\mu)]
                    \le C'\cR^\Phi(\mu),\\
 \cR^\Phi(\mu)&\le C[\cR^\Psi(\mu)+B_\eps(\mu)]
                    \le C'\cR^\Psi(\mu).
\end{align*}
The equality in \eqref{eq:area-transfer} is the exact
pullback of the expanded capped costs through the injective
map $L$, first for finite expansions and then by monotone
convergence.
\end{proof}

\subsection{Comparing the variables}

Cost comparison alone does not compare the random sources. We therefore
compare their moments and then convert the resulting moment estimates into
convex order. In particular, once \(\|H_i\|_p\asymp q_i(p)\) is
established uniformly for \(p\ge2\), the next lemma gives the required
two-sided convex-order comparison.

\begin{lemma}\label{lem:quantile-convex}
Let $Y$ be a normalized symmetric log-concave-tailed variable,
with quantile $q(1)=1$. Suppose $A$ is symmetric and, for
some positive constants $c,C,c_0$,
\[
 cq(p)\le\norm A_p\le Cq(p)\quad(p\ge2),\qquad
 \E\abs A\ge c_0.
\]
Then there is a constant $C_1>0$, depending on $c,C,c_0$, such that
$C_1^{-1}Y\cx A\cx C_1Y$.
\end{lemma}
\begin{proof}
Concavity and $q(0)\ge0$ give $q(\lambda p)\le\lambda q(p)$ for $\lambda\ge1$.
Markov's inequality gives
\[
 \Prob(|A|>eCq(p))\le\frac{\E|A|^p}{(eCq(p))^p}\le e^{-p}.
\]
Also $\|A\|_{2p}\le(2C/c)\|A\|_p$. Apply Paley--Zygmund
to $|A|^p$ at level $2^{-p}\E|A|^p$:
\begin{align*}
 \Prob\{|A|\ge\tfrac12\|A\|_p\}
 &\ge(1-2^{-p})^2
                 \left(\frac{\|A\|_p}{\|A\|_{2p}}\right)^{2p}\\
 &\ge\tfrac14(c/(2C))^{2p}\ge e^{-\kappa p},
\end{align*}
where $\kappa\ge1$ depends only on $c,C$. With $\delta_0=c/4$ we obtain
\begin{equation}\label{eq:moment-tail}
 \Prob(|A|>eCq(p))\le e^{-p},\qquad
 \Prob(|A|>\delta_0q(p))\ge e^{-\kappa p}\quad(p\ge2).
\end{equation}
The strict lower threshold is valid because
$\frac12\|A\|_p\ge\frac c2q(p)>\delta_0q(p)$.

Let $b=\lim_{p\to\infty}q(p)$. The inverse tail satisfies
\begin{equation}\label{eq:quantile-strict-tail}
 \Prob(|Y|>q(p))\le e^{-p},
 \qquad q(N(y))=y\quad(0\le y<b,\ N(y)>0).
\end{equation}
At $y<b$ the tail exponent is continuous, so
$\Prob(|Y|>y)=e^{-N(y)}$. The first inequality in
\eqref{eq:quantile-strict-tail} also holds at a finite endpoint.
It is a strict-tail assertion, so an endpoint atom causes no loss.

We first prove $|A|\st M(1+|Y|)$ for a sufficiently large constant
$M$. Choose $M\ge\max(2eC,C)$. For $y\ge0$ there are three cases.
If $y<b$ and $N(y)\ge2$, use \eqref{eq:moment-tail} at
$p=N(y)$ and $M(1+y)\ge eCy$:
\[
 \Prob(|A|>M(1+y))\le e^{-N(y)}=\Prob(|Y|>y).
\]
If $y<b$ and $N(y)<2$, then $q(2)\le2$ and
\[
 \Prob(|A|>M(1+y))\le\Prob(|A|>eCq(2))
                 \le e^{-2}\le\Prob(|Y|>y).
\]
If $b<\infty$ and $y\ge b$, the moment upper bound gives
$\|A\|_\infty=\lim_{p\to\infty}\|A\|_p\le Cb$.
The left side is then zero. Thresholds below $M$ require no
estimate since $M(1+|Y|)\ge M$. Thus
\begin{equation}\label{eq:stochastic-plus-one}
 |A|\st M(1+|Y|).
\end{equation}

For the reverse comparison, enlarge $M$ so that
$M\ge\max(2\kappa,\kappa/\delta_0,1/\delta_0)$. If $0\le r<\delta_0q(2)$, then
$M(1+r)\ge q(2\kappa)$ and
\[
 \Prob(|Y|>M(1+r))\le e^{-2\kappa}
                       \le\Prob(|A|>\delta_0q(2))\le\Prob(|A|>r).
\]
If $\delta_0q(2)\le r<\delta_0 b$, continuity gives $r=\delta_0q(p)$ for some $p\ge2$.
Since $q(\kappa p)\le \kappa q(p)=(\kappa/\delta_0)r\le M(1+r)$,
\[
 \Prob(|Y|>M(1+r))\le e^{-\kappa p}\le\Prob(|A|>r).
\]
If $b<\infty$ and $r\ge \delta_0 b$, then $M(1+r)\ge b$ and the
left side vanishes. This proves $|Y|\st M(1+|A|)$.

It remains to remove the added one. For a symmetric variable $B$
with $\E|B|\ge b_0>0$, Lemma~\ref{lem:convex-pairing} gives
$b_0\eps\cx B$. Put $C_*=2M\max(1,b_0^{-1})$. For every
convex scalar function $f$, an independent sign $\eps$ gives
\begin{align*}
 \E f\bigl(M\eps(1+|B|)\bigr)
 &\le\tfrac12\E f(2M\eps)+\tfrac12\E f(2M\eps|B|)\\
 &\le\tfrac12\E f(C_*B)+\tfrac12\E f(C_*B)=\E f(C_*B).
\end{align*}
In the last line we used $\eps|B|\stackrel{\rm d}=B$ and
$uB\cx vB$ for $0\le u\le v$, which follows from convexity
and $\E B=0$. Apply this calculation to $B=Y$ and to $B=A$.
The required lower bounds are $\E|Y|\ge e^{-1}$ and
$\E|A|\ge c_0$. The two stochastic comparisons now give the
two convex orders in the lemma.
\end{proof}

\begin{proposition}\label{prop:source-transfer}
The auxiliary coordinates $H_i$ in \eqref{eq:source-expansion}
satisfy, coordinatewise,
\begin{equation}\label{eq:source-transfer}
 C^{-1}Y_i\cx H_i\cx CY_i.
\end{equation}
Thus $B_H(\mu)\asympc B_Y(\mu)$ for every finite index law,
with numerical constants independent of the tails.
\end{proposition}
\begin{proof}
Fix a coordinate and omit its index. We proceed in three steps:
compare the one-dimensional support values, compare the moments
of the expanded sum with those values, and apply
Lemma~\ref{lem:quantile-convex}.

\emph{Support values.} Let
\[
 Q(p)=\sup\{u\ge0:V(u)\le p\},\qquad
 Q_\Phi(p)=\int_0^\infty\min\{p,\Phi(1/t)\}\,dt.
\]
The area identity and \eqref{eq:cost-transfer} give
\begin{equation}\label{eq:profile-comparison}
 Q(p)\le2Q_\Phi(p)+4,\qquad Q_\Phi(p)\le2Q(p)+2.
\end{equation}
Both profiles are at least one for $p\ge1$, hence are
comparable there. If $b>1$, then
$Q(p)=q((p+1)/2)\asympc q(p)$ for $p\ge1$.
If $b=1$, $Q(p)=q(p)=1$ for $p\ge1$.
Also \eqref{eq:quantile-representation} gives directly
\begin{equation}\label{eq:quantile-moment}
 e^{-1}q(p)\le\norm Y_p
 \le q(p)\left(1+\frac{\Gamma(p+1)}{p^p}\right)^{1/p}
 \le Cq(p),\qquad p\ge1.
\end{equation}
For the lower bound restrict to $E\ge p$; for the upper
bound use $q(E)\le q(p)$ when $E\le p$ and
$q(E)\le(E/p)q(p)$ otherwise.

\emph{Moments of a finite sum.}
We prove $\norm H_p\asympc Q_\Phi(p)$ for $p\ge2$.
First take a finite expansion of this coordinate and denote its nonnegative
coefficients, including the core coefficient one, by $c_j$.
During this finite-expansion argument, $\Phi$ denotes the
corresponding partial cost.
Their sum is at most three. Thus
$H=\sum_jc_j\xi_j^{[\beta_j]}$ with $c_j\ge0$ and
$\sum_jc_j\le3$. Let $V_\beta(u)$ equal $u^2$
for $\abs u\le1$, $2\abs u-1$ for $1<\abs u\le\beta$,
and infinity outside the cap. The area identity identifies
$Q_\Phi(p)$ with the support value
\begin{equation}\label{eq:column-support}
 \sup\left\{\sum_jc_ju_j:\sum_jV_{\beta_j}(u_j)\le p\right\}.
\end{equation}
For $f_\beta=V_\beta^*$ and $z\ge0$,
\begin{equation}\label{eq:capped-mgf}
 \log\E e^{z\xi^{[\beta]}}\le f_\beta(8z).
\end{equation}
Indeed, if $z\le1/4$, use
$-\log(1-z^2)\le2z^2$ and $f_\beta(8z)=16z^2$.
For $1/4\le z\le1/2$ the left side is at most
$\log(4/3)$ and the right side at least one.
For $z\ge1/2$ and finite $\beta$, the left side is at
most $\beta z$ and
$f_\beta(8z)=8\beta z-2\beta+1\ge\beta z$.
For an infinite cap the right side is infinite in this
last range.

To justify the Chernoff step without an attainment
assumption, define
\[
 J(v)=\inf\left\{\sum_jV_{\beta_j}(u_j):\sum_jc_ju_j=v\right\}.
\]
For every finite $p$, compactness of the scalar sublevels gives
\[
 \{v:J(v)\le p\}
 =\left\{\sum_jc_ju_j:\sum_jV_{\beta_j}(u_j)\le p\right\}.
\]
The right side is a compact linear image. This proves that $J$
is closed. Mixing feasible vectors proves convexity, and $J(0)=0$.
Its conjugate is
\begin{align*}
 J^*(z)
 &=\sup_{(u_j)}\left\{z\sum_jc_ju_j-\sum_jV_{\beta_j}(u_j)\right\}\\
 &=\sum_j\sup_{u_j}\{c_jzu_j-V_{\beta_j}(u_j)\}
   =\sum_jf_{\beta_j}(c_jz).
\end{align*}
For every $z\ge0$, independence and \eqref{eq:capped-mgf} give
\[
 \log\E e^{zH}\le\sum_jf_{\beta_j}(8c_jz)=J^*(8z).
\]
Markov's inequality and the substitution $r=8z$ now yield
\begin{align*}
 \Prob(H>v)
 &\le\inf_{z\ge0}\exp[-zv+J^*(8z)]\\
 &=\exp\left[-\sup_{r\ge0}\{r(v/8)-J^*(r)\}\right]
 =e^{-J(v/8)},\qquad v\ge0.
\end{align*}
The last equality uses $J^{**}=J$, the biconjugacy theorem for
a proper closed convex function \cite[Theorem~12.2]{Rockafellar}.
In this formula
$J^{**}(v)=\sup_{r\in\R}\{rv-J^*(r)\}$.
Since $J$ and $J^*$ are even, $v\ge0$ permits restriction to
$r\ge0$. By \eqref{eq:column-support},
$J(2Q_\Phi(p))\ge p$ or the argument is outside the
effective domain. Hence
\[
 \Prob(\abs H>16Q_\Phi(p))\le2e^{-p}.
\]
The profile is concave and vanishes at zero, so
$Q_\Phi(vp)\le vQ_\Phi(p)$ for $v\ge1$.
Applying the last tail bound at budget $vp$ and
integrating gives
\[
 \E\abs H^p\le(16Q_\Phi(p))^p
       \left(1+2p\int_1^\infty v^{p-1}e^{-pv}\,dv\right)
 \le(CQ_\Phi(p))^p.
\]

For the lower bound take nonnegative feasible $u_j$ in
\eqref{eq:column-support}. The coordinates with $u_j\le1$
contribute at most three. On $I=\{j:u_j>1\}$,
$\abs I\le p$ and $\sum_{j\in I}u_j\le p$.
For $1<u_j\le\beta_j$,
\[
 \Prob(\xi_j^{[\beta_j]}\ge u_j)=\tfrac12e^{-u_j}.
\]
This includes equality at the cap, since the truncated variable
has an atom there. Force these events for all $j\in I$.
The sum over the remaining coordinates is independent and
symmetric, so it is nonnegative with probability at least $1/2$.
The intersection therefore has probability at least
\[
 \tfrac12\prod_{j\in I}(\tfrac12e^{-u_j})
 \ge\tfrac12\exp[-(1+\log2)p].
\]
On this event $H\ge\sum_{j\in I}c_ju_j$.
It follows that $\norm H_p\ge c\sum_{j\in I}c_ju_j$.
Let $a_0=1-e^{-1}$. Conditional Jensen on the independent core gives
\[
 \|H\|_p\ge\E|H|
 \ge\E\bigl|\E[H\mid\xi^{[1]}]\bigr|
 =\E|\xi^{[1]}|=a_0.
\]
For every feasible $(u_j)$,
\[
 \sum_jc_ju_j\le3+\sum_{j\in I}c_ju_j
                  \le(3/a_0+c^{-1})\|H\|_p.
\]
Taking the supremum in \eqref{eq:column-support} proves
$Q_\Phi(p)\le C\|H\|_p$.

\emph{Passage to the full coordinate.}
Let $H^{(J)}$ retain the core and the kept bins with $j\le J$.
For each fixed $p\ge1$,
\[
 \|H-H^{(J)}\|_p
 \le\|Z_1\|_p\sum_{j>J}2^{-j}\longrightarrow0,
 \qquad \Phi^{(J)}(r)\uparrow\Phi(r).
\]
Monotone convergence gives $Q_{\Phi^{(J)}}(p)\uparrow Q_\Phi(p)$.
Passing to the limit in
$cQ_{\Phi^{(J)}}(p)\le\|H^{(J)}\|_p\le CQ_{\Phi^{(J)}}(p)$
proves the full-coordinate comparison. Together with \eqref{eq:profile-comparison} it
gives $\norm H_p\asympc q(p)$ for $p\ge2$ and a numerical
first-moment lower bound. Lemma~\ref{lem:quantile-convex}
proves \eqref{eq:source-transfer}. Since the coordinates of $H$ and
$Y$ are independent, Lemma~\ref{lem:convex-pairing}(ii) gives
\[
 C^{-1}Y\cx H\cx CY.
\]
Lemma~\ref{lem:convex-pairing}(i) then yields
$B_H(\mu)\asympc B_Y(\mu)$ for every prescribed finite index law.
\end{proof}

\begin{proof}[Proof of the fixed-law estimate \eqref{eq:main-fixed}]
For $J\ge0$, consider the finite expansion obtained by keeping the
core coordinate together with the components satisfying
$j\le J$ and $w_{ij}\ge1$, with the same caps as in the full expansion.
These caps are allowed to be infinite. Applying
Theorem~\ref{thm:capped-geometry} to this finite expansion gives
\[
\cR^{\Phi^{(J)}}(\mu)
=
\cR^{\rm cap}((L_J)\push\mu)
\asympc
B_{\Xi^{(J)}}((L_J)\push\mu)
=
B_{H^{(J)}}(\mu).
\]

We verify the last equality. By construction,
\[
H^{(J)}=L_J^{\mathsf T}\Xi^{(J)}.
\]
Hence every coupling of $X\sim\mu$ and $\Xi^{(J)}$ satisfies
\[
\E\ip{L_JX}{\Xi^{(J)}}
=
\E\ip X{H^{(J)}}
\le
B_{H^{(J)}}(\mu).
\]
Taking the supremum gives
\[
B_{\Xi^{(J)}}((L_J)\push\mu)
\le
B_{H^{(J)}}(\mu).
\]

For the converse, let $K(h,d\xi)$ be the conditional law of
$\Xi^{(J)}$ given $H^{(J)}=h$. Given any coupling $\pi(dx,dh)$ of
$X\sim\mu$ and $H^{(J)}$, define
\[
Q(dx,dh,d\xi)=\pi(dx,dh)K(h,d\xi).
\]
Then $Q$ preserves the marginal law of $\Xi^{(J)}$ and satisfies
\[
L_J^{\mathsf T}\xi=h
\]
almost surely. Therefore
\[
\E_Q\ip{L_JX}{\Xi^{(J)}}
=
\E_\pi\ip X{H^{(J)}}.
\]
Taking the supremum over couplings of $X$ and $H^{(J)}$ yields
\[
B_{H^{(J)}}(\mu)
\le
B_{\Xi^{(J)}}((L_J)\push\mu),
\]
and hence
\[
B_{\Xi^{(J)}}((L_J)\push\mu)
=
B_{H^{(J)}}(\mu).
\]

We now let $J\to\infty$. The induced costs increase pointwise to
$\Phi$, while the corresponding random vectors converge in $L^1$,
as shown in the proof of Proposition~\ref{prop:source-transfer}.
Lemma~\ref{lem:limits} therefore yields
\[
\cR^\Phi(\mu)\asympc B_H(\mu)
\]
with the same universal constants. Finally,
Propositions~\ref{prop:cost-transfer} and
\ref{prop:source-transfer} give
\[
\cR^\Psi(\mu)
\asympc
\cR^\Phi(\mu)
\asympc
B_H(\mu)
\asympc
B_Y(\mu).
\qedhere
\]
\end{proof}
\section{The two measure formulas and examples}\label{sec:one-measure}
In the preceding section, we established the first assertion of
Theorem~\ref{thm:main}. We now derive the remaining deterministic
formulations using finite-dimensional minimax. We then specialize the
result to the Gaussian and Bernoulli cases and discuss examples beyond
the $\Delta_2$ setting.
\subsection{From rate--distortion to $\Gamma^\Psi$ and $\cM^\Psi$}

Recall $F_t$ from \eqref{eq:F}. To connect $\RD^\Psi_\mu$ with
$F_t$, we first remove the constraint on the second marginal. For
$t>0$ and $\mu\in\cP(T)$, define
\[
\overline{\RD}_\mu(t)
=
\inf\left\{
  I(X;Y)+\E D_t(X,Y):
  X\sim\mu,\; Y\text{ is \(T\)-valued}
\right\}.
\]
Thus only the law of $X$ is prescribed; the law of $Y$ may be any
probability measure on $T$.
The next lemma identifies this relaxed
quantity in terms of $F_t$ and compares it with the original
rate--distortion functional.

\begin{lemma}\label{lem:relaxed-RD}
For every $t>0$ and $\mu\in\cP(T)$,
\begin{equation}\label{eq:free-RD}
\overline{\RD}_\mu(t)
=
\inf_{\nu\in\cP(T)}
\sum_{x\in T}\mu(x)F_t(x,\nu).
\end{equation}
Moreover,
\begin{equation}\label{eq:free-replica}
\overline{\RD}_\mu(t)
\le
\RD^\Psi_\mu(t),
\qquad
\RD^\Psi_\mu(2t)
\le
2\overline{\RD}_\mu(t).
\end{equation}
\end{lemma}

\begin{proof}
The first inequality in \eqref{eq:free-replica} is immediate, since
the infimum defining $\overline{\RD}_\mu(t)$ is taken over a larger
class of couplings.

We next prove \eqref{eq:free-RD}. Fix $\nu\in\cP(T)$. If $K_x$
denotes the conditional law of $Y$ given $X=x$, then
\[
\inf_{K_x}
\left\{
\KL(K_x\Vert\nu)
+
\sum_{y\in T}K_x(y)D_t(x,y)
\right\}
=
F_t(x,\nu).
\]
Indeed, when
\[
\sum_y\nu(y)e^{-D_t(x,y)}>0,
\]
the minimum is attained at
\[
K_x(y)
=
\frac{\nu(y)e^{-D_t(x,y)}}
{\sum_z\nu(z)e^{-D_t(x,z)}}.
\]
If the denominator vanishes, both sides are infinite.

Now let
\[
\Prob(X=x,Y=y)=\mu(x)K_x(y).
\]
For every $\nu\in\cP(T)$,
\[
\sum_x\mu(x)\KL(K_x\Vert\nu)
=
I(X;Y)+\KL(\Law(Y)\Vert\nu).
\]
Hence, after minimizing over $\nu$, the optimal choice is
$\nu=\Law(Y)$. Minimizing next over the kernels $K_x$ gives
\eqref{eq:free-RD}.

It remains to prove the second inequality in
\eqref{eq:free-replica}. Let $(X,Y)$ be admissible in the definition
of $\overline{\RD}_\mu(t)$. Conditional on $Y$, draw $X'$ independently from the conditional
law of $X$ given $Y$. Then $X'\sim\mu$, so $(X,X')$ is admissible
for $\RD^\Psi_\mu(2t)$. By data processing,
\[
I(X;X')\le I(X;Y).
\]
The scale triangle inequality gives
\[
D_{2t}(X,X')
\le
D_t(X,Y)+D_t(Y,X').
\]
Since $(X',Y)$ has the same law as $(X,Y)$,
\[
\E D_{2t}(X,X')
\le
2\E D_t(X,Y).
\]
Therefore
\[
\RD^\Psi_\mu(2t)
\le
I(X;Y)+2\E D_t(X,Y)
\le
2\bigl(I(X;Y)+\E D_t(X,Y)\bigr).
\]
Taking the infimum over $(X,Y)$ proves
\[
\RD^\Psi_\mu(2t)
\le
2\overline{\RD}_\mu(t).
\]
\end{proof}
\begin{theorem}\label{thm:single-measure}
For the intrinsic costs $\Psi$ defined above and every nonempty finite
set $T\subset\R^n$,
\[
\Gamma^\Psi(T)
\asympc
\cR^\Psi(T)
\asympc
\cM^\Psi(T).
\]
The constants are universal.
\end{theorem}

\begin{proof} Set
\[
r_k=2^{-k},\qquad k\in\mathbb Z.
\]

We first compare \(\Gamma^\Psi(T)\) with \(\cR^\Psi(T)\).
For \(N\ge0\), let
\[
 G_N
 =
 \inf_{(\nu_k)_{|k|\le N}}
 \max_{x\in T}
 \sum_{|k|\le N}r_kF_{r_k}(x,\nu_k).
\]
Because \(F_t(x,\nu)\) may be infinite on the boundary of the
simplex, we first regularize the measures. Let \(u_T\) be the uniform
law on \(T\), and for \(0<\varepsilon<1\) set
\[
 \cP_\varepsilon(T)
 =
 \{(1-\varepsilon)\nu+\varepsilon u_T:\nu\in\cP(T)\}.
\]
On \(\cP_\varepsilon(T)\), every \(F_t(x,\nu)\) is finite,
continuous, and convex in \(\nu\). For any function $g:T\to\R$,
\[
 \max_{x\in T}g(x)
 =
 \sup_{\mu\in\cP(T)}\sum_{x\in T}\mu(x)g(x).
\]
Thus the finite-scale objective is linear in $\mu$ and convex in each
$\nu_k$. Sion's minimax theorem therefore gives
\begin{align*}
 &\inf_{(\nu_k)_{|k|\le N}\in\cP_\varepsilon(T)^{2N+1}}
 \max_{x\in T}
 \sum_{|k|\le N}r_kF_{r_k}(x,\nu_k)\\
 &\qquad=
 \sup_{\mu\in\cP(T)}
 \sum_{|k|\le N}r_k
 \inf_{\nu\in\cP_\varepsilon(T)}
 \sum_{x\in T}\mu(x)F_{r_k}(x,\nu).
\end{align*}

If \(\nu^\varepsilon=(1-\varepsilon)\nu+\varepsilon u_T\), then
\[
 F_t(x,\nu^\varepsilon)
 \le
 F_t(x,\nu)-\log(1-\varepsilon).
\]
Thus restricting the minimization to \(\cP_\varepsilon(T)\) changes
both sides by at most
\[
 -\log(1-\varepsilon)\sum_{|k|\le N}r_k.
\]
Letting \(\varepsilon\downarrow0\) and using
Lemma~\ref{lem:relaxed-RD}, we obtain the exact finite-scale identity
\begin{equation}\label{eq:finite-Gamma-free-RD}
 G_N
 =
 \sup_{\mu\in\cP(T)}
 \sum_{|k|\le N}r_k\overline{\RD}_\mu(r_k).
\end{equation}

We now pass to all dyadic scales. The sequence \(G_N\) is
nondecreasing and \(G_N\le\Gamma^\Psi(T)\). Let
\(G=\lim_NG_N\), and for each \(N\) choose measures
\((\nu_k^{(N)})_{|k|\le N}\) within \(1/N\) of \(G_N\). By
compactness of \(\cP(T)\) and a diagonal subsequence, for every fixed
\(k\) we may assume
\[
 \nu_k^{(N)}\longrightarrow\nu_k.
\]
Since \(F_t(x,\cdot)\) is lower semicontinuous, for every fixed
\(M\) and \(x\in T\),
\[
 \sum_{|k|\le M}r_kF_{r_k}(x,\nu_k)
 \le
 \liminf_{N\to\infty}
 \sum_{|k|\le M}r_kF_{r_k}(x,\nu_k^{(N)})
 \le G.
\]
Letting \(M\to\infty\) shows that the limiting sequence
\((\nu_k)_{k\in\mathbb Z}\) is admissible for \(\Gamma^\Psi(T)\)
with value at most \(G\). Hence
\[
 \Gamma^\Psi(T)=\lim_{N\to\infty}G_N.
\]
On the right-hand side of \eqref{eq:finite-Gamma-free-RD}, all terms
are nonnegative, so
\[
 \lim_{N\to\infty}
 \sup_{\mu\in\cP(T)}
 \sum_{|k|\le N}r_k\overline{\RD}_\mu(r_k)
 =
 \sup_{\mu\in\cP(T)}
 \sum_{k\in\mathbb Z}r_k\overline{\RD}_\mu(r_k).
\]
Therefore
\begin{equation}\label{eq:Gamma-free-RD}
 \Gamma^\Psi(T)
 =
 \sup_{\mu\in\cP(T)}
 \sum_{k\in\mathbb Z}r_k\overline{\RD}_\mu(r_k).
\end{equation}

Lemma~\ref{lem:relaxed-RD} and \eqref{eq:dyadic-area} imply
\[
 \frac14\,\cR^\Psi(\mu)
 \le
 \sum_{k\in\mathbb Z}r_k\overline{\RD}_\mu(r_k)
 \le
 2\,\cR^\Psi(\mu).
\]
Taking the supremum over \(\mu\) in \eqref{eq:Gamma-free-RD} gives
\begin{equation}\label{eq:Gamma-RD}
 \frac14\,\cR^\Psi(T)
 \le
 \Gamma^\Psi(T)
 \le
 2\,\cR^\Psi(T).
\end{equation}

We next compare $\Gamma^\Psi(T)$ with $\cM^\Psi(T)$. If the same
measure $\nu\in\cP(T)$ is used at every dyadic scale, then monotonicity
of $t\mapsto F_t(x,\nu)$ gives
\[
\sum_{k\in\mathbb Z}r_kF_{r_k}(x,\nu)
\le
2\int_0^\infty F_t(x,\nu)\,dt.
\]
Taking the maximum over $x$ and then the infimum over $\nu$ gives
\begin{equation}\label{eq:Gamma-M-upper}
\Gamma^\Psi(T)\le2\cM^\Psi(T).
\end{equation}

For the reverse inequality, let
\[
D=\diam_2(T).
\]
There is nothing to prove when $D=0$. Choose $k_0\in\mathbb Z$ such
that
\[
D\le R_0:=r_{k_0}<2D.
\]
Fix $\varepsilon>0$ and choose $(\nu_k)_{k\in\mathbb Z}$ such that
\[
\max_{x\in T}
\sum_{k\in\mathbb Z}r_kF_{r_k}(x,\nu_k)
\le
\Gamma^\Psi(T)+\varepsilon.
\]
For $j\ge0$, set
\[
t_j=R_0\,2^{-j}=r_{k_0+j},
\qquad
\eta_j=\nu_{k_0+j},
\]
and define
\[
\nu=\sum_{j\ge0}2^{-j-1}\eta_j.
\]
Since $\nu\ge2^{-j-1}\eta_j$,
\[
F_{t_j}(x,\nu)
\le
F_{t_j}(x,\eta_j)+(j+1)\log2.
\]
Using the monotonicity of $F_t(x,\nu)$ on each interval
$[t_{j+1},t_j]$, we obtain
\[
\begin{aligned}
\int_0^{R_0}F_t(x,\nu)\,dt
&\le
\sum_{j\ge1}t_jF_{t_j}(x,\nu)\\
&\le
\Gamma^\Psi(T)+\varepsilon
+
\log2\sum_{j\ge1}t_j(j+1)\\
&\le
\Gamma^\Psi(T)+\varepsilon+CR_0.
\end{aligned}
\]

For $t\ge R_0\ge D$, all coordinate errors lie in the quadratic range
of $\Psi$. Hence
\[
D_t(x,y)
=
\frac{\|x-y\|_2^2}{4t^2}
\le
\frac{D^2}{4t^2}.
\]
Moreover, Jensen's inequality gives
\[
F_t(x,\nu)
\le
\sum_{y\in T}\nu(y)D_t(x,y).
\]
Therefore
\[
\int_{R_0}^\infty F_t(x,\nu)\,dt
\le
\frac{D^2}{4R_0}
\le
CD.
\]
Since $R_0<2D$, letting $\varepsilon\downarrow0$ gives
\begin{equation}\label{eq:measure-error}
\cM^\Psi(T)
\le
\Gamma^\Psi(T)+CD.
\end{equation}

It remains to absorb the diameter term. Choose $x,y\in T$ with
\[
\|x-y\|_2=D
\]
and put
\[
\mu_0=\frac12(\delta_x+\delta_y).
\]
Put \(\delta_t=D_t(x,y)\). Every self-coupling of \(\mu_0\) is
determined by its total off-diagonal mass \(q\in[0,1]\). Its mutual
information is \(\log2-h(q)\), where \(h\) is the binary entropy, and
its expected distortion is \(q\delta_t\). Thus
\[
 \RD^\Psi_{\mu_0}(t)
 =
 \min_{0\le q\le1}\{\log2-h(q)+q\delta_t\}
 =
 \log2-\log(1+e^{-\delta_t})
 \ge
 c\min\{1,\delta_t\}.
\]
Hence Lemma~\ref{lem:area-identity}, with $p=1$, gives
\[
\begin{aligned}
\cR^\Psi(T)
&\ge
\cR^\Psi(\mu_0)\\
&\ge
c\int_0^\infty\min\{1,D_t(x,y)\}\,dt\\
&=
c\|x-y\|_2
=
cD.
\end{aligned}
\]
Combining this with \eqref{eq:Gamma-RD}, we obtain
\[
D\le C\Gamma^\Psi(T).
\]
Substituting this estimate into \eqref{eq:measure-error} yields
\[
\cM^\Psi(T)\le C\Gamma^\Psi(T).
\]
Together with \eqref{eq:Gamma-M-upper} and \eqref{eq:Gamma-RD}, this
proves the theorem.
\end{proof}
\begin{proof}[Completion of Theorem~\ref{thm:main}]
The first estimate in  Theorem~\ref{thm:main}
\[
B_Y(\mu)\asympc \cR^\Psi(\mu)
\]
was proved in Section~\ref{sec:general-tails}. Taking the supremum over
$\mu\in\cP(T)$ and using \eqref{eq:set-pairing}, we obtain
\[
S_Y(T)
=
\sup_{\mu\in\cP(T)}B_Y(\mu)
\asympc
\sup_{\mu\in\cP(T)}\cR^\Psi(\mu)
=
\cR^\Psi(T).
\]
Theorem~\ref{thm:single-measure} gives
\[
\cR^\Psi(T)
\asympc
\Gamma^\Psi(T)
\asympc
\cM^\Psi(T).
\]
Combining the two estimates completes the proof of Theorem~\ref{thm:main}.
\end{proof}

\subsection{Independent coordinates without symmetry}

For a nonsymmetric source we apply the theorem to its difference
with an independent copy. We give the coordinate rescaling explicitly.

\begin{corollary}\label{cor:nonsymmetric}
Let $A_1,\ldots,A_n$ be independent centered nondegenerate variables
with log-concave densities, and let $A'$ be an independent copy.
Put $Y^0=A-A'$. Let $q_i^0$ be the inverse tail of $|Y_i^0|$,
as in the definition of $q_i$ above, and set
\[
 \lambda_i=q_i^0(1),\qquad
 \widetilde Y_i=Y_i^0/\lambda_i,
 \qquad \Lambda x=(\lambda_i x_i)_{i=1}^n.
\]
Construct $\Psi$ from the normalized source $\widetilde Y$.
For every finite nonempty $T\subset\R^n$,
\begin{equation}\label{eq:nonsymmetric}
 \E\max_{x,y\in T}\ip{x-y}A
       \asympc\Gamma^\Psi(\Lambda T)\asympc\cM^\Psi(\Lambda T).
\end{equation}
Moreover, for every finite index law $\mu$,
\[
 B_{Y^0}(\mu)\asympc\cR^\Psi(\Lambda\push\mu).
\]
\end{corollary}
\begin{proof}
If $f_i$ is the density of $A_i$, then
\[
 g_i(s)=\int_\R f_i(s+t)f_i(t)\,dt
\]
is the density of $Y_i^0$. Its integrand is log-concave in $(s,t)$,
so Theorem~\ref{thm:prekopa-input} makes $g_i$ log-concave.
It is even. The same theorem applied to a half-line gives
\[
 \Prob(|Y_i^0|\ge u)=2\int_0^\infty g_i(u+t)\,dt
 \quad(u\ge0),
\]
which is a log-concave function of $u$. The coordinates of $Y^0$
are independent and symmetric.

Set $f_T(z)=\max_{x,y\in T}\ip{x-y}z$. This function is convex,
even, and subadditive. Since $\E A'=0$,
\begin{align*}
 \E f_T(A)
 &=\E f_T\bigl(\E[A-A'\mid A]\bigr)
 \le\E f_T(A-A')\\
 &\le\E f_T(A)+\E f_T(-A')=2\E f_T(A).
\end{align*}
The rescaling satisfies
$\ip x{Y^0}=\ip{\Lambda x}{\widetilde Y}$, so symmetry gives
\[
 \E f_T(Y^0)=2S_{\widetilde Y}(\Lambda T).
\]
Apply Theorem~\ref{thm:main} to $\widetilde Y$ and $\Lambda T$.
\end{proof}
\subsection{Classical cases and examples beyond the $\Delta_2$ condition}

We first verify that Theorem~\ref{thm:main} recovers the Gaussian and
Bernoulli cases.

\begin{proof}[Proof of Corollary~\ref{cor:gaussian-bernoulli}]
For a symmetric sign, the endpoint is $b=1$, and by our convention
\[
V(u)=u^2,\qquad |u|\le1.
\]
Hence
\[
\Psi(r)
=
\int_0^1 2u\,\ind_{\{2u<r\}}\,du
=
\min\left(1,\frac{r^2}{4}\right)
=
\psi_{\rm B}(r/2).
\]
Thus Theorem~\ref{thm:main} gives the Bernoulli assertions. The fixed
change of scale in the argument of $\psi_{\rm B}$ changes the
corresponding functionals only by numerical factors.

We next consider a standard Gaussian variable $g$. Choose $a>0$ such
that
\[
\Prob(|g|\ge a)=e^{-1},
\]
and normalize by $Y=g/a$. Its tail exponent satisfies
\[
N(1)=1,
\qquad
N(u)=-\log\Prob(|g|\ge au).
\]
For $u\ge1$,
\[
N'(u)
=
\frac{a\varphi(au)}
{\displaystyle\int_{au}^{\infty}\varphi(v)\,dv},
\]
where $\varphi$ is the standard Gaussian density. The standard Gaussian
tail estimates
\[
c\,\frac{\varphi(x)}{1+x}
\le
\int_x^\infty\varphi(v)\,dv
\le
\frac{\varphi(x)}{x},
\qquad x>0,
\]
show that
\[
N'(u)\asympc u,
\qquad u\ge1.
\]
Together with the quadratic definition of $V$ on $[0,1]$, this gives
numerical constants $c,C>0$ such that
\[
cu\le V'_+(u)\le Cu,
\qquad u>0.
\]
It follows directly from the definition of $\Psi$ that
\[
c'r^2\le\Psi(r)\le C'r^2,
\qquad r\ge0.
\]
Indeed,
\[
\int_0^{r/C}cu\,du
\le
\Psi(r)
\le
\int_0^{r/c}Cu\,du.
\]
Thus the intrinsic Gaussian cost is equivalent to
\[
\psi_{\rm G}(r)=r^2.
\]
Theorem~\ref{thm:main}, together with the fixed normalization above,
therefore gives the Gaussian assertions.

Finally, Talagrand's majorizing-measure theorem and the Bernoulli
theorem of Bednorz and Lata\l a give
\[
S_G(T)\asympc\gamma_2(T),
\qquad
S_\eps(T)\asympc\mathfrak b(T),
\]
respectively \cite{Talagrand,BL}. Hence our functionals recover the
classical Gaussian and Bernoulli descriptions.
\end{proof}

We next give an example for which the $\Delta_2$ condition fails. It
also shows why the intrinsic cost cannot in general be replaced by the
usual chain of increment moments.

Consider the tail exponent
\begin{equation}\label{eq:superexp-example}
N(u)=2^u-1,
\qquad
q(p)=\frac{\log(1+p)}{\log2}.
\end{equation}
Then
\[
\frac{N(2u)}{N(u)}
=
2^u+1
\longrightarrow\infty,
\]
so no $\Delta_2$ bound holds. Since
\[
V(u)
=
\begin{cases}
u^2,&0\le u\le1,\\
2^{u+1}-3,&u>1,
\end{cases}
\]
we have
\[
V'_+(u)
=
\begin{cases}
2u,&0<u<1,\\
2^{u+1}\log2,&u>1.
\end{cases}
\]
Consequently,
\[
\Psi(r)
=
\begin{cases}
r^2/4,&0\le r\le2,\\
1,&2<r\le4\log2,\\
r/\log2-3,&r>4\log2.
\end{cases}
\]

Let $Y_1,\ldots,Y_m$ be independent copies of the corresponding
symmetric variable and set
\[
T_m=\{0,e_1,\ldots,e_m\}.
\]
Then
\begin{equation}\label{eq:basis-example}
S_Y(T_m)
\asympc
\log\log(m+e)
\asympc
\cR^\Psi(T_m)
\asympc
\cM^\Psi(T_m).
\end{equation}
Only the first comparison needs to be checked. By
\eqref{eq:quantile-representation}, we may write
\[
|Y_i|\stackrel{\rm d}=q(E_i),
\]
where $E_1,\ldots,E_m$ are independent mean-one exponential variables.
If
\[
M_m=\max_{i\le m}E_i,
\]
then
\[
S_Y(T_m)
\le
\E q(M_m)
\le
q(\E M_m)
\le
q(\log m+1)
\lesssim
\log\log(m+e),
\]
where we used the concavity of $q$.

For the reverse inequality, put
\[
u_m=q(\log(m+1)).
\]
For each $i$,
\[
\Prob(Y_i\ge u_m)
=
\frac1{2(m+1)}.
\]
Hence, with probability bounded below by a numerical constant, at least
one of $Y_1,\ldots,Y_m$ is at least $u_m$. Therefore
\[
S_Y(T_m)
\ge
c\,u_m
\asympc
\log\log(m+e).
\]
The remaining comparisons in \eqref{eq:basis-example} follow from
Theorem~\ref{thm:main}.

This example also separates our functional from the usual full-moment
chain. Indeed,
\[
\norm{Y_i}_p\asympc\log(p+1),
\qquad p\ge2,
\]
and, for $i\ne j$, conditional Jensen gives
\[
\norm{Y_i-Y_j}_p
\ge
\norm{Y_i}_p
\gtrsim
\log(p+1).
\]
Let $K$ be the largest integer such that
\[
2^{2^K}<m+1.
\]
For any admissible sequence of partitions, some cell at level $K$
contains at least two points of $T_m$. Choose a point in such a cell.
By nestedness, its cells at every level $k\le K$ are also
nonsingletons. Their $L^{2^k}$-diameters are therefore at least
$c k$. Hence
\[
\gamma_Y(T_m)
\ge
c\sum_{k=1}^K k
\ge
c'K^2
\asympc
\bigl(\log\log(m+e)\bigr)^2.
\]
Thus
\[
\frac{\gamma_Y(T_m)}{S_Y(T_m)}
\longrightarrow\infty.
\]
In particular, once the $\Delta_2$ condition is removed, the usual
moment chain need not have the correct order, whereas
$\Gamma^\Psi(T_m)\asympc S_Y(T_m)$ by Theorem~\ref{thm:main}.

The failure of $\Delta_2$ may also occur through arbitrarily large
changes of slope. To illustrate this, let
\[
b_j=4^j,\qquad j\ge0,
\]
and
\[
L_j=\exp(16^j),\qquad j\ge1.
\]
Set $N(u)=u$ on $[0,1]$, and for $j\ge1$ let $N$ have constant slope
$L_j$ on $[b_{j-1},b_j]$. Then $N$ is convex and normalized by
$N(1)=1$. Moreover,
\[
2b_j<b_{j+1}
\]
and
\[
N(b_j)\le L_jb_j.
\]
Therefore
\[
\frac{N(2b_j)}{N(b_j)}
=
1+\frac{L_{j+1}b_j}{N(b_j)}
\ge
1+\frac{L_{j+1}}{L_j}
\longrightarrow\infty.
\]
Thus this tail has infinitely many arbitrarily large violations of the
$\Delta_2$ condition. Nevertheless, the expansion
\eqref{eq:source-expansion} still has coefficient sum at most three in
each coordinate. In particular, the constants in
Theorem~\ref{thm:main} remain universal even for such tails.
\appendix
\section{An estimate for polyhedral partitions}

The contraction theorem also has a geometric consequence for convex
polyhedral partitions. This appendix records the corresponding surface
representation and contraction estimate. These results are not used
in the proof of Theorem~\ref{thm:main}.

Let
\[
 G(z)=\max_{j\le m}(\ip{y_j}z-b_j),
\]
where identical affine functions have been merged and only cells of
positive probability are retained. Write
\[
 C_j=\{z:\ip{y_j}z-b_j=G(z)\}.
\]
Cell boundaries have probability zero, and off these boundaries
\(\nabla G(z)=y_j\) on \(C_j\).

We first record the surface representation used below. Set
\[
 p_d(z)=2^{-d}e^{-\norm z_1},
 \qquad
 \tau(s)=1+\abs s,
\]
and, for \(L>0\),
\[
 \tau_L(s)=
 \begin{cases}
 e^{\abs s-L},&\abs s\le L,\\
 1+\abs s-L,&\abs s>L.
 \end{cases}
\]
If \(C_j,C_k\) share a facet of dimension \(d-1\), denote its
relative interior by \(\Gamma_{jk}\) and put \(d_{jk}=y_j-y_k\).
Pairs without such a facet contribute zero. The measure
\(\mathcal H^{d-1}\) is the ordinary surface measure on these facets;
for \(d=1\) it counts boundary points.

\begin{lemma}[Surface representation]
With
\[
 \mathcal E(G)=\E\ip{\nabla G(Z)}Z,
 \qquad
 \mathcal E_{K\beta}(G)=
 \E\ip{\nabla G(Z)}{Z-\clip_{K\beta}Z},
\]
one has
\begin{align}
 \mathcal E(G)
 &=
 \sum_{j<k}\int_{\Gamma_{jk}}
 \frac{\sum_i\tau(z_i)d_{jk,i}^2}{\norm{d_{jk}}_2}
 p_d(z)\,d\mathcal H^{d-1}(z),
 \label{eq:surface-full}\\
 \mathcal E_{K\beta}(G)
 &=
 \sum_{j<k}\int_{\Gamma_{jk}}
 \frac{\sum_i\tau_{K\beta_i}(z_i)d_{jk,i}^2}{\norm{d_{jk}}_2}
 p_d(z)\,d\mathcal H^{d-1}(z).
 \label{eq:surface-tail}
\end{align}
In particular, both quantities are nonnegative.
\end{lemma}

\begin{proof}
In one dimension,
\[
 (\tau p_1)'=-sp_1,
 \qquad
 (\tau_Lp_1)'=-(s-\clip_Ls)p_1.
\]
The weighted densities are continuous at zero and at \(\pm L\), so
these derivatives have no additional atoms. Across a common facet the
gradient changes from \(y_k\) to \(y_j\). Its unit normal in that
direction is \(d_{jk}/\norm{d_{jk}}_2\), and therefore, in the sense
of distributions,
\[
 D^2G\big|_{\Gamma_{jk}}
 =
 \frac{d_{jk}\otimes d_{jk}}{\norm{d_{jk}}_2}
 \mathcal H^{d-1}\big|_{\Gamma_{jk}}.
\]
Here
\((d_{jk}\otimes d_{jk})_{i\ell}=d_{jk,i}d_{jk,\ell}\).
Multiplying the diagonal entries by the one-dimensional weights gives
the facet terms in \eqref{eq:surface-full}--\eqref{eq:surface-tail}.

For completeness, we justify the integration at infinity. Write
\(a_i(z_i)\) for either \(\tau(z_i)\) or
\(\tau_{K\beta_i}(z_i)\), and \(b_i(z_i)\) for the corresponding
\(z_i\) or \(z_i-\clip_{K\beta_i}(z_i)\). Choose smooth cutoffs
\(\eta_R\) equal to one on \(\|z\|_2\le R\), zero on
\(\|z\|_2\ge2R\), with \(|\nabla\eta_R|\le C/R\). Since
\(\partial_i(a_ip_d)=-b_ip_d\), integration by parts gives
\[
 \int\eta_R b_i\partial_iG\,p_d\,dz
 =
 \int\eta_Ra_ip_d\,d(D_{ii}G)
 +
 \int(\partial_i\eta_R)(\partial_iG)a_ip_d\,dz.
\]
If \(M_G=\max_j\|y_j\|_2\), the last integral has absolute value at
most
\[
 \frac{CM_G}{R}
 \int_{\|z\|_2\ge R}(1+|z_i|)p_d(z)\,dz
 \longrightarrow0.
\]
The finite facet sum is integrable against the exponentially decaying
weights. Letting \(R\to\infty\) proves the two identities, whose
integrands are nonnegative.
\end{proof}

\begin{theorem}\label{thm:surface}
Suppose
\[
 \abs{\E[Z_i\mid Z\in C_j]}\le\beta_i,
 \qquad 1\le\beta_i<\infty,
\]
for every positive-probability cell \(C_j\). Then, for sufficiently
large numerical \(K\),
\begin{equation}\label{eq:surface-contraction}
 0\le\mathcal E_{K\beta}(G)
 \le
 Ce^{-cK}\mathcal E(G).
\end{equation}
\end{theorem}

\begin{proof}
Let \(J\) be the cell label and set
\[
 U=\E[Z\mid J],
 \qquad
 R=\E[Z-\clip_{K\beta}Z\mid J].
\]
By the hypothesis on the cells,
\[
 |U_i|\le\beta_i\qquad\text{for every }i.
\]
The conditional coordinate density on a convex cell has the form
\(e^{-\abs s}h(s)\) with \(h\) log-concave, by Pr\'ekopa's theorem.
Lemma~\ref{lem:conditional-tail} therefore gives
\[
 |R_i|\le\epsilon_K(1+|U_i|).
\]
Conditional Jensen gives \(U\cx Z\), and the residual is dominated in
convex order by the independent sparse overshoot law whose selector
probabilities are at most \(e^{-K}\). Theorem~\ref{thm:amplitude}
then yields
\[
 R\cx\theta_KZ,
 \qquad
 \theta_K\le Ce^{-cK}.
\]

Let \(\nu=\Law(y_J)\). The actual cell assignment is optimal among
couplings with these label masses. Indeed, for every competing
coupling,
\[
 \ip{y_J}Z-b_J\le G(Z),
\]
while the expected intercept depends only on the law of \(J\), and
the cell coupling attains equality. Hence
\[
 B_Z(\nu)=\E\ip{y_J}Z=\mathcal E(G).
\]
Consequently,
\[
 \mathcal E_{K\beta}(G)
 =
 \E\ip{y_J}R
 \le
 B_R(\nu)
 \le
 \theta_KB_Z(\nu)
 \le
 Ce^{-cK}\mathcal E(G).
\]
Here the estimate $B_R(\nu)\le\theta_KB_Z(\nu)$ follows from
$R\cx\theta_KZ$ and Lemma~\ref{lem:convex-pairing}(i).
The lower bound follows from the nonnegative surface representation
\eqref{eq:surface-tail}.
\end{proof}

Equation~\eqref{eq:surface-contraction} compares the complete facet
sums. It does not assert a separate contraction estimate for each
individual facet.

\section*{Acknowledgments}
This work was supported by the National Key R\&D Program of China (No.2024YFA1013501), the National Natural Science Foundation of China (No. 12571162 and 12371148), and Shandong Provincial Natural Science Foundation (No. ZR2024MA082).

\end{document}